\documentclass[a4paper,11pt]{article}
\usepackage{amsmath,amssymb,amsthm}
\usepackage{mathrsfs}
\usepackage{cases}
\usepackage{a4wide}
\usepackage[hang,small,bf,justification=centering]{caption}
\usepackage[subrefformat=parens]{subcaption}
\usepackage{here}
\usepackage{comment}
\usepackage{ulem}
\usepackage{color}

\usepackage{graphicx} % Required for inserting images
\numberwithin{equation}{section}

\newtheorem{prop}{Proposition}[section]
\newtheorem{thm}{Theorem}[section]

\newtheorem{lem}{Lemma}[section]
\newtheorem{rem}{Remark}[section]

\newtheorem{asm}{Assumption}[section]

\newcommand{\Section}[1]{Section~\ref{#1}}
\newcommand{\Figure}[1]{Figure~\ref{#1}}

\newcommand{\Proposition}[1]{Proposition~\ref{#1}}

\newcommand{\Lemma}[1]{Lemma~\ref{#1}}

\newcommand{\Assumption}[1]{Assumption~\ref{#1}}

\def\bv#1{\mbox{\boldmath{$#1$}}}    % bold symbol to distingush from R^2

\newcommand{\dL}[1]{{\rm d}\mathscr{L}^{#1}}
\newcommand{\dH}[1]{{\rm d}\mathscr{H}^{#1}}
\newcommand{\cpindexp}[1]{b_{#1}^+}
\newcommand{\cpindexm}[1]{b_{#1}^-}
\newcommand{\tpindex}[2]{s^{#1}_{#2}}

\newcommand{\refsurface}[1]{\Upsilon_{#1}}
\newcommand{\simplex}[2]{\sigma^{#1}_{#2}}
\newcommand{\naturalindex}[1]{\mathbb{N}_{\leq #1}}
\newcommand{\mL}{\mathscr{L}}
\newcommand{\mT}{\mathcal{T}}
\newcommand{\mR}{\mathcal{R}}
\newcommand{\mV}{\mathcal{V}}
\newcommand{\subdomains}[1]{\mR_{#1}[\Gamma(t)]}
\newcommand{\subdomainsDiscrete}[2]{\mR_{#1}[\Gamma^{#2}]}
\newcommand{\mH}{\mathscr{H}}

\newcommand{\id}{{\rm id}}
\newcommand{\dd}[1]{\frac{\rm d}{\rm d#1}}
\newcommand{\ddt}{\dd{t}}

\newcommand{\sgrad}[0]{\nabla_{\Gamma}}
\newcommand{\sLaplacian}[0]{\Delta_{\Gamma}}

\newcommand{\inprod}[2]{\left<#1,\,#2\right>}
\newcommand{\inprodMassLamped}[4]{\left<#1,\,#2\right>^{#3}_{#4}}

\newcommand{\jump}[3]{\left[#3\right]_{#1}^{#2}}
\newcommand{\tension}[1]{\varsigma_{#1}}

\newcommand{\closure}[1]{\overline{#1}}

\newcommand{\smallChemical}[0]{w}
\newcommand{\largeChemical}[0]{W}

\newcommand{\triangulation}[1]{\mathscr{T}^{#1}}
\newcommand{\vVertex}[2]{\Vec{X}^{#1}_{#2}}
\newcommand{\vVertexSynonym}[2]{\Vec{q}^{#1}_{#2}}
\newcommand{\femspaceChemical}[1]{S^{#1}}
\newcommand{\femspaceVertex}[2]{V^{#1}_{#2}}
\newcommand{\femspaceVertexVector}[2]{\underline{V}^{#1}_{#2}}
\newcommand{\femNormalVertex}[2]{\Vec{\omega}^{#1}_{#2}}
\newcommand{\femspaceGamma}[1]{V^h(\Gamma^{#1})}
\newcommand{\femspaceGammaVector}[1]{\underline{V}^h(\Gamma^{#1})}
\newcommand{\femspaceBasisBulk}[2]{\Psi^{#1}_{#2}}
\newcommand{\femspaceBasisCurve}[2]{\Phi^{#1}_{#2}}
\newcommand{\metrics}[1]{\uuline{#1}}
\newcommand{\interpolation}[2]{\pi^{#1}_{#2}}

\newcommand{\normal}[2]{\Vec{\nu}^{#1}_{#2}}
\newcommand{\Normal}[1]{\Vec{A}\{#1\}}
\newcommand{\conormal}[2]{\Vec{\mu}^{#1}_{#2}}
\newcommand{\polygonCurve}[2]{\Gamma^{#1}_{#2}}
\newcommand{\meanCurvature}[1]{\varkappa_{#1}}
\newcommand{\curvature}[1]{\varkappa_{#1}}
\newcommand{\discreteCurvature}[2]{\kappa^{#1}_{#2}}
\newcommand{\numDomVertex}[1]{K^{#1}_{\Omega}}
\newcommand{\numSurfVertex}[2]{K^{#1}_{#2}}
\newcommand{\numSurfVertexFull}[1]{K^{#1}_\Gamma}
\newcommand{\numSimplex}[2]{J^{#1}_{#2}}
\newcommand{\numTjVertex}[1]{Z_{#1}}
\newcommand{\zerovec}[0]{\Vec{0}}
\newcommand{\phaseContent}[0]{\bv \beta}

\newcommand{\discreteVol}[1]{v^{#1}}
\newcommand{\polyhedral}[2]{\vec{\mathfrak{X}}^{#1}_{#2}}
\newcommand{\DtN}[3]{\Lambda^{#1}_{#2}{#3}}
\newcommand{\orderedSeq}[2]{\vec{\varrho}^{#1}_{#2}}
\newcommand{\volume}[1]{\operatorname{vol}\left(#1\right)}

\newcommand{\bR}{{\mathbb R}}
\newcommand{\errorXx}{\|\Gamma^h - \Gamma\|_{L^\infty}}
\newcommand{\errorUu}{\|\largeChemical^h - \smallChemical\|_{L^\infty}}
\newcommand{\revised}[1]{#1}
\title{
A parametric finite element method for a degenerate multi-phase Stefan problem with triple junctions
}
\author{Tokuhiro Eto\thanks{Universit\'{e} Claude Bernard Lyon 1, CNRS, Centrale Lyon, INSA Lyon, Universit\'{e} Jean Monnet, ICJ UMR5208, 69622 Villeurbanne, France. E-mail: eto@math.univ-lyon1.fr}\and Harald Garcke\thanks{Fakult\"{a}t f\"{u}r Mathematik, Universit\"{a}t Regensburg, 93040 Regensburg, Germany. E-mail: harald.garcke@ur.de}\and Robert N\"{u}rnberg\thanks{Dipartimento di Mathematica, Universit\`{a} di Trento, 38123 Trento, Italy. E-mail: robert.nurnberg@unitn.it}}
\date{}

\begin{document}

\maketitle

\begin{abstract}
    In this study, we propose a parametric finite element method for a degenerate multi-phase Stefan problem with triple junctions.
    This model describes the energy-driven motion of a surface cluster whose distributional solution was studied by Garcke and Sturzenhecker.
    We approximate the weak formulation of this sharp interface model by an unfitted finite element method that uses parametric elements
    for the representation of the moving interfaces.
    We establish existence and uniqueness of the discrete solution and prove unconditional stability of the proposed scheme.
    Moreover, a modification of the original scheme leads to a structure-preserving variant, in that it conserves the discrete analogue
    of a quantity that is preserved by the classical solution. Some numerical
    results demonstrate the applicability of our introduced schemes.
\end{abstract}

\section{Introduction}\label{sec:intro}
In this study, we consider the evolution of a surface cluster $\Gamma(t)$ in $\mathbb{R}^d$, $d=2,3$, which is governed by the following system of equations for $(w(\cdot,t),\Gamma(t))$:
\begin{equation}\label{eq:sharp_p}
  \begin{cases}
    \Delta\smallChemical = 0&\qquad\mbox{in}\quad\Omega\backslash\Gamma(t),\quad t \in (0,T),\\
    \jump{\cpindexm{i}}{\cpindexp{i}}{\beta}\smallChemical = -\tension{i}\meanCurvature{i}&\qquad\mbox{on}\quad\Gamma_{i}(t),\quad t \in (0,T),\quad i\in\naturalindex{I_S},\\
    \jump{\cpindexm{i}}{\cpindexp{i}}{\beta}V_{i} = \jump{\cpindexm{i}}{\cpindexp{i}}{\nabla\smallChemical}\cdot\normal{}{i}&\qquad\mbox{on}\quad\Gamma_{i}(t),\quad t \in (0,T),\quad i\in\naturalindex{I_S},\\
    \sum_{j=1}^{3}\tension{\tpindex{k}{j}}\conormal{}{\tpindex{k}{j}} = \zerovec&\qquad\mbox{on}\quad\mT_k(t),\quad t \in (0,T),\quad k\in\naturalindex{I_T},\\
    \nabla\smallChemical\cdot\normal{}{\Omega} = 0&\qquad\mbox{on}\quad\partial\Omega,\quad t \in (0,T),\\
    \Gamma(0) = \Gamma_0,
  \end{cases}
\end{equation}
where $T > 0$ is a time horizon; $\naturalindex{K} = \{1,\ldots,K\}$ denotes the set of all natural numbers which are not larger than $K$ for each $K\in\mathbb{N}$,
and $\Omega\subset\mathbb{R}^d$ is a smooth bounded domain.
The problem \eqref{eq:sharp_p} represents a degenerate multi-phase Stefan
problem first studied by Garcke and Sturzenhecker in \cite{GarckeSturzenhecker1998}. In particular, the domain $\Omega$ is divided into several phases, which are separated by a surface cluster consisting of $I_S$ surfaces that meet at
$I_T$ triple junctions.

\begin{figure}[H]
    \centering
    \includegraphics[keepaspectratio, scale=0.25]{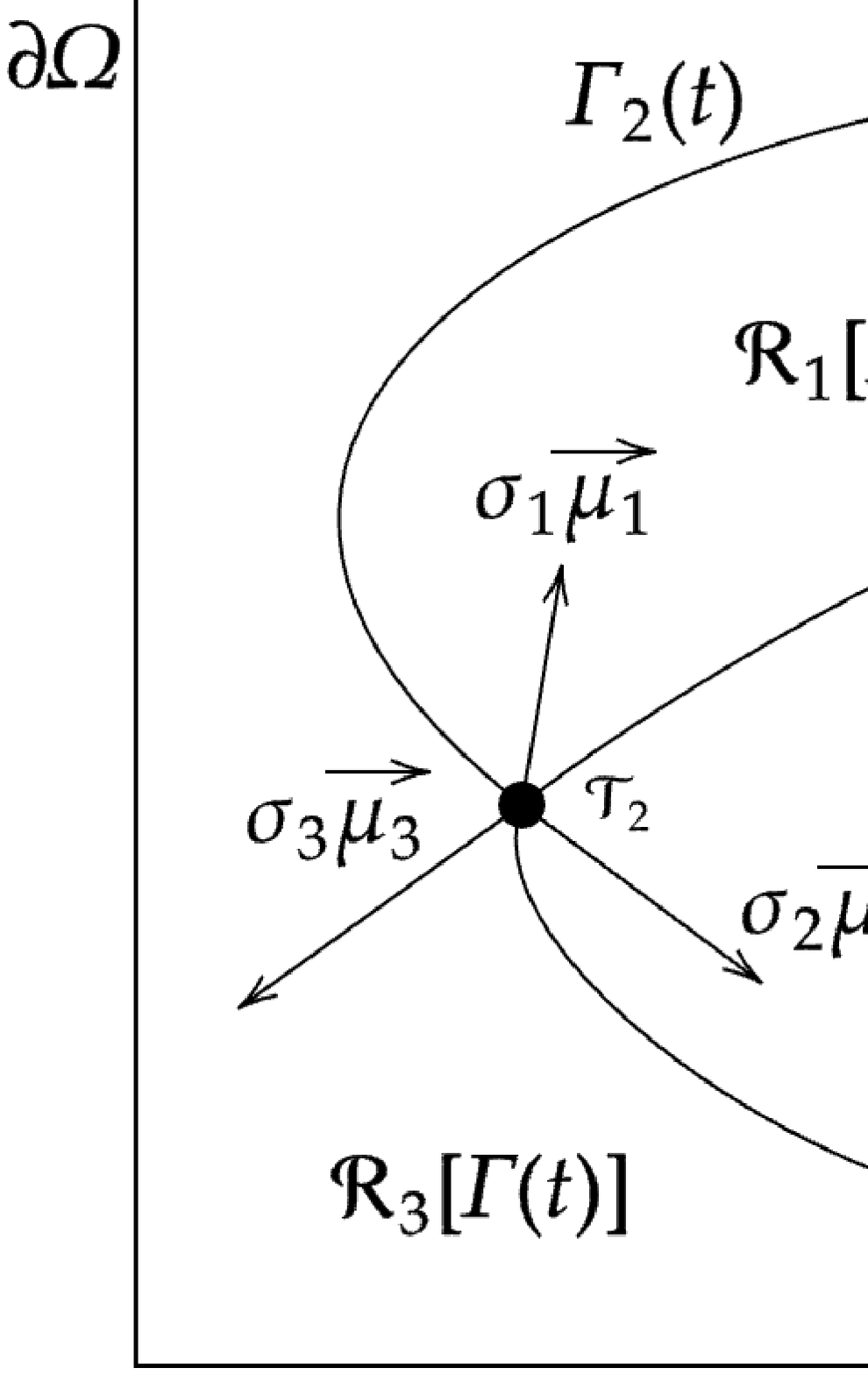}
    \caption{A surface cluster in 2d made up of three open curves and 
two triple junctions.}\label{fig:3p}
\end{figure}

For the representation of the evolving surface cluster $\Gamma(t)$, we adopt
the notation from \cite{GarckeNuernbergZhao2024}. See Figure~\ref{fig:3p}
for the sketch of an example in $\bR^2$.

The domain $\Omega$ is split into subdomains $\subdomains{\ell}\ (\ell\in\naturalindex{I_R})$ with $I_R\geq 2$
by a surface cluster $\Gamma(t) = \bigcup_{i=1}^{I_S}\Gamma_i(t)$ with $I_S\geq 1$,
i.e., $\Omega = \Gamma(t)\cup\bigcup_{\ell=1}^{I_R}\subdomains{\ell}$ with
$\subdomains{\ell_1}\cap\subdomains{\ell_2}=\emptyset$ if $\ell_1\neq \ell_2$.
Each $\Gamma_i(t)\ (i\in\naturalindex{I_S})$ is supposed to be either 
a closed hypersurface without boundary or a hypersurface with boundary in $\Omega$.
$\mT_k(t)\ (k\in\naturalindex{I_T})$ denotes triple junctions at which exactly three hypersurfaces
$\Gamma_{\tpindex{k}{1}}(t)$, $\Gamma_{\tpindex{k}{2}}(t)$ and $\Gamma_{\tpindex{k}{3}}(t)$
with $1\leq \tpindex{k}{1} < \tpindex{k}{2} < \tpindex{k}{3} \leq I_S$ meet at their boundary $\bigcap_{j=1}^3\partial\Gamma_{\tpindex{k}{j}}(t)$.
We define a mapping $b:\naturalindex{I_S}\ni i\mapsto(\cpindexp{i},\cpindexm{i})\in(\naturalindex{I_R})^2$
which means that the unit normal vector field $\normal{}{i}$ on $\Gamma_i(t)$
points from the region $\subdomains{\cpindexm{i}}$
into the region $\subdomains{\cpindexp{i}}$.
We note that the regions $\subdomains{\ell}$ are not necessarily connected.
To cope with this restriction, we require the mapping $b$ to satisfy
$(\cpindexp{i},\cpindexm{i}) = (\cpindexp{j},\cpindexm{j})$ if $\{\cpindexp{i},\cpindexm{i}\} = \{\cpindexp{j},\cpindexm{j}\}$.
This means that the normal direction is consistent even if two different interfaces separate the same two phases.

The parameters $\tension{i}\ (i\in\naturalindex{I_S})$ are positive constants
which denote the surface tensions of $\Gamma_{i}(t)$.
The parameters $\beta_\ell\ (\ell\in\naturalindex{I_R})$ denote constants
which, depending on the application, are related to the mass or energy
contents of the phases and satisfy $\beta_{\ell_1}\neq\beta_{\ell_2}$ if $\ell_1\neq \ell_2$.
Moreover, 
$\normal{}{\Omega}$ denotes the outward unit normal vector field on $\partial\Omega$;
$\meanCurvature{i}\,(i\in\naturalindex{I_S})$ is the scalar mean curvature of $\Gamma_i(t)$ in the direction of $\normal{}{i}$,
where we use the sign convention that $\meanCurvature{i} < 0$ if $\Gamma_i(t)$ is a closed convex surface with outer normal $\normal{}{i}$.
For a quantity $q$, $\jump{\cpindexm{i}}{\cpindexp{i}}{q}$ denotes
the jump of $q$ across $\Gamma_{i}(t)$ from the $\cpindexm{i}$-th phase to the $\cpindexp{i}$-th phase,
say $\jump{\cpindexm{i}}{\cpindexp{i}}{q}(x) := \lim_{\varepsilon\downarrow 0}\{q(x + \varepsilon\normal{}{i}(x)) - q(x - \varepsilon\normal{}{i}(x))\}$ for each $x\in\Gamma_{i}(t)$.
In addition, $\conormal{}{i}$ denotes the conormal, i.e.\ the intrinsic outer 
unit normal to $\partial\Gamma_i$, the boundary of $\Gamma_i$, 
that lies within the tangent plane of $\Gamma_i$,
 while $V_i$ is the normal velocity of $\Gamma_i(t)$ in the direction
$\vec\nu_i$. 
Finally,
the cluster $\Gamma_0$ represents initial data to close the system \eqref{eq:sharp_p}.

The system \eqref{eq:sharp_p} can be obtained as a gradient flow
of the surface energy of $\Gamma(t)$ as explained in \Section{sec:gs}.
In the physical point of view, each interface $\Gamma_{i}(t)$ is a free boundary between two phases, while $\smallChemical$ is a so-called chemical potential.
The fourth equation in \eqref{eq:sharp_p} describes force balance conditions on the triple junction lines $\mathcal{T}_k$, $k\in\naturalindex{I_T}$.
In the case of equal energy densities, they lead to the well known $120^\circ$ angle conditions at triple junction lines.
We observe that only differences of $\beta_\ell$ appear in the model
\eqref{eq:sharp_p}, and so from now on we assume without loss of generality 
that these coefficients are normalized to satisfy 
$\sum_{\ell=1}^{I_R}\beta_\ell = 0$.
%on replacing $\beta_i$ by $\beta_i - \sum_{\ell=1}^{I_R}\beta_\ell/I_R$;

We note that the system \eqref{eq:sharp_p} reduces to the two-phase Mullins--Sekerka model in the case $I_R = 2$.
Indeed, on letting $I_R = 2$, $I_S = 1$, $I_T = 0$, $(\cpindexp{1},\cpindexm{1}) = (1,2)$, $\beta_1 = -0.5$ and $\beta_2 = 0.5$,
we obtain that
\begin{equation}\label{eq:two_p}
    \begin{cases}
        \Delta\smallChemical = 0&\qquad\mbox{in}\quad\Omega\backslash\Gamma(t),\quad t \in (0,T),\\
        \jump{2}{1}{\beta}\smallChemical = -\tension{1}\meanCurvature{1}&\qquad\mbox{on}\quad\Gamma_{1}(t),\quad t \in (0,T),\\
        \jump{2}{1}{\beta}V_{1} = \jump{2}{1}{\nabla\smallChemical}\cdot\normal{}{1}&\qquad\mbox{on}\quad\Gamma_{1}(t),\quad t \in (0,T),\\
        \nabla\smallChemical\cdot\normal{}{\Omega} = 0&\qquad\mbox{on}\quad\partial\Omega,\quad t \in (0,T),\\
        \Gamma(0) = \Gamma_0,
    \end{cases}
\end{equation}
where $\subdomains{1}$ is the region enclosed by $\Gamma_1(t)$, and $\normal{}{1}$ is the inward unit normal vector field of $\Gamma_1(t)$.
Meanwhile, $\subdomains{2}$ is the exterior domain, say $\subdomains{2} = \Omega\backslash\closure{\subdomains{1}}$.
Since $\jump{2}{1}{\beta} = -1$, \eqref{eq:two_p} collapses to the classical two-phase Mullins--Sekerka problem, often also called Hele--Shaw problem,
\cite{EScherSimonett}.

% The review of the model
The sharp interface model that we consider in this paper is a degenerate version of the multi-phase Stefan problem,
and its distributional solution was first studied by Garcke and Sturzenhecker \cite{GarckeSturzenhecker1998}.
Following Luckhaus and Sturzenhecker \cite{LuckhausSturzenhecker}, they showed existence of a distributional solution to the problem \eqref{eq:sharp_p}
using a minimizing movement scheme of an energy functional tailored to the multi-phase flow.
The purpose of this paper is to establish a numerical scheme to compute their week solutions.

% SP-PFEM and the BGN method
The numerical scheme proposed in this paper is based on the parametric finite element method which is tailored
to deal with the appearance of triple junctions, which is the so-called BGN method \cite{BarrettGarckeRobert2007, BarrettGarckeNuernbergSurfaceDiffusion2007, ejam3d, clust3d}.
In the latter papers, the evolution governed by the surface diffusion and the mean curvature flow with triple junction was computed.
As discussed in \Section{sec:gs}, the multi-phase model that we consider possesses a gradient flow structure and 
conserves an energy-type quantity. It would be good if we could retain these structure and properties even for approximate solutions generated by
a fully discrete numerical scheme. For this purpose, Bao and Zhao \cite{BaoZhao2021} introduced a structure-preserving parametric finite element method (SP-PFEM)
for the surface diffuse flow. After that, the BGN method was combined with this structure-preserving numerical method by Bao, Garcke, N\"{u}rnberg and Zhao \cite{BaoGarckeNuernbergZhao2022},
and the surface diffuse flow with triple junctions could be handled without loss of the weighted mass conservation property.
This methodology was adapted to the two-phase Mullins--Sekerka problem and the multi-phase Mullins--Sekerka problem by N\"{u}rnberg \cite{Nurnberg202203} and the authors \cite{EtoGarckeNurnberg2024}, respectively.
We note that the mathematical model that the authors \cite{EtoGarckeNurnberg2024} dealt with is different from the model \eqref{eq:sharp_p} since
the chemical potential of each phase is expressed as a single scalar-valued function in \eqref{eq:sharp_p}.
For further recent works which analyzed geometric flows using the SP-PFEM,
we refer the reader to \cite{BaoLi2023,BaoLi2024,GarckeNuenbergZhao2023,GarckeNuernbergZhao2024,LiZhao2024}.

% Another numerical schemes than SP-PFEM and the BGN method to compute Mullins--Sekerka model
Various numerical methods have been developed to address the Mullins--Sekerka problem.
For the boundary integral method for the two-phase case, we refer the reader to \cite{BatesBrown,Bates1995ANS,ChenKublikTsai2017,Mayer2000,ZhuChenHou}.
A level-set formulation of moving boundaries together with the finite difference method was proposed in \cite{ChenMettimanOsherSmereka1997}.
For an implementation of the method of fundamental solutions for the two-phase Mullins--Sekerka problem in $\mathbb{R}^2$, see \cite{Eto2023}.
Combination of the method of fundamental solutions and a deep-learning based approach can be found in \cite{IzsakDjebbar2023},
although the main focus of the paper is approximating the Dirichlet-to-Neumann map with the help of neural networks.

% Numerical schemes for the diffuse interface model
We also mention some numerical methods for the diffuse interface model.
Numerical analysis of the scalar Cahn--Hilliard equation is dealt with in the works \cite{BarrettBlowey1995,BarrettBloweyGarcke1999,BloweyElliott1992,ElliottFrench1987}.
Feng and Prohl \cite{FengProhl2004} proposed a mixed fully discrete finite element method for the Cahn--Hilliard equation 
and provided error estimates between the solution of the Mullins--Sekerka problem and the approximate solution of the Cahn--Hilliard equation which are computed by their scheme.
The established error bounds yielded a convergence result in \cite{FengProhl2005}.
Aside from the sharp interface model, the Cahn--Hilliard equation for the multi-component case
has been computed in several works, see \cite{BarrettBloweyGarcke2001, BloweyCopettiElliott1996,Eyre1993,LiChoiKim2016,Nurnberg09}.
The multi-component Cahn--Hilliard equation on surfaces has recently been
considered in \cite{LiLiuXiaHeLi2022}. See \cite{BretinRolandMasnouSengersTerii2023}
for a study of the multi-component Cahn--Hilliard model involving concentration-dependent mobilities, and it adopted
the semi-implicit scheme together with the Fourier transform to compute a discrete solution to the geometric flow.

% analytical researches
We briefly review some analytical researches on the Mullins--Sekerka problem and its diffuse interface model.
For classical solutions, 
Chen, Hong and Yi \cite{ChenHongYi1996} showed the existence 
of a classical solution to the Mullins--Sekerka problem local-in-time in the two-dimensional case, while
Escher and Simonett \cite{EScherSimonett1998} gave a similar result in the general dimensional case.
For the notion of weak solutions, Luckhaus and Sturzenhecker \cite{LuckhausSturzenhecker}
established the existence of weak solutions to the two-phase Mullins--Sekerka problem in a distributional sense.
Therein, the weak solution was obtained as a limit of a sequence of time discrete approximate solutions
under the no mass loss assumption. The time implicit scheme is the basis of the approach in \cite{BronsardGarckeStoth1998}.
After that, R\"{o}ger \cite{Roger2005} removed the technical assumption of no mass loss in the case when the Dirichlet-Neumann
boundary condition is imposed by using geometric measure theory.
For this direction, we also find a convergence result of a minimizing movement scheme involving the Wasserstein distance of sets for the one-phase Mullins--Sekerka problem in the work by Chambolle and Laux \cite{ChambolleLaux2021}.
They regarded the gradient of the chemical potential as a solenoidal vector field and constructed a sequence of approximate set-theoretic solutions.
This sequence turned out to converge (up to a subsequence) to a distributional solution which was transported by the vector field.

Recently, researches which treat the boundary contact case
gradually appear. For the ninety-degree contact angle case, Abels, Rauchecker and Wilke \cite{AbelsMaxWilke} showed
well-posedness of the two-phase Mullins--Sekerka problem with the help of linearization approach.
Meanwhile, Garcke and Rauchecker \cite{GarckeRauchecker2022} carried out a stability analysis
in a curved domain in $\mathbb{R}^2$ via a linearization approach.
The ninety-degree contact angle condition was numerically computed by Eto \cite{Eto2023} in $\mathbb{R}^2$,
and its outcome exhibited that the curve shortening property and the area-preserving property were satisfied even in the discrete level.
Whereas, Hensel and Stinson \cite{HenselStinson2024} 
proposed a varifold solution to the two-phase Mullins--Sekerka problem with fixed contact angle condition
focussing on the energy dissipation property of the model.
Their definition of weak solutions was based on the maximum slope of the area functional in some topology.
Fischer, Hensel, Laux and Simon  \cite{FischerHenselLauxSimon2024} announced that this weak solution is unique
as long as a smooth solution to the two-phase Mullins--Sekerka problem exists, although their result was restricted to the two-dimensional case, and the contact angle condition was not treated.
For a gradient flow aspect of the Mullins--Sekerka flow, see e.g. \cite[\S 3.2]{Serfaty2011}.
A recent study by Escher, Matioc and Matioc \cite{EscherMatiocMatioc2024} showed well-posedness of the two-phase Mullins--Sekerka problem in the whole domain $\mathbb{R}^2$
in terms of the potential method. Their result supported the mathematical model which was considered in \cite{Bates1995ANS, Eto2023}.

% structure of the paper
This paper is organized as follows.
In \Section{sec:gs}, we focus on the gradient flow structure of the problem \eqref{eq:sharp_p}.
We also confirm that a classical solution to the problem \eqref{eq:sharp_p} satisfies
the dissipation property of surface energy and the total energy content preserving property.
In \Section{sec:wfm}, we will construct a weak formulation of the problem \eqref{eq:sharp_p}.
In \Section{sec:pfem}, we will define a discrete solution to the problem \eqref{eq:sharp_p} by using a parametric finite element method.
We also prove existence and uniqueness of the discrete solution and unconditional stability of the proposed scheme.
In \Section{sec:mf}, we will formulate the matrix form of the weak formulation to solve the linear system for the three-phase case.
In \Section{sec:spns}, we will introduce a vertex normal so that the fully discrete scheme conserves the total energy content on the discrete level.
In \Section{sec:ne}, we will carry out several numerical experiments to confirm feasibility of the proposed scheme and its accuracy.
Finally, we will conclude this paper in \Section{sec:con} by summarizing the main results.

\section{Gradient Flow and Energy Dissipation Properties}\label{sec:gs}
In a spirit of Garcke and Sturzenhecker \cite{GarckeSturzenhecker1998},
we begin by focusing on a gradient flow structure of the underlying problem \eqref{eq:sharp_p}.
Let $E(\Gamma(t))$ be the surface energy of $\Gamma(t)$ defined by
\begin{equation}\label{eq:surf_eg}
  E(\Gamma(t)) := \sum_{i=1}^{I_S}\int_{\Gamma_{i}(t)}\tension{i}\,\dH{d-1},
\end{equation}
where $\mH^{d-1}$ denotes the $(d-1)$-dimensional Hausdorff measure.

Let $m\in\mathbb{R}$ be a constant which indicates the total energy content.
Then, we define a function space $\mathcal{M}_m$ and its tangential space $T_{\Gamma}\mathcal{M}_m$ by
\begin{align*}
  \mathcal{M}_m &:= \left\{(\chi_1,\cdots,\chi_{I_R})\in BV(\Omega;\{0,1\})^{I_R}\biggm |\sum_{\ell=1}^{I_R}\int_{\subdomains{\ell}}\beta_{\ell}\chi_\ell\,\dL{d} =m\right\},\\
  T_{\Gamma}\mathcal{M}_m &:= \left\{V:\Gamma\to\mathbb{R}\biggm |\sum_{i=1}^{I_S} \int_{\Gamma_{i}}\jump{\cpindexm{i}}{\cpindexp{i}}{\beta}V_i\,\dH{d-1} = 0\right\},
\end{align*}
where $BV(\Omega;\{0,1\})$ denotes the set of all functions of bounded variation on $\Omega$ taking the values in $\{0,1\}$,
and $\mL^d$ denotes the $d$-dimensional Lebesgue measure.
Then, we define an inner product on $T_{\Gamma}\mathcal{M}_m$ by
\begin{equation*}
  \left<V_1,V_2\right>_{T_\Gamma\mathcal{M}_m} := \int_\Omega \nabla \smallChemical_1 \cdot\nabla \smallChemical_2\,\dL{d}
  = - \sum_{i=1}^{I_S}\int_{\Gamma_{i}}\smallChemical_1\jump{\cpindexm{i}}{\cpindexp{i}}{\beta}V_2\,\dH{d-1},
\end{equation*}
where $\smallChemical_k \in H^1(\Omega)\,(k=1,2)$ solves
\begin{equation}\label{eq:inner_p}
  \begin{cases}
    \Delta \smallChemical_k = 0&\qquad\mbox{in}\quad\Omega\backslash\Gamma,\\
    \jump{\cpindexm{i}}{\cpindexp{i}}{\nabla \smallChemical_k}\cdot\normal{}{i} = \jump{\cpindexm{i}}{\cpindexp{i}}{\beta}V_{k,i}&\qquad\mbox{on}\quad\Gamma_{i},\quad \forall i\in\naturalindex{I_S},\\
    \nabla \smallChemical_k\cdot\normal{}{\Omega} = 0&\qquad\mbox{on}\quad\partial\Omega.
  \end{cases}
\end{equation}
For the definition of this inner product in the two-phase case, we refer the reader to Garcke \cite{Garcke2013}.
The first variation of the energy $E(\Gamma)$ with respect to this metric on $T_{\Gamma}\mathcal{M}_m$ is given by
\begin{equation}\label{eq:var_eg}
  \nabla_{T_{\Gamma}\mathcal{M}_m}E(\Gamma) = -\DtN{\beta}{\Gamma}{\tension{}\curvature{}},
\end{equation}
where $\DtN{\beta}{\Gamma}{}$ denotes the Dirichlet-to-Neumann operator on $\Gamma$ defined by
\begin{equation*}
    \DtN{\beta}{\Gamma}{f} := g\qquad\mbox{for}\quad f\in H^{\frac{1}{2}}(\Gamma).
\end{equation*}
The operator $\DtN{\beta}{\Gamma}{f}$ is defined as follows.
For $f\in H^{\frac{1}{2}}(\Gamma)$, define $w\in H^1(\Omega)$ such that
\begin{equation}\label{eq:dn}
    \begin{cases}
        \Delta w = 0&\qquad\mbox{in}\quad\Omega\backslash\Gamma,\\
        \jump{\cpindexm{i}}{\cpindexp{i}}{\beta}w = f_i&\qquad\mbox{on}\quad\polygonCurve{}{i},\quad\forall i\in\naturalindex{I_S},\\
        \nabla w\cdot\normal{}{\Omega} = 0&\qquad\mbox{on}\quad\partial\Omega.
    \end{cases}
\end{equation}
Then we let $g\in H^{-\frac{1}{2}}(\Gamma)$ be such that
\begin{equation*}
    \jump{\cpindexm{i}}{\cpindexp{i}}{\nabla w}\cdot\normal{}{i} = \jump{\cpindexm{i}}{\cpindexp{i}}{\beta}g_i\qquad\mbox{on}\quad\polygonCurve{}{i},\quad\forall i\in\naturalindex{I_S}.
\end{equation*}
In the above, we have used the slight abuse of notation
$H^{\pm\frac12}(\Gamma) = \bigotimes_{i=1}^{I_S}H^{\pm\frac12}(\Gamma_i)$.
By the definition of $\DtN{\beta}{\Gamma}{}$, we see that the system \eqref{eq:sharp_p} can be interpreted as the geometric evolution equation:
\begin{equation*}
    V = -\DtN{\beta}{\Gamma(t)}{\tension{}\curvature{}}\qquad\mbox{on}\quad \polygonCurve{}{}(t),\quad t \in (0,T).
\end{equation*}
To observe \eqref{eq:var_eg}, suppose that $\smallChemical$ solves \eqref{eq:inner_p} with this $V$.
Assume that $\widetilde{\smallChemical}$ is a solution to \eqref{eq:dn} with $f = -\tension{}\curvature{}$,
and set $g = -\DtN{\beta}{\polygonCurve{}{}(t)}{\tension{}\curvature{}}$.
Then, we deduce from the transport theorem (see e.g., \cite[Theorem 3.4.2]{BanschDeckelnickGarckePozzi2023}) that
\begin{align*}
    \ddt E(\Gamma(t)) &= -\sum_{i=1}^{I_S} \int_{\polygonCurve{}{i}(t)}\tension{i}\curvature{i}V_i\,\dH{d-1} = \sum_{i=1}^{I_S}\int_{\polygonCurve{}{i}(t)}\jump{\cpindexm{i}}{\cpindexp{i}}{\beta}\widetilde{\smallChemical} V_i\,\dH{d-1}\\
    &= \sum_{i=1}^{I_S}\int_{\polygonCurve{}{i}(t)}\widetilde{\smallChemical}\jump{\cpindexm{i}}{\cpindexp{i}}{\nabla\smallChemical}\cdot\normal{}{i}\,\dH{d-1} = -\int_{\Omega}\nabla\widetilde{\smallChemical}\cdot\nabla\smallChemical\,\dL{d} = \inprod{\DtN{\beta}{\polygonCurve{}{}(t)}{\tension{}\curvature{}}}{V}_{T_{\polygonCurve{}{}(t)}\mathcal{M}_m}.
\end{align*}
Therefore, we obtain the first variation \eqref{eq:var_eg}.
In other words, the system \eqref{eq:sharp_p} also describes the gradient flow of the surface energy:
\begin{equation}\label{eq:gf}
    V = -\nabla_{T_{\Gamma(t)}\mathcal{M}_m}E(\Gamma(t))\qquad\mbox{on}\quad \Gamma(t),\quad t \in (0,T).
\end{equation}

We explore some properties of classical solutions of the multi-phase Stefan problem \eqref{eq:sharp_p}.
\begin{prop}[Energy dissipation]
    Let $(\smallChemical,\Gamma(t))$ be a classical solution of \eqref{eq:sharp_p}.
    Then, the surface energy of $\Gamma(t)$ is not increasing in time. Precisely, it holds that
    \begin{equation*}
        \ddt E(\Gamma(t)) = -\int_\Omega|\nabla\smallChemical|^2\,\dL{d}.
    \end{equation*}
\end{prop}
\begin{proof}
    We deduce from the transport theorem and the homogeneous Neumann boundary condition that
    \begin{align*}
        \ddt E(\Gamma(t)) &= -\sum_{i=1}^{I_S} \int_{\Gamma_i(t)}\tension{i}\meanCurvature{i}V_i\,\dH{d-1} = \sum_{i=1}^{I_S}\int_{\Gamma_i(t)}\jump{\cpindexp{i}}{\cpindexm{i}}{\beta}\smallChemical V_i\,\dH{d-1}\\
        &= \sum_{i=1}^{I_S}\int_{\Gamma_{i}(t)}\smallChemical\jump{\cpindexm{i}}{\cpindexp{i}}{\nabla \smallChemical}\cdot\normal{}{i}\,\dH{d-1} = -\sum_{\ell=1}^{I_R}\int_{\subdomains{\ell}}\left|\nabla\smallChemical\right|^2\,\dL{d} \\ & 
= -\int_\Omega|\nabla\smallChemical|^2\,\dL{d}.
    \end{align*}
\end{proof}
\begin{prop}[Total energy content preservation]\label{prop:mass_p}
    Let $(w,\Gamma(t))$ be a classical solution of \eqref{eq:sharp_p}.
    Then, the total energy content is preserved in time. Namely, it holds that
    \begin{equation*}
        \ddt\sum_{\ell=1}^{I_R}\beta_\ell\volume{\subdomains{\ell}} = 0,
    \end{equation*}
    where $\volume{\cdot} := \mL^d({\cdot})$.
\end{prop}
\begin{proof}
    For each $\ell\in\naturalindex{I_R}$, let $\mV_\ell$ be the normal velocity of $\partial\subdomains{\ell}$ inward to $\subdomains{\ell}$.
    Then, we have
    \begin{align}\label{eq:mass_p_1}
        \ddt \beta_\ell\volume{\subdomains{\ell}} &= -\beta_\ell\int_{\partial\subdomains{\ell}}\mV_\ell\,\dH{d-1}\nonumber \\
        &= -\beta_\ell\sum_{\substack{i\in\naturalindex{I_S}\\\cpindexp{i} = \ell}} \int_{\Gamma_i(t)}V_i\,\dH{d-1} + \beta_\ell\sum_{\substack{i\in\naturalindex{I_S}\\\cpindexm{i} = \ell}} \int_{\Gamma_i(t)}V_i\,\dH{d-1}.
    \end{align}
    Summing up \eqref{eq:mass_p_1} over all $\ell\in\naturalindex{I_R}$ together with the homogeneous Neumann boundary condition, we have
    \begin{align*}
        \ddt \sum_{\ell=1}^{I_R}\beta_\ell\volume{\subdomains{\ell}} &= -\sum_{i=1}^{I_S}\int_{\Gamma_i(t)}\jump{\cpindexm{i}}{\cpindexp{i}}{\beta} V_i\,\dH{d-1} = -\sum_{i=1}^{I_S}\int_{\Gamma_i(t)}\jump{\cpindexm{i}}{\cpindexp{i}}{\nabla\smallChemical}\cdot\normal{}{i}\,\dH{d-1}\\
        &= \sum_{\ell=1}^{I_R}\int_{\subdomains{\ell}}\Delta\smallChemical\,\dL{d} = 0.
    \end{align*}
\end{proof}
\begin{rem}\label{rem:ap}
In contrast to the multi-phase model which was dealt with in \cite{EtoGarckeNurnberg2024},
we observe that the system \eqref{eq:sharp_p} does not conserve the mass of each of the phases.
In other words, the mathematical model \eqref{eq:sharp_p} is essentially different from that in \cite[Eq.(1.2)]{EtoGarckeNurnberg2024}.
\end{rem}

\section{Weak formulation}\label{sec:wfm}
In this section, we establish a weak formulation for the degenerate multi-phase Stefan problem \eqref{eq:sharp_p}.
We take a test function $\varphi\in H^1(\Omega)$ and multiply the first equation of \eqref{eq:sharp_p} by $\varphi$.
For each $\ell\in\naturalindex{I_R}$, we have
\begin{align}\label{eq:wkf_1}
    \int_{\subdomains{\ell}}\Delta \smallChemical\,\varphi\,\dL{d} &= \sum_{\substack{i\in\naturalindex{I_S}\\\cpindexp{i} = \ell}}\int_{\Gamma_i(t)}-(\nabla\smallChemical\cdot\normal{}{i})\varphi\,\dH{d-1} + \sum_{\substack{i\in\naturalindex{I_S}\\\cpindexm{i} = \ell}}\int_{\Gamma_i(t)}(\nabla\smallChemical\cdot\normal{}{i})\varphi\,\dH{d-1}\nonumber\\
    &\quad - \int_{\subdomains{\ell}}\nabla\smallChemical\cdot\nabla\varphi\,\dL{d}.
\end{align}
Summing up \eqref{eq:wkf_1} over all $\ell\in\naturalindex{I_R}$, we have
\begin{align*}
    0 &= \sum_{i=1}^{I_S} \int_{\Gamma_i(t)}(-\jump{\cpindexm{i}}{\cpindexp{i}}{\nabla\smallChemical}\cdot\normal{}{i})\varphi\,\dH{d-1} - \int_\Omega\nabla\smallChemical\cdot\nabla\varphi\,\dL{d}\\
    &= \sum_{i=1}^{I_S}\int_{\Gamma_i(t)}-\jump{\cpindexm{i}}{\cpindexp{i}}{\beta}V_i\varphi\,\dH{d-1} - \int_\Omega\nabla\smallChemical\cdot\nabla\varphi\,\dL{d},
\end{align*}
where we have used the third equation of the system \eqref{eq:sharp_p} to transform $\jump{\cpindexm{i}}{\cpindexp{i}}{\nabla w}\cdot\normal{}{i}$ to $\jump{\cpindexm{i}}{\cpindexp{i}}{\beta}V_i$.
%Thus, we have the following weak formulation of the diffuse equation together with the motion law:
%\begin{equation*}
%    \sum_{i=1}^{I_S}\int_{\Gamma_{i}(t)}\jump{\cpindexm{i}}{\cpindexp{i}}{\beta}V_{i}\varphi\,\dH{d-1} + \int_\Omega\nabla \smallChemical\cdot\nabla\varphi\,\dL{d} = 0\qquad\forall\varphi\in H^1(\Omega).
%\end{equation*}
Meanwhile, a weak formulation for the Gibbs--Thomson law of the system \eqref{eq:sharp_p} is straightforward:
\begin{equation*}
    \sum_{i=1}^{I_S}\int_{\Gamma_{i}(t)}(\jump{\cpindexm{i}}{\cpindexp{i}}{\beta}\smallChemical + \tension{i}\meanCurvature{i})\xi\,\dH{d-1} = 0\qquad\forall\xi\in L^2(\Gamma(t)).
\end{equation*}
Finally, let us obtain a weak formulation for the curvature vector.
Let $\Vec\id:\mathbb{R}^d\to\mathbb{R}^d$ be the identity function.
Then, we know the following formula for the curvature vector:
\begin{equation}\label{eq:wfm_curv}
  \sLaplacian\Vec\id = \meanCurvature{i}\normal{}{i}\qquad\mbox{on}\quad\Gamma_i(t),
\end{equation}
where $\sLaplacian$ denotes the Laplace--Beltrami operator defined by $\sgrad\cdot\sgrad$ by using the surface gradient $\sgrad$.
We now introduce a test function space for \eqref{eq:wfm_curv} by
\begin{equation*}
    \underline{V}(\polygonCurve{}{}(t)) := \left\{(\vec{\eta}_1,\ldots,\vec{\eta}_{I_S})\in \bigotimes_{i=1}^{I_S}H^1(\polygonCurve{}{i}(t);\mathbb{R}^d)\biggm| \vec{\eta}_{\tpindex{k}{1}} = \vec{\eta}_{\tpindex{k}{2}} = \vec{\eta}_{\tpindex{k}{3}}\quad\mbox{on}\quad\mathcal{T}_k(t),\quad\forall k\in\naturalindex{I_T}\right\}.
\end{equation*}
Multiplying \eqref{eq:wfm_curv} with $\tension{i}\vec\eta_i$, for
an arbitrary $\Vec{\eta}\in \underline{V}(\Gamma(t))$, integrating over
$\Gamma_i(t)$, summing and performing integration by parts yields
\begin{align*}
    &\sum_{i=1}^{I_S}\int_{\Gamma_i(t)}\tension{i}\meanCurvature{i}\normal{}{i}\cdot\Vec{\eta}_i\,\dH{d-1} = \sum_{i=1}^{I_S}\int_{\Gamma_i(t)}\tension{i}\sLaplacian\Vec\id\cdot\Vec{\eta}_i\,\dH{d-1}\\
    &= \sum_{i=1}^{I_S}\left(\int_{\partial\Gamma_i(t)}\tension{i}(\sgrad\Vec\id\,\Vec{\eta}_i)\cdot\conormal{}{i}\,\dH{d-2} - \int_{\Gamma_i(t)}\tension{i}\sgrad\Vec\id\cdot\sgrad\Vec{\eta}_i\,\dH{d-1}\right)\\
    &= \sum_{k=1}^{I_T}\int_{\mT_k(t)}\sum_{j=1}^3\tension{\tpindex{k}{j}}(\sgrad\Vec\id\,\Vec{\eta}_{\tpindex{k}{j}})\cdot\conormal{}{\tpindex{k}{j}}\,\dH{d-2} - \sum_{i=1}^{I_S}\int_{\Gamma_i(t)}\tension{i}\sgrad\Vec\id\cdot\sgrad\Vec{\eta}_i\,\dH{d-1}\\
    &= \sum_{k=1}^{I_T}\int_{\mT_k(t)}(\sgrad\Vec\id\,\Vec{\eta}_{\tpindex{k}{1}})\sum_{j=1}^3\tension{\tpindex{k}{j}}\conormal{}{\tpindex{k}{j}}\,\dH{d-2} - \sum_{i=1}^{I_S}\int_{\Gamma_i(t)}\tension{i}\sgrad\Vec\id\cdot\sgrad\Vec{\eta}_i\,\dH{d-1}\\
    &= -\sum_{i=1}^{I_S}\int_{\Gamma_i(t)}\tension{i}\sgrad\Vec\id\cdot\sgrad\Vec{\eta}_i\,\dH{d-1}.
\end{align*}
For the integrals on the triple junctions $\mT_k(t)$ we have first used that
$\Vec{\eta}\in \underline{V}(\Gamma(t))$, and then noticed that they vanish due
to Young's law, which is the fourth condition of the system \eqref{eq:sharp_p}.
%Hence, we have the following weak formulation for the mean curvature vector:
%For every $\Vec{\eta}\in H^1(\Gamma(t);\mathbb{R}^d)$ which is continuous on $\bigcup_{k=1}^{I_T}\mT_k(t)$,
%\begin{equation*}
%    \int_{\Gamma(t)}\meanCurvature{\tension{}}\normal{}{}\cdot\Vec{\eta}\,\dH{d-1} + \int_{\Gamma(t)}\tension{}\sgrad\Vec\id\cdot\sgrad\Vec{\eta}\,\dH{d-1} = 0.
%\end{equation*}
For later use, we define the inner products in $\Omega$ and on $\Gamma(t)$ by
\begin{equation*}
    \inprod{u}{v}_\Omega := \int_\Omega u\,v\,\dL{d}
%\qquad\mbox{for} \quad u,v\in L^2(\Omega)
%\end{equation*}
\quad\text{and}\quad
%\begin{equation*}
    \inprod{u}{v}_{\Gamma(t)} := \sum_{i=1}^{I_S} \int_{\Gamma_i(t)}u_i\,v_i\,\dH{d-1}.
%\qquad\mbox{for}\quad u,v\in L^2(\Gamma).
\end{equation*}
We conclude this section by summarizing the weak formulation of the system,
where from now on we use the short-hand notation $\meanCurvature{\tension{}}$
defined by $\meanCurvature{\tension{}}\mid_{\Gamma_i(t)} = \tension{i}\meanCurvature{i}$ for $i\in\naturalindex{I_S}$.
The weak formulation of problem \eqref{eq:sharp_p}
amounts to finding a time-dependent pair $(\smallChemical,\Gamma(t))$
of a diffuse function $\smallChemical$ and a surface cluster $\Gamma(t)$
such that for each $t \in (0,T)$:
\newline\newline
\begin{subequations}\label{eq:wkf}
    \textbf{Motion law: } For all $\varphi\in H^1(\Omega)$, it holds that
    \begin{equation}\label{eq:wkf_dm}
        \inprod{\nabla\smallChemical}{\nabla\varphi}_\Omega + \sum_{i=1}^{I_S}\jump{\cpindexm{i}}{\cpindexp{i}}{\beta}\inprod{V_i}{\varphi}_{\Gamma_i(t)} = 0.
    \end{equation}
    \textbf{Gibbs--Thomson law: } For all $\xi\in L^2(\Gamma(t))$, it holds that
    \begin{equation}\label{eq:wkf_ml}
        \sum_{i=1}^{I_S}\inprod{\jump{\cpindexm{i}}{\cpindexp{i}}{\beta}\smallChemical + \meanCurvature{\tension{}}}{\xi}_{\Gamma_i(t)} = 0.
    \end{equation}
    \textbf{Curvature vector: } For all $\vec{\eta} \in \underline{V}(\Gamma(t))$,
    it holds that
    \begin{equation}\label{eq:wkf_cv}
        \inprod{\meanCurvature{\tension{}}\normal{}}{\vec{\eta}}_{\Gamma(t)} + \inprod{\tension{}\sgrad\vec\id}{\sgrad\vec{\eta}}_{\Gamma(t)} = 0.
    \end{equation}
\end{subequations}

\section{Parametric finite element method}\label{sec:pfem}
In this section, we propose a parametric finite element method for the degenerate multi-phase Stefan problem \eqref{eq:sharp_p} based on the weak formulation \eqref{eq:wkf}.
For the presentation of the necessary finite element spaces we closely follow
the notation from \cite{BaoGarckeNuernbergZhao2022}, see also \cite{clust3d}. 

We split the interval $[0,T]$ into $M$ sub-intervals $[t_{m-1},t_m]$ with $1\leq m\leq M$ whose lengths are equal to $\tau_m$.
Then, we construct a fully discrete scheme which yields an approximate solution $(\largeChemical^{m+1},\vVertex{m+1}{},\discreteCurvature{m+1}{\tension{}})$,
provided that the previous discrete cluster $\Gamma^m$ is given,
with $\Gamma^{m+1} = \vVertex{m+1}{}(\Gamma^m)$.
In order to describe $\Gamma^m$, and the discrete matching conditions that have
to hold on the triple junction, we 
let $\refsurface{i}^h\,(i\in \naturalindex{I_S})$ be polyhedral reference surfaces with $\closure{\refsurface{i}^h} = \bigcup_{j=1}^{\numSimplex{}{i}}\closure{\simplex{}{i,j}}$, where $\{\simplex{}{i,j}\}_{j=1}^{J_i}$ is a family of
 mutually disjoint open $(d-1)$-simplices with vertices
$\{\vVertexSynonym{}{i,k}\}_{k=1}^{K_i}$ vertices.

Moreover, we assume that each boundary $\partial\refsurface{i}^h$ is split into $I_P^i$ sub-boundaries $\partial_p\refsurface{i}^h\,(p\in\naturalindex{I_P^i})$,
and each sub-boundary $\partial_p\refsurface{i}^h$ corresponds to the 
parameterization of a triple junction.
In particular, we will let 
$\Gamma^m_i = \Vec{\mathfrak{X}}^m_i(\refsurface{i}^h)$, so that the
triple junction $\mathcal{T}_k$, where the surfaces
$\Gamma_{s^k_1}(t)$, $\Gamma_{s^k_2}(t)$, $\Gamma_{s^k_3}(t)$ meet, 
is approximated by the images of $\Vec{\mathfrak{X}}^m$ on
$\partial_{p_1^k} \refsurface{s_1^k}^h$, 
$\partial_{p_2^k} \refsurface{s_2^k}^h$, 
$\partial_{p_3^k} \refsurface{s_3^k}^h$.
To this end, we have to ensure that these sub-boundaries perfectly match up on 
the triple junctions, and in particular contain the same number of vertices.
Hence, we assume that for every $k\in\naturalindex{I_T}$, it holds that
\begin{equation}\label{eq:q}
    \numTjVertex{k} := \#Q_{s_1^k,p^k_1} = \#Q_{s_2^k,p^k_2} = 
\#Q_{s_3^k, p^k_3},
\end{equation}
where
$Q_{s,p} := \left\{\vVertexSynonym{}{s,\ell}\right\}_{\ell=1}^{\numSurfVertex{}{s}}\cap \partial_p\refsurface{s}^h$
denotes the set of vertices belonging to the boundary patch
$\partial_p\refsurface{s}^h$.
Then we assume in addition that there exist bijections $\orderedSeq{k}{r}:\,\naturalindex{\numTjVertex{k}}\to Q_{s_r^k,p^k_r}\,(r=1,2,3)$
such that $(\orderedSeq{k}{r}(1),\ldots,\orderedSeq{k}{r}(\numTjVertex{k}))\,(r=1,2,3)$ are ordered sequences of the vertices.

\begin{figure}[H]
\centering
\includegraphics[angle=-90,width=0.38\textwidth]{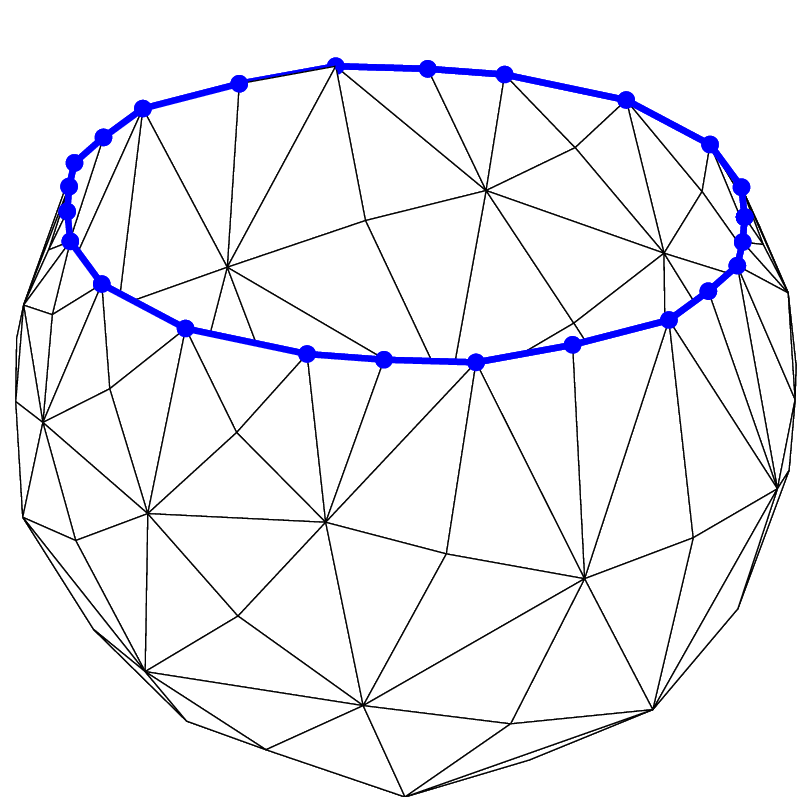}
%$\blue \curvearrowright$
\includegraphics[angle=-90,width=0.19\textwidth]{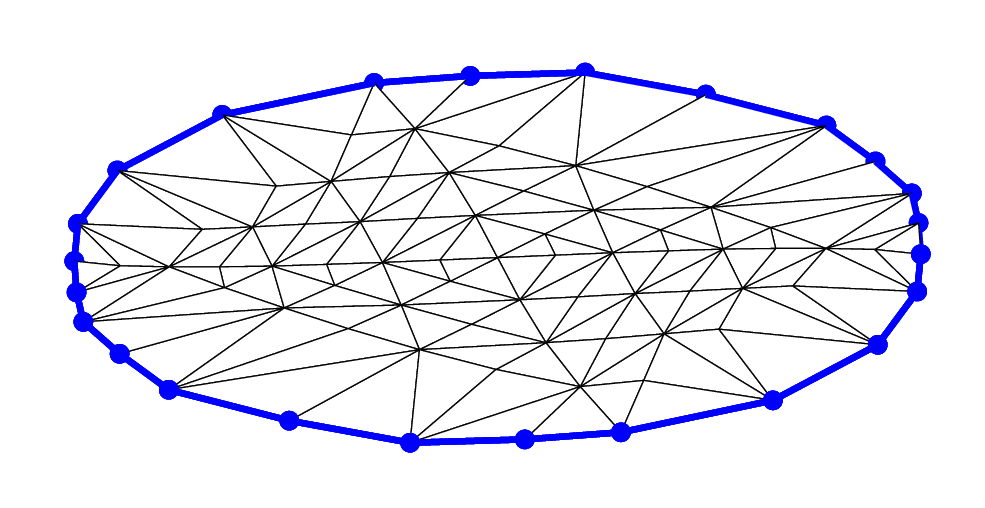}
%\boldmath{$\blue \curvearrowleft$}
\includegraphics[angle=-90,width=0.38\textwidth]{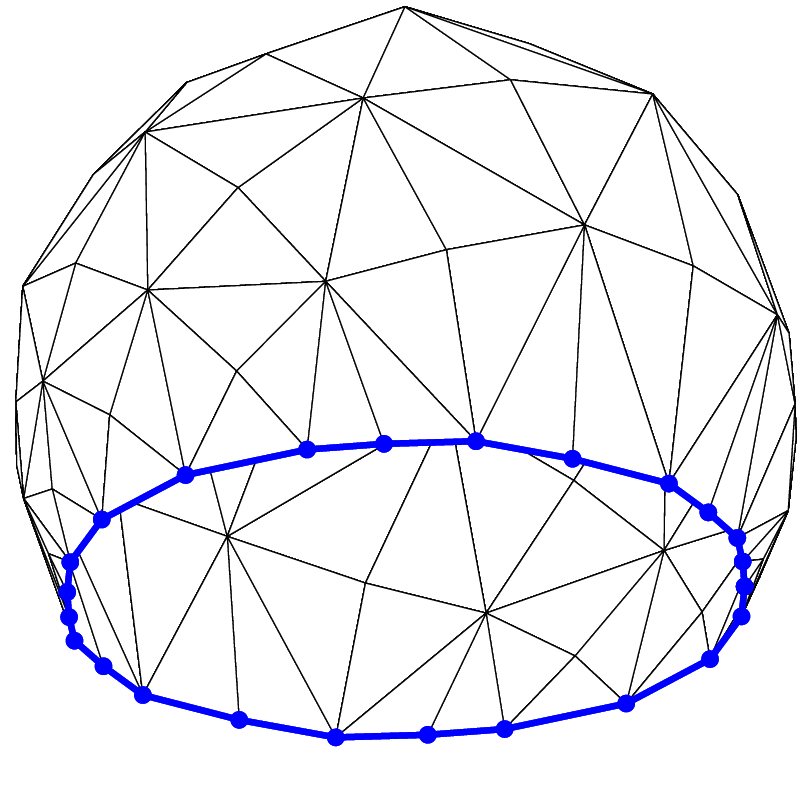}
\setlength{\abovecaptionskip}{20pt}
\caption{Example triangulation $\Gamma^m$ of a triple bubble with polygonal 
triple junction line $\mathcal{T}^m_1$, where here $Z_1 = 24$.}
\label{fig:kirchzarten}
\end{figure}

Let
\begin{multline*}
    \underline{V}^h(\refsurface{}^h) := \bigg\{(\polyhedral{}{1},\cdots,\polyhedral{}{I_S})\in\bigotimes_{i=1}^{I_S}C(\closure{\refsurface{i}^h};\bR^d)\biggm|\polyhedral{}{i}\mid_{\simplex{}{i,j}}\mbox{ is affine}\quad\forall j\in \naturalindex{\numSimplex{}{i}},\,\forall i\in\naturalindex{I_S},\\
    \mbox{and}\qquad\polyhedral{}{\tpindex{k}{1}}(\orderedSeq{k}{1}(z)) = \polyhedral{}{\tpindex{k}{2}}(\orderedSeq{k}{2}(z)) = \polyhedral{}{\tpindex{k}{3}}(\orderedSeq{k}{3}(z))\qquad\forall z\in\naturalindex{\numTjVertex{k}},\,\forall  k\in\naturalindex{I_T}\bigg\}.
\end{multline*}
Then, for each $m\geq 0$ and $\Vec{\mathfrak{X}}^m\in V^h(\refsurface{}^h)$,
we define $\Gamma^m := \Vec{\mathfrak{X}}^m(\refsurface{}^h)$ with
$\polygonCurve{m}{i} = \Vec{\mathfrak{X}}^m_i(\refsurface{i}^h)$,
$\simplex{m}{i,j} := \polyhedral{m}{i}(\simplex{}{i,j})$ and
$\vVertexSynonym{m}{i,k} := \polyhedral{m}{i}(\vVertexSynonym{}{i,k})$.
The discrete triple junctions $\mT^m_k\,(k\in\naturalindex{I_T})$
are defined by
%\begin{align*}
$\mT^m_k := \left\{ \polyhedral{m}{\tpindex{k}{1}}(\orderedSeq{k}{1}(z)) \biggm| z\in \naturalindex{\numTjVertex{k}}\right\}$. %\\
%    &\left(=  \left\{ \polyhedral{m}{\tpindex{k}{2}}(\orderedSeq{k}{2}(z)) \biggm| z\in \naturalindex{\numTjVertex{k}}\right\} =  \left\{ \polyhedral{m}{\tpindex{k}{3}}(\orderedSeq{k}{3}(z)) \biggm| z\in \naturalindex{\numTjVertex{k}}\right\}\right).
%\end{align*}
We visualize an example discrete surface cluster in Figure~\ref{fig:kirchzarten}.

Let us define finite element spaces from where we will seek approximate solutions.
Let $\triangulation{m}$ be a triangulation of $\closure{\Omega}$. Then, an approximation solution $\largeChemical^{m+1}$ is sought in the finite element space defined by
\begin{equation*}
    \femspaceChemical{m} := \left\{v\in C(\closure{\Omega})\biggm| v\mid_o\ \mbox{is affine}\ \forall o\in\triangulation{m}\right\}.
\end{equation*}

On the polyhedral surface $\Gamma^m$ we introduce finite element spaces defined by
\begin{align*}
    V^h(\polygonCurve{m}{i}) &:= \left\{v\in C(\polygonCurve{m}{i})\biggm|v\mid_{\simplex{m}{i,j}}\mbox{ is affine}\quad \forall j\in\naturalindex{\numSimplex{}{i}}\right\},\qquad i\in\naturalindex{I_S},\\
    \underline{V}^h(\polygonCurve{m}{i}) &:= \left\{\bv v = (v_1,\cdots,v_d)\in C(\polygonCurve{m}{i};\mathbb{R}^d)\biggm|v_k\in V^h(\polygonCurve{m}{i})\quad \forall k\in \naturalindex{d}\right\},\qquad i\in\naturalindex{I_S}.
\end{align*}
For later use, we also define 
\begin{equation} \label{eq:defV0}
    V^h_0(\polygonCurve{m}{i}) := \left\{v\in V^h(\polygonCurve{m}{i})\biggm|\, v = 0\quad\mbox{on}\quad \partial\polygonCurve{m}{i}\right\},
\end{equation}
and let $\{\femspaceBasisCurve{m,0}{i,k}\}_{k=1}^{K^0_i}$ be the standard basis of $V^h_0(\polygonCurve{m}{i})$, so that 
$\femspaceBasisCurve{m,0}{i,k}(\vec q^m_{i,\ell}) = \delta_{k\ell}$, 
$k,\ell \in \naturalindex{K^0_i}$.

Then, the approximate solutions $\vVertex{m+1}{}$ and $\kappa_{\tension{}}^{m+1}$ are respectively sought in the finite element spaces defined by
\begin{align*}
    \femspaceVertexVector{h}{}(\Gamma^m) &:= \left\{(\vVertex{}{1},\cdots,\vVertex{}{I_S})\in\bigotimes_{i=1}^{I_S}\underline{V}^h(\polygonCurve{m}{i})\biggm|\vVertex{}{\tpindex{k}{1}} = \vVertex{}{\tpindex{k}{2}} = \vVertex{}{\tpindex{k}{3}}\quad\mbox{on}\quad\mT_k^m,\ \forall k\in\naturalindex{I_T}\right\},\\
    \femspaceVertex{h}{}(\Gamma^m) &:= \bigotimes_{i=1}^{I_S}V^h(\polygonCurve{m}{i}).
\end{align*}
We now define the normal vector of each simplex $\simplex{m}{i,j}$.
To this end, let
$\left\{\vec q_{i,j,\ell}^{m}\right\}_{\ell=1}^{d}$ be the vertices of 
$\sigma_{i,j}^{m}$, and ordered with the same orientation for all 
$\sigma_{i,j}^{m}$, $j\in\naturalindex{J_i}$. Then we define
\begin{equation*}
    \normal{m}{i,j} := \frac{\Normal{\simplex{m}{i,j}}}{|\Normal{\simplex{m}{i,j}}|}\quad\mbox{with}\quad\Normal{\simplex{m}{i,j}} := \begin{cases}
       (\vVertexSynonym{m}{i,j,2} - \vVertexSynonym{m}{i,j,1})^\perp\quad\mbox{if}\quad d = 2,\\
       (\vVertexSynonym{m}{i,j,2} - \vVertexSynonym{m}{i,j,1})\wedge(\vVertexSynonym{m}{i,j,3} - \vVertexSynonym{m}{i,j,1})\quad\mbox{if}\quad d = 3,
    \end{cases}
\end{equation*}
where $|\cdot| = \mH^{d-1}(\cdot)$;
the symbol $\wedge$ denotes the wedge product, and $v^\perp := (-v_2,v_1)^T$ for $v = (v_1,v_2)^T\in\mathbb{R}^2$.
Let $\normal{m}{i}$ be the normal vector on $\polygonCurve{m}{i}$ which equals $\normal{m}{i,j}$ on $\simplex{m}{i,j}$.

Let us define the mass lumped inner product of two piecewise continuous functions $u$ and $v$ on $\polygonCurve{m}{i}$ by
\begin{equation*}
    \inprodMassLamped{u}{v}{h}{\polygonCurve{m}{i}} := \frac{1}{d}\sum_{j=1}^{\numSimplex{}{i}}|\simplex{m}{i,j}|\sum_{k=1}^d\lim_{\substack{\simplex{m}{i,j}}\ni\vec{q}\to\vVertexSynonym{m}{i,j,k}}(uv)(\vec{q}) \qquad i\in\naturalindex{I_S}.
\end{equation*}
Using this, we define the mass lumped inner product on $\Gamma^m$ by
\begin{equation}\label{eq:mlip}
    \inprodMassLamped{u}{v}{h}{\Gamma^m} := \sum_{i=1}^{I_S}\inprodMassLamped{u_i}{v_i}{h}{\polygonCurve{m}{i}}.
\end{equation}
Meanwhile, we will write the natural $L^2$--inner product as follows:
\begin{equation*}
    \inprodMassLamped{u}{v}{}{\polygonCurve{m}{}} = \sum_{i=1}^{I_S}\inprodMassLamped{u_i}{v_i}{}{\polygonCurve{m}{i}} =\sum_{i=1}^{I_S} \int_{\polygonCurve{m}{i}} u_i v_i\,\dH{d-1}.
\end{equation*}
The notion of these inner products can be extended
for two vector- and tensor-valued functions.
The vertex normal $\femNormalVertex{m}{i}\in \underline V^h(\polygonCurve{m}{i})$ on $\polygonCurve{m}{i}$ 
is defined in terms of the $L^2$ projection as follows:
\begin{equation*}
    \inprodMassLamped{\femNormalVertex{m}{i}}{\Vec{\xi}}{h}{\polygonCurve{m}{i}} = \inprodMassLamped{\normal{m}{i}}{\Vec{\xi}}{}{\polygonCurve{m}{i}},\quad\forall\Vec{\xi}\in \underline V^h(\polygonCurve{m}{i})\qquad i\in\naturalindex{I_S}.
\end{equation*}

We are now in the position to introduce our finite element approximation
for the weak formulation \eqref{eq:wkf}. 

For $m\geq 0$, given $\Gamma^m$,
find $(\largeChemical^{m+1},\vVertex{m+1}{},\discreteCurvature{m+1}{\tension{}})\in\femspaceChemical{m}\times\femspaceGammaVector{m}\times\femspaceGamma{m}$
such that the following three conditions hold:\newline\newline
    \noindent
    \begin{subequations}\label{eq:fma}
    \textbf{Motion law: } For all $\varphi\in\femspaceChemical{m}$, it holds that
    \begin{equation}\label{eq:fma_ml}
        \inprodMassLamped{\nabla\largeChemical^{m+1}}{\nabla\varphi}{}{\Omega} + \sum_{i=1}^{I_S} \jump{\cpindexm{i}}{\cpindexp{i}}{\beta} \inprodMassLamped{\interpolation{h}{i}\left[\frac{\vVertex{m+1}{i} - \Vec\id}{\tau_m}\cdot\femNormalVertex{m}{i}\right]}{\varphi}{(h)}{\polygonCurve{m}{i}} = 0.
    \end{equation}
    \textbf{Gibbs--Thomson law: } For all $\xi\in\femspaceGamma{m}$, it holds that
    \begin{equation}\label{eq:fma_gtl}
        \inprodMassLamped{\discreteCurvature{m+1}{\tension{}}}{\xi}{h}{\polygonCurve{m}{}} + \sum_{i=1}^{I_S}\jump{\cpindexm{i}}{\cpindexp{i}}{\beta}\inprodMassLamped{\largeChemical^{m+1}}{\xi}{(h)}{\polygonCurve{m}{i}} = 0.
    \end{equation}
    \textbf{Curvature vector: } For all $\Vec{\eta}\in\femspaceGammaVector{m}$, it holds that
    \begin{equation}\label{eq:fma_cv}
        \inprodMassLamped{\discreteCurvature{m+1}{\tension{}}\femNormalVertex{m}{}}{\Vec{\eta}}{h}{\Gamma^m} + \inprodMassLamped{\tension{}\sgrad\vVertex{m+1}{}}{\sgrad\Vec{\eta}}{}{\polygonCurve{m}{}} = 0,
    \end{equation}
    \end{subequations}
where $\interpolation{h}{i}:C(\polygonCurve{m}{i})\to V^h(\polygonCurve{m}{i})$ denotes the standard interpolation operator, and set 
$\Gamma^{m+1} = \vVertex{m+1}{}(\Gamma^m)$.
Above, and in what follows, the notation $\cdot^{(h)}$ means an expression with or
without the superscript $h$. That is, \eqref{eq:fma} represents two different schemes.
When the superscript $h$ appears, we use the mass lumped inner product introduced in \eqref{eq:mlip},
while the inner products without the superscript $h$ are computed by true integration.

Before proving existence and uniqueness of \eqref{eq:fma},
we impose a technical assumption on the surface cluster.
This assumption is analogous to that in \cite[Assumption~1]{EtoGarckeNurnberg2024}.
\begin{asm}\label{asm:1}
    For every $i\in\naturalindex{I_S}$, it holds that
    \begin{equation*}
        \operatorname{span}\left\{\femNormalVertex{m}{i}(\vec q^m_{i,k}) \biggm| j\in\naturalindex{K^0_i}\right\} \neq \{\zerovec\}.
    \end{equation*}
    Moreover, we assume that
    \begin{equation*}
        \operatorname{span}\left\{\sum_{i=1}^{I_S} \jump{\cpindexm{i}}{\cpindexp{i}}{\beta}\inprodMassLamped{\femNormalVertex{m}{i}}{\varphi}{(h)}{\polygonCurve{m}{i}}\biggm|\varphi\in\femspaceChemical{m}\right\}=\mathbb{R}^d.
    \end{equation*}
\end{asm}
The first assumption says that, for each $i\in\naturalindex{I_S}$, amongst all the non-boundary vertex normals on 
$\Gamma^m_i$, there should be at least one nonzero vector.
The second condition in Assumption~\ref{asm:1}, on the other hand,
is a very mild constraint on the interaction between bulk and interface meshes.
See \cite[Ass.~1]{EtoGarckeNurnberg2024} and 
\cite[Ass.~108]{BarrettGarckeRobertBook2020}
for very similar assumptions.

\begin{thm}[Existence and uniqueness]\label{thm:eu}
Let $\vVertex{m}{}\in\femspaceGammaVector{m}$ and $S^m$ be such that
Assumption~\ref{asm:1} holds.
Then, there exists a unique solution
$(\largeChemical^{m+1},\vVertex{m+1}{},\discreteCurvature{m+1}{\tension{}})\in\femspaceChemical{m}\times\femspaceGammaVector{m}\times\femspaceGamma{m}$
to \eqref{eq:fma}. 
\end{thm}
\begin{proof}
    Since the system \eqref{eq:fma} is linear and finite dimensional, with as
many unknowns as equations,
the existence of a solution is directly deduced from the uniqueness of the solution.
    To show the uniqueness, we need to prove that only the trivial solution
satisfies the homogeneous system.
    
    Suppose that
    $(\largeChemical,\vVertex{}{},\discreteCurvature{}{\tension{}})\in\femspaceChemical{m}\times\femspaceGammaVector{m}\times\femspaceGamma{m}$  satisfies
    \begin{subequations}\label{eq:hfma}
        \begin{equation}\label{eq:hfma_ml}
            \tau_m\inprodMassLamped{\nabla\largeChemical}{\nabla\varphi}{}{\Omega} + \sum_{i=1}^{I_S} \jump{\cpindexm{i}}{\cpindexp{i}}{\beta}\inprodMassLamped{\interpolation{h}{i}\left[\vVertex{}{}\cdot\femNormalVertex{m}{i}\right]}{\varphi}{(h)}{\polygonCurve{m}{i}} =0\quad\forall\varphi\in\femspaceChemical{m},
        \end{equation}
        \begin{equation}\label{eq:hfma_gtl}
            \inprodMassLamped{\discreteCurvature{}{\tension{}}}{\xi}{h}{\polygonCurve{m}{}} + \sum_{i=1}^{I_S}\jump{\cpindexm{i}}{\cpindexp{i}}{\beta}\inprodMassLamped{\largeChemical}{\xi}{(h)}{\polygonCurve{m}{i}} = 0\quad\forall\xi\in\femspaceGamma{m},
        \end{equation}
        \begin{equation}\label{eq:hfma_cv}
            \inprodMassLamped{\discreteCurvature{}{\tension{}}\femNormalVertex{m}{}}{\Vec{\eta}}{h}{\polygonCurve{m}{}} + \inprodMassLamped{\tension{}\sgrad\vVertex{}{}}{\sgrad\Vec{\eta}}{}{\polygonCurve{m}{}} = 0\quad\forall\Vec{\eta}\in\femspaceGammaVector{m}.
        \end{equation}
    \end{subequations}
    We choose $\varphi = \largeChemical$ in \eqref{eq:hfma_ml}, $\xi=\interpolation{h}{}\left[\vVertex{}{}\cdot\femNormalVertex{m}{}\right]$ in \eqref{eq:hfma_gtl} and $\Vec{\eta}=\vVertex{}{}$ in \eqref{eq:hfma_cv} to obtain
    \begin{align}
        \tau_m\|\nabla \largeChemical\|^2_2 + \sum_{i=1}^{I_S}\jump{\cpindexm{i}}{\cpindexp{i}}{\beta}\inprodMassLamped{\interpolation{h}{i}\left[\vVertex{}{}\cdot\femNormalVertex{m}{i}\right]}{\largeChemical}{(h)}{\polygonCurve{m}{i}} &= 0,\label{eq:eu_1}\\
        \inprodMassLamped{\discreteCurvature{}{\tension{}}}{\interpolation{h}{}\left[\vVertex{}{}\cdot\femNormalVertex{m}{}\right]}{h}{\polygonCurve{m}{}} + \sum_{i=1}^{I_S}\jump{\cpindexm{i}}{\cpindexp{i}}{\beta}\inprodMassLamped{\largeChemical}{\interpolation{h}{i}\left[\vVertex{}{}\cdot\femNormalVertex{m}{i}\right]}{(h)}{\polygonCurve{m}{i}} &= 0,\label{eq:eu_2}\\
        \inprodMassLamped{\discreteCurvature{}{\tension{}}\femNormalVertex{m}{}}{\vVertex{}{}}{h}{\polygonCurve{m}{}} + \inprodMassLamped{\tension{}\sgrad\vVertex{}{}}{\sgrad\vVertex{}{}}{}{\polygonCurve{m}{}} = 0.\label{eq:eu_3}
    \end{align}
    Combining \eqref{eq:eu_1}, \eqref{eq:eu_2} and \eqref{eq:eu_3}, we have
    \begin{equation*}
        \tau_m\|\nabla\largeChemical\|^2_2 + \inprodMassLamped{\tension{}\sgrad\vVertex{}{}}{\sgrad\vVertex{}{}}{}{\polygonCurve{m}{}} = 0.
    \end{equation*}
    The left-hand side of the above equation is non-negative, and hence we deduce that $\largeChemical$ and $\vVertex{}{}$ are constant,
    say $\largeChemical\equiv W_c$ and $\vVertex{}{}\equiv\vVertex{}{c}$. We deduce from \eqref{eq:hfma_ml} that
    \begin{equation*}
        \vVertex{}{c}\cdot\sum_{i=1}^{I_S}\jump{\cpindexm{i}}{\cpindexp{i}}{\beta}\inprodMassLamped{\femNormalVertex{m}{i}}{\varphi}{(h)}{\polygonCurve{m}{i}} = 0\quad\forall\varphi\in\femspaceChemical{m}.
    \end{equation*}
    By the second hypothesis of \Assumption{asm:1}, we have $\vVertex{}{} \equiv \vVertex{}{c} = \zerovec$.
    Meanwhile, we deduce from \eqref{eq:hfma_gtl} that
    \begin{equation}\label{eq:eu_4}
        \discreteCurvature{}{\tension{},i}\equiv W_c\jump{\cpindexm{i}}{\cpindexp{i}}{\beta}\quad\forall i\in\naturalindex{I_S},
    \end{equation}
    and hence $\discreteCurvature{}{\tension{}}$ is constant on each surface $\polygonCurve{m}{i}$. On recalling \eqref{eq:defV0}, we now choose in
\eqref{eq:hfma_cv} the test function $\Vec\eta \in\femspaceGammaVector{m}$ 
defined by $\vec\eta_i = \discreteCurvature{}{\tension{},i}\vec{z}^m_i$, where
$z^m_i := \sum_{k=1}^{\numSurfVertex{0}{i}}\femspaceBasisCurve{m,0}{i,k}\femNormalVertex{m}{i}(\vec q^{m}_{i,k})$ for $i\in\naturalindex{I_S}$.
On recalling that $\vec X = \vec 0$, we hence obtain
\begin{equation*}
    0 = \sum_{i=1}^{I_S} (\discreteCurvature{}{\tension{},i})^2 \inprodMassLamped{\femNormalVertex{m}{i}}{\vec{z}^m_i}{h}{\polygonCurve{m}{i}} 
= \sum_{i=1}^{I_S} (\discreteCurvature{}{\tension{},i})^2 \inprodMassLamped{\vec z^m_i}{\vec{z}^m_i}{h}{\polygonCurve{m}{i}}.
\end{equation*}
Therefore, the first condition in \Assumption{asm:1} yields that
$\discreteCurvature{}{\tension{},i}= 0$ for every $i\in\naturalindex{I_S}$, and hence $\discreteCurvature{}{\tension{}}\equiv 0$.
    In addition, \eqref{eq:eu_4} implies that $\largeChemical\equiv W_c = 0$ since $\jump{\cpindexm{i}}{\cpindexp{i}}{\beta}\neq 0$ for all $i\in\naturalindex{I_S}$.
    This completes the proof.
\end{proof}

\begin{thm}[Unconditional stability]\label{thm:us}
Let $(\largeChemical^{m+1},\vVertex{m+1}{},\discreteCurvature{m+1}{\tension{}})\in\femspaceChemical{m}\times\femspaceGammaVector{m}\times\femspaceGamma{m}$ be a solution to \eqref{eq:fma}. Then
    \begin{equation*}
        E(\polygonCurve{m+1}{}) + \tau_m\|\nabla\largeChemical^{m+1}\|^2_2\leq E(\polygonCurve{m}{}),
    \end{equation*}
 where we recall from \eqref{eq:surf_eg} that
        $E(\polygonCurve{m}{}) = \sum_{i=1}^{I_S}\tension{i}|\polygonCurve{m}{i}|$.
\end{thm}
\begin{proof}
    Substituting $\varphi = \largeChemical^{m+1}$ in \eqref{eq:fma_ml} and $\xi := \interpolation{h}{}\left[\frac{\vVertex{m+1}{} - \Vec\id}{\tau_m}\cdot\femNormalVertex{m}{}\right]$ in \eqref{eq:fma_gtl},
    we have
    \begin{align}\label{eq:us_1}
        \|\nabla\largeChemical^{m+1}\|^2_2 &= \frac{1}{\tau_m}\inprodMassLamped{\discreteCurvature{m+1}{\tension{}}}{\interpolation{h}{}\left[\left(\vVertex{m+1}{}-\Vec\id\right)\cdot\femNormalVertex{m}{}\right]}{h}{\polygonCurve{m}{}}\nonumber\\
         &= \frac{1}{\tau_m}\inprodMassLamped{\discreteCurvature{m+1}{\tension{}}\femNormalVertex{m}{}}{\vVertex{m+1}{} - \Vec\id}{h}{\polygonCurve{m}{}}.
    \end{align}
    Meanwhile, we choose $\Vec{\eta} = \vVertex{m+1}{} - \Vec\id\in\femspaceGammaVector{m}$ in \eqref{eq:fma_cv} to obtain
    \begin{equation}\label{eq:us_2}
        \inprodMassLamped{\discreteCurvature{m+1}{\tension{}}\femNormalVertex{m}{}}{\vVertex{m+1}{} - \Vec\id}{h}{\polygonCurve{m}{}} = -\inprodMassLamped{\tension{}\sgrad\vVertex{m+1}{}}{\sgrad(\vVertex{m+1}{} - \Vec\id)}{}{\polygonCurve{m}{}}.
    \end{equation}
    We now recall \cite[Lemma 57]{BarrettGarckeRobertBook2020}, which states that
    \begin{equation}\label{eq:us_3}
        \inprodMassLamped{\sgrad\vVertex{m+1}{i}}{\sgrad(\vVertex{m+1}{i} - \Vec\id)}{}{\polygonCurve{m}{i}}\geq |\polygonCurve{m+1}{i}| - |\polygonCurve{m}{i}| ,
    \end{equation}
for $i\in\naturalindex{I_S}$. Multiplying \eqref{eq:us_3} with $\tension{i}$, summing
and combining with \eqref{eq:us_1} and \eqref{eq:us_2} yields
    \[
        \tau_m\|\nabla\largeChemical^{m+1}\|^2_2 = -\inprodMassLamped{\tension{}\sgrad\vVertex{m+1}{}}{\sgrad(\vVertex{m+1}{} - \Vec\id)}{}{\polygonCurve{m}{}}
%        &\leq -E(\polygonCurve{m+1}{}) + E(\polygonCurve{m}{}) \\ 
        \leq -E(\polygonCurve{m+1}{}) + E(\polygonCurve{m}{}).
    \]
\end{proof}

\section{Matrix form of the linear system}\label{sec:mf}
In this section, we derive a matrix form of the linear system induced by the finite element approximation \eqref{eq:fma} at each time step.
To facilitate understanding, we first concentrate on the two-dimensional
three-phase example shown in Figure~\ref{fig:3p}. Namely, we consider
a curve network with $I_R = 3$, $I_S = 3$, $I_T = 2$, $(\cpindexp{1},\cpindexm{1}) = (3,2)$, $(\cpindexp{2},\cpindexm{2})=(1,3)$ and $(\cpindexp{3},\cpindexm{3})=(2,1)$.

Given $\Gamma^{m}$, let $(\largeChemical^{m+1}, \discreteCurvature{m+1}{\tension{}},\vVertex{m}{} + \delta\vVertex{m+1}{})\in \femspaceChemical{m}\times\femspaceGamma{m}\times\femspaceGammaVector{m}$
be the unique solution of \eqref{eq:fma}.
Let $\numDomVertex{m}$ and $\numSurfVertex{m}{c}\,(c\in\naturalindex{I_S})$ be the number of all vertices on $\triangulation{m}$ and $\polygonCurve{m}{c}$, respectively.
We let $\numSurfVertexFull{m} := \sum_{c=1}^{I_S} \numSurfVertex{m}{c}$.
In the sequel, we will regard $(\largeChemical^{m+1},\discreteCurvature{m+1}{\tension{}},\delta\vVertex{m+1}{})$ as their
coefficients with respect to the basis functions $\{\femspaceBasisBulk{m}{i}\}_{1\leq i\leq \numDomVertex{m}}$ and $\{\{\femspaceBasisCurve{m}{c,\ell}\}_{1\leq \ell\leq \numSurfVertex{m}{c}}\}_{c=1}^{I_S}$
of $\femspaceChemical{m}$ and $\femspaceGamma{m}$, respectively.
We also define the orthogonal projection
$\metrics{P}:(\bR^d)^K \to \underline{\mathbb X}$
onto the Euclidean space associated with $\femspaceGammaVector{m}$,
see \cite[p.~211]{clust3d}.
%We now introduce the orthogonal projection $\metrics{\mathcal{P}}:\femspaceGamma{m}\to\femspaceGammaVector{m}$ and its discrete variant
%\begin{equation*}
%    \metrics{P}:(\mathbb{R}^2)^N\to \left\{(\vec{z_1},\vec{z_2},\vec{z_3})\in(\mathbb{R}^2)^N\mid [\vec{z_1}]_0=[\vec{z_2}]_0=[\vec{z_3}]_0,\ [\vec{z_1}]_{N_1} = [\vec{z_2}]_{N_2} = [\vec{z_3}]_{N_3} \right\}.
%\end{equation*}
Then, the solution of \eqref{eq:fma} can be written as $(\largeChemical^{m+1},\discreteCurvature{m+1}{\tension{}},\vVertex{m}{} + \metrics{P}\delta\vVertex{m+1}{})$,
where $(\largeChemical^{m+1},\discreteCurvature{m+1}{\tension{}},\delta\vVertex{m+1}{})$ solves the following linear system:

\begin{equation}\label{eq:mf}
    \begin{pmatrix}
        \tau_m A_\Omega & O & \Vec{N}_{\Omega,\Gamma}^T\metrics{P}\\
        B_{\Omega,\Gamma} & C_\Gamma & O\\
        O & \metrics{P}\Vec{D}_\Gamma & \metrics{P}\metrics{E_\Gamma}\metrics{P}
    \end{pmatrix}\begin{pmatrix}
        \largeChemical^{m+1}\\
        \discreteCurvature{m+1}{\tension{}}\\
        \delta\vVertex{m+1}{}
    \end{pmatrix} = \begin{pmatrix}
        O \\
        O \\
        \bv -\metrics{P}\metrics{E_\Gamma}\metrics{P}\vVertex{m}{}
    \end{pmatrix},
\end{equation}
where $A_\Omega\in\mathbb{R}^{\numDomVertex{m}\times\numDomVertex{m}}$, $\Vec{N}_{\Omega,\Gamma}\in(\mathbb{R}^d)^{\numSurfVertexFull{m}\times\numDomVertex{m}}$, $B_{\Omega,\Gamma}\in\mathbb{R}^{\numSurfVertexFull{m}\times\numDomVertex{m}}$,
$C_\Gamma\in\mathbb{R}^{\numSurfVertexFull{m}\times \numSurfVertexFull{m}}$, $\vec{D}_\Gamma\in(\mathbb{R}^d)^{\numSurfVertexFull{m}\times \numSurfVertexFull{m}}$ and $E_\Gamma\in(\mathbb{R}^{d\times d})^{\numSurfVertexFull{m}\times \numSurfVertexFull{m}}$ are defined by
\begin{align*}
    \vec{N}_{\Omega,\Gamma} := \begin{pmatrix}
        \jump{\cpindexm{1}}{\cpindexp{1}}{\beta}\Vec{N}_1 \\
        \jump{\cpindexm{2}}{\cpindexp{2}}{\beta}\Vec{N}_2 \\
        \jump{\cpindexm{3}}{\cpindexp{3}}{\beta}\Vec{N}_3 \\
    \end{pmatrix},\quad
    &
    B_{\Omega,\Gamma} := \begin{pmatrix}
        \jump{\cpindexm{1}}{\cpindexp{1}}{\beta}B_1 \\
        \jump{\cpindexm{2}}{\cpindexp{2}}{\beta}B_2 \\
        \jump{\cpindexm{3}}{\cpindexp{3}}{\beta}B_3 \\
    \end{pmatrix},\quad
    C_\Gamma := \begin{pmatrix}
        C_1 & O & O \\
        O & C_2 & O \\
        O & O & C_3
        \end{pmatrix},
\end{align*}
\begin{align*}
    \Vec{D}_\Gamma := \begin{pmatrix}
        \Vec{D}_1 & O & O \\
        O & \Vec{D}_2 & O \\
        O & O & \Vec{D}_3
        \end{pmatrix},\quad
    &
    \metrics{E_\Gamma} := \begin{pmatrix}
        \tension{1}\metrics{E_1} & O & O \\
        O & \tension{2}\metrics{E_2} & O \\
        O & O & \tension{3}\metrics{E_3}
    \end{pmatrix}
\end{align*}
with
\begin{align*}
    \left[A_\Omega\right]_{i,j} := \inprodMassLamped{\nabla\femspaceBasisBulk{m}{j}}{\nabla\femspaceBasisBulk{m}{i}}{}{\Omega},\qquad
    &\left[\vec{N}_{c}\right]_{\ell,i} := \inprodMassLamped{\femspaceBasisCurve{m}{c,\ell}}{\femspaceBasisBulk{m}{i}}{(h)}{\polygonCurve{m}{c}}\femNormalVertex{m}{c,\ell},\\
    \left[B_c\right]_{k,j} := \inprodMassLamped{\femspaceBasisBulk{m}{j}}{\femspaceBasisCurve{m}{c,k}}{(h)}{\polygonCurve{m}{c}},\qquad
    &\left[C_c\right]_{k,\ell} := \inprodMassLamped{\femspaceBasisCurve{m}{c,\ell}}{\femspaceBasisCurve{m}{c,k}}{(h)}{\polygonCurve{m}{c}},\\
    \left[\Vec{D}_c\right]_{k,\ell} := \inprodMassLamped{\femspaceBasisCurve{m}{c,\ell}}{\femspaceBasisCurve{m}{c,k}}{h}{\polygonCurve{m}{c}}\femNormalVertex{m}{c,\ell},\qquad
    & \left[\metrics{E_c}\right]_{k,\ell} := \inprodMassLamped{\sgrad\femspaceBasisCurve{m}{c,\ell}}{\sgrad\femspaceBasisCurve{m}{c,k}}{}{\polygonCurve{m}{c}}
\end{align*}
for each $c\in \naturalindex{I_S}$.

Generalizing \eqref{eq:mf} to general surface clusters in $\bR^d$, $d=2,3$, is
now straightforward.
The only required changes are to adapt the definition of $\metrics{P}$ 
to the cluster at hand, recall \cite{clust3d}, and to let the block matrices
be given by
$(\Vec{N}_{\Omega,\Gamma})_{i,1} := \jump{\cpindexm{i}}{\cpindexp{i}}{\beta}\vec{N}_i\,(i\in\naturalindex{I_S})$,
$(B_{\Omega,\Gamma})_{i,1} := \jump{\cpindexm{i}}{\cpindexp{i}}{\beta}B_i\,(i\in\naturalindex{I_S})$,
$C_\Gamma := \operatorname{diag}(C_i)_{i=1,\cdots,I_S}$, $\Vec{D}_\Gamma := \operatorname{diag}(\vec{D}_i)_{i=1,\cdots,I_S}$
and $\metrics{E_\Gamma} := \operatorname{diag}(\tension{i}\metrics{E_i})_{i=1,\cdots,I_S}$.

The linear system \eqref{eq:mf} 
can either be solved with a Krylov subspace type iterative
solver, together with a suitable preconditioner, or it can be reduced with the
help of a Schur complement approach. In fact, on following the approach
from \cite{BarrettGarckeRobert2010}, we can eliminate first 
$\discreteCurvature{m+1}{\tension{}}$ and then $\largeChemical^{m+1}$,
to obtain a linear system just in terms of $\delta\vVertex{m+1}{}$, which can
then be solved with preconditioned iterative solvers as in \cite{clust3d}.
Similarly to \cite{BarrettGarckeRobert2010}, when eliminating 
$\largeChemical^{m+1}$ some attention must be placed on the fact that 
$A_\Omega$ has a one-dimensional kernel. We refer to 
\cite{BarrettGarckeRobert2010} for the precise details.

\section{Structure-preserving numerical scheme}\label{sec:spns}
In this section, we refine the previously proposed 
fully discrete scheme \eqref{eq:fma} such that the discrete total energy
content is preserved at each time step.
To this end, we follow the idea by Bao, Garcke, N\"{u}rnberg and Zhao \cite{BaoGarckeNuernbergZhao2022}
which can interpolate between discrete surface clusters at intermediate timestamps.
We define $\Gamma_i^h(t)=\bigcup_{j=1}^{J_i} \overline{\sigma_{i,j}^{h}(t)}$ as the interpolating polyhedral surfaces
where $\{\sigma_{i,j}^{h}(t)\}_{j=1}^{J_i}$ are mutually disjoint $(d-1)$-simplices with vertices $\{\vec q_{i,k}^{h}(t)\}_{k=1}^{K_i}$ given by
\begin{equation}
\label{eq:qh}
\vec q_{i,k}^{h}(t) = \frac{t_{m+1} - t}{\tau_m}\vec q_{i,k}^{m} + \frac{t - t_m}{\tau_m}\,\vec q_{i,k}^{m+1},\quad t\in[t_m,~t_{m+1}],\quad k = 1,\ldots, K_i.
\end{equation}
We now define a discrete normal vector $\normal{m+\frac{1}{2}}{}$ by
\begin{equation*}
    \normal{m+\frac{1}{2}}{i}\mid_{\simplex{m}{i,j}} := \frac{1}{\tau_m|\Normal{\simplex{m}{i,j}}|}\int_{t_m}^{t_{m+1}}\Normal{\simplex{h}{i,j}(t)}\,{\rm d}t\qquad\mbox{for}\quad j\in\naturalindex{\numSimplex{}{i}}\quad\mbox{and}\quad i\in\naturalindex{I_S},
\end{equation*}
and the associated vertex normal vector $\femNormalVertex{m+\frac{1}{2}}{}$ 
through
\begin{equation*}
    \inprodMassLamped{\femNormalVertex{m+\frac{1}{2}}{i}}{\Vec{\xi}}{h}{\polygonCurve{m}{i}} = \inprodMassLamped{\normal{m+\frac{1}{2}}{i}}{\Vec{\xi}}{}{\polygonCurve{m}{i}}\qquad\forall\ \Vec{\xi}\in V^h(\polygonCurve{m}{i}),
\qquad i \in \naturalindex{I_S}.
\end{equation*}
Using the vertex normal vector $\femNormalVertex{m+\frac{1}{2}}{}$, 
we introduce another finite element approximation of \eqref{eq:wkf}.
For each $m \geq 0$, given $\Gamma^{m}$,
we find $(\largeChemical^{m+1},\vVertex{m+1}{},\discreteCurvature{m+1}{\tension{}})\in\femspaceChemical{m}\times\femspaceGammaVector{m}\times\femspaceGamma{m}$ such that\newline\newline

\begin{subequations}\label{eq:fma2}
\textbf{Motion law: } For all $\varphi\in\femspaceChemical{m}$, it holds that
\begin{equation}\label{eq:fma_ml2}
    \inprodMassLamped{\nabla\largeChemical^{m+1}}{\nabla\varphi}{}{\Omega} + \sum_{i=1}^{I_S}\jump{\cpindexm{i}}{\cpindexp{i}}{\beta} \inprodMassLamped{\interpolation{h}{i}\left[\frac{\vVertex{m+1}{i} - \Vec\id}{\tau_m}\cdot\femNormalVertex{m+\frac{1}{2}}{i}\right]}{\varphi}{(h)}{\polygonCurve{m}{i}} = 0.
\end{equation}
\textbf{Gibbs--Thomson law: } For all $\xi\in\femspaceGamma{m}$, it holds that
\begin{equation}\label{eq:fma_gtl2}
    \inprodMassLamped{\discreteCurvature{m+1}{\tension{}}}{\xi}{h}{\polygonCurve{m}{}} + \sum_{i=1}^{I_S}\jump{\cpindexm{i}}{\cpindexp{i}}{\beta}\inprodMassLamped{\largeChemical^{m+1}}{\xi}{(h)}{\polygonCurve{m}{i}} = 0.
\end{equation}
\textbf{Curvature vector: } For all $\Vec{\eta}\in\femspaceGammaVector{m}$, it holds that
\begin{equation}\label{eq:fma_cv2}
    \inprodMassLamped{\discreteCurvature{m+1}{\tension{}}\femNormalVertex{m+\frac{1}{2}}{}}{\Vec{\eta}}{h}{\Gamma^m} + \inprodMassLamped{\tension{}\sgrad\vVertex{m+1}{}}{\sgrad\Vec{\eta}}{}{\polygonCurve{m}{}} = 0.
\end{equation}
\end{subequations}
We observe that the only change from \eqref{eq:fma} to \eqref{eq:fma2} is that
$\vec\omega^m$ is now replaced by $\vec\omega^{m+\frac12}$. Since the latter
depends on $\Gamma^{m+1} = \vec X^{m+1}(\Gamma^m)$, the scheme \eqref{eq:fma2} 
leads to a nonlinear system of equations, in contrast to \eqref{eq:fma}. 
We shall show that, \revised{in addition to being unconditionally stable}, 
\eqref{eq:fma2} \revised{also} 
preserves the total energy content at each time step.
Before proving this, we recall the following lemma from \cite[Lemma 3.1]{BaoGarckeNuernbergZhao2022} without its proof (see also \cite[Lemma 4.4]{GarckeNuernbergZhao2024}).

\begin{lem}\label{lem:vpsub}
    Let $\vVertex{m+1}{}\in\femspaceGammaVector{m}$.
    Then, for every $\ell\in\naturalindex{I_R}$, it holds that
    \begin{equation*}
        \volume{\subdomainsDiscrete{\ell}{m+1}} - \volume{\subdomainsDiscrete{\ell}{m}} = -\sum_{i=1}^{I_S}\inprodMassLamped{\left(\vVertex{m+1}{} - \Vec\id\right)\cdot\normal{m+\frac{1}{2}}{}}{\chi_{\ell,i}}{}{\polygonCurve{m}{i}},
    \end{equation*}
    where
    \begin{equation*}
        \chi_{\ell,i} :=
        \begin{cases}
            1&\qquad\mbox{if}\quad \cpindexp{i} = \ell,\\
            -1&\qquad\mbox{if}\quad \cpindexm{i} = \ell,\\
            0&\qquad\mbox{otherwise}.
        \end{cases}
    \end{equation*}
\end{lem}
With the help of \Lemma{lem:vpsub}, we can establish the following theorem,
where in preparation we define the discrete total energy content as
    \begin{equation} \label{eq:vm}
        v^m := \sum_{\ell=1}^{I_R}\beta_\ell\volume{\subdomainsDiscrete{\ell}{m}}\qquad\mbox{for}\quad m=0,\ldots,M.
    \end{equation}

\begin{thm}\label{thm:vp}
    Let $(\largeChemical^{m+1},\vVertex{m+1}{},\discreteCurvature{m+1}{\tension{}})\in\femspaceChemical{m}\times\femspaceGammaVector{m}\times\femspaceGamma{m}$ be a solution to \eqref{eq:fma2}. 
\revised{%
Then
    \begin{equation*}
        E(\polygonCurve{m+1}{}) + \tau_m\|\nabla\largeChemical^{m+1}\|^2_2\leq E(\polygonCurve{m}{}).
    \end{equation*}
Moreover,} it holds that $v^{m+1} = v^m$.
\end{thm}
\begin{proof}
\revised{%
The proof of the stability bound is identical to the proof of
Theorem~\ref{thm:us}. In order to prove the energy content preservation, we}
deduce from \Lemma{lem:vpsub} that
    \begin{align*}
        &v^{m+1} - v^m= \sum_{\ell=1}^{I_R}\beta_\ell\left\{-\sum_{\substack{i\in\naturalindex{I_S}\\\cpindexp{i} = \ell}}\inprodMassLamped{\vVertex{m+1}{} - \Vec\id}{\normal{m+\frac{1}{2}}{i}}{}{\polygonCurve{m}{i}} + \sum_{\substack{i\in\naturalindex{I_S}\\\cpindexm{i}=\ell}}\inprodMassLamped{\vVertex{m+1}{} - \Vec\id}{\normal{m+\frac{1}{2}}{i}}{}{\polygonCurve{m}{i}}\right\}\\
        &= \sum_{\ell=1}^{I_R}\beta_\ell\left\{-\sum_{\substack{i\in\naturalindex{I_S}\\\cpindexp{i} = \ell}}\inprodMassLamped{\vVertex{m+1}{} - \Vec\id}{\femNormalVertex{m+\frac{1}{2}}{i}}{h}{\polygonCurve{m}{i}} + \sum_{\substack{i\in\naturalindex{I_S}\\\cpindexm{i}=\ell}}\inprodMassLamped{\vVertex{m+1}{} - \Vec\id}{\femNormalVertex{m+\frac{1}{2}}{i}}{h}{\polygonCurve{m}{i}}\right\}\\
        &= -\tau_m\sum_{i=1}^{I_S}\jump{\cpindexm{i}}{\cpindexp{i}}{\beta}\inprodMassLamped{\interpolation{h}{}\left[\frac{\vVertex{m+1}{}-\Vec\id}{\tau_m}\cdot\femNormalVertex{m+\frac{1}{2}}{}\right]}{1}{(h)}{\polygonCurve{m}{i}} 
%= \tau_m\inprodMassLamped{\nabla\largeChemical^{m+1}}{\nabla 1}{}{\Omega} 
= 0.
    \end{align*}
    Here, we have invoked \eqref{eq:fma_ml2} with $\varphi=1$ to derive the last equality.
\end{proof}

\begin{rem}
We observe that as $\femNormalVertex{m+\frac{1}{2}}{}$
depends on $\vVertex{m+1}{}$, the scheme \eqref{eq:fma2} 
is no longer linear. In practice the nonlinear systems of equations arising 
at each time level of \eqref{eq:fma2} can be solved with the aid of a 
simple lagged iteration as mentioned in \cite[\S3]{Nurnberg202203}.
\end{rem}

\section{Numerical experiments}\label{sec:ne}
In this section, we present numerical experiments to verify the feasibility of the proposed scheme together with its accuracy and efficiency.
We implemented the fully discrete finite element approximations 
\eqref{eq:fma} and \eqref{eq:fma2} within the
finite element toolbox ALBERTA, see \cite{Alberta}.
The arising linear systems of the form \eqref{eq:mf} are either solved with 
a GMRes iterative solver, applying as preconditioner a least squares solution 
of the block matrix in \eqref{eq:mf} without the projection matrices
$\metrics{P}$. For the computation of the least squares solution we employ
the sparse factorization package SPQR, see \cite{Davis11}.
Alternatively, and especially in 3d, we solve \eqref{eq:mf} with the help of
the previously discussed Schur complement approach. Here we use either
UMFPACK, see \cite{Davis04}, or a multigrid algorithm for the solution of the
arising linear subproblems. We refer to \cite{BarrettGarckeRobert2010} for more
details in the two-phase case.

\begin{figure}[H]
    \centering
    \includegraphics[angle=-0,keepaspectratio,scale=0.4]{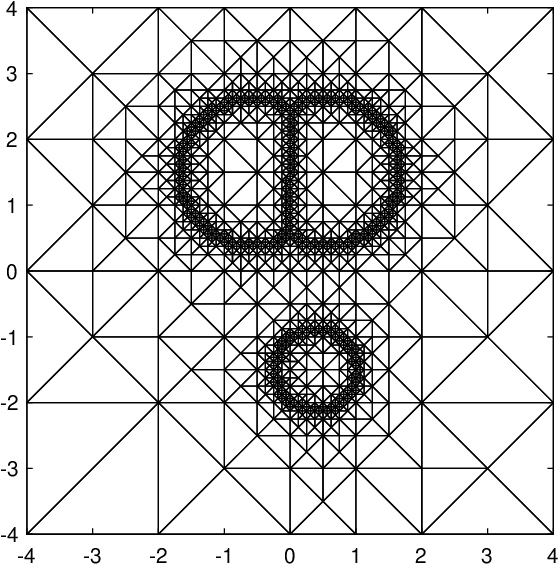}
    \caption{Adaptive bulk mesh in two space dimensions.}\label{fig:mesh}
\end{figure}

For the triangulation $\triangulation{m}$ of the bulk domain $\Omega$, that is used
for the bulk finite element space $\femspaceChemical{m}$, we use an adaptive mesh with
fine elements close to the interface $\Gamma^m$ and coarser elements away from it.
The precise strategy is as described in \cite[\S5.1]{BarrettGarckeRobert2010}
and for a domain $\Omega=(-H,H)^d$, for some $H>0$, and two integer parameters $N_c < N_f$ 
results in elements with maximal diameter approximately equal to $h_f = \frac{2H}{N_f}$ close to $\Gamma^m$,
and elements with maximal diameter approximately equal to $h_c = \frac{2H}{N_c}$ far away from it.
Throughout this section, we always set $H = 4$.
See Figure~\ref{fig:mesh} for an example adaptive bulk mesh in two space
dimensions.
\revised{Observe that employing adaptive bulk meshes does not invalidate the
stability results in Theorems~\ref{thm:us} and \ref{thm:vp}.}

We stress that due to the unfitted nature of our finite element approximations,
special quadrature rules need to be employed in order to assemble terms that
feature both bulk and surface finite element functions. For all the
computations presented in this section, \revised{unless otherwise stated}, 
we use true integration for these terms, and we refer to
\cite{BarrettGarckeRobert2010,Nurnberg202203} for details on the practical
implementation.
Throughout this section we use (almost) uniform time steps, in that
$\tau_m=\tau$ for $m=0,\ldots, M-2$ and 
$\tau_{M-1} = T - t_{m-1} \leq \tau$. 
In the sequel, we set all the physical parameters to unity, and in particular
$\tension{i} = 1$ for all $i\in\naturalindex{I_S}$, unless otherwise stated.

\subsection{Convergence experiment in 2d}
\begin{figure}[H]
    \centering
    \includegraphics[keepaspectratio, scale=0.5]{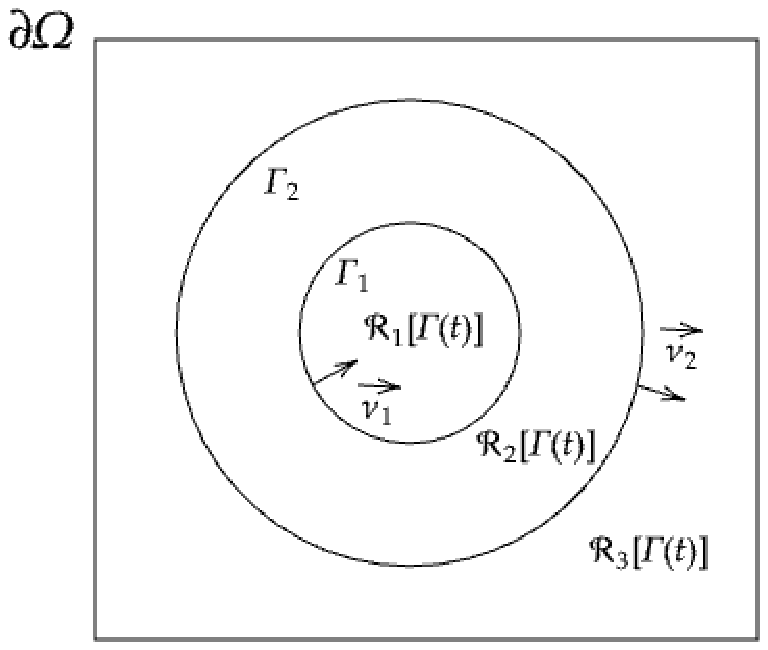}
    \caption{Two concentric circles.}\label{fig:2c}
\end{figure}

% for old annular repeat see
% /home/rn/hpc_cluster/data/alberta/egn2/2d.convergence_experiment
We check the accuracy of the proposed scheme using a rigorous solution to the degenerate multi-phase Stefan problem \eqref{eq:sharp_p} and compare it with the numerical solution.
To this end, we consider a curve network which consists of two concentric circles in $\mathbb{R}^2$.
See Figure~\ref{fig:2c} for the setting.

Then, we can obtain the formula of the rigorous solution to the above system as follows:
\begin{subequations} \label{eq:truesol2d}
\begin{equation} \label{eq:w2d}
    w(x,t) = \begin{cases}
        -\frac{1}{\jump{2}{1}{\beta}R_1(t)}&\qquad\mbox{if}\quad |x| < R_1(t),\\
        -\frac{1}{\jump{2}{1}{\beta}R_1(t)} + \alpha_2(R_1(t),R_2(t))\log{\frac{|x|}{R_1(t)}},&\qquad\mbox{if}\quad R_1(t)\leq |x| < R_2(t),\\
        \frac{1}{\jump{2}{3}{\beta}R_2(t)}&\qquad\mbox{if}\quad R_2(t)\leq |x|,
    \end{cases}
\end{equation}
where
\begin{equation*}
    \alpha_2(R_1,R_2) := \frac{\frac{1}{\jump{2}{1}{\beta}R_1} + \frac{1}{\jump{2}{3}{\beta}R_2}}{\log{\frac{R_2}{R_1}}}\qquad\mbox{for}\quad R_1, R_2 > 0.
\end{equation*}
We deduce from the motion law that the radii $R_1(t)$ and $R_2(t)$ satisfy the following ordinary differential equations:
% \begin{equation}\label{eq:rd2}
%     \begin{cases}
%         \dot{R}_1 = -\frac{\alpha_2(R_1,R_2)}{\jump{2}{1}{\beta}R_1},\\
%         \dot{R}_2 = -\frac{\alpha_2(R_1,R_2)}{\jump{2}{3}{\beta}R_2},
%     \end{cases}
% \end{equation}
\begin{align}\label{eq:rd2}
    \begin{split}
        \dot{R}_1(t) &= -\frac{\alpha_2(R_1(t),R_2(t))}{\jump{2}{1}{\beta}R_1(t)},\\
        \dot{R}_2(t) &=
-\frac{\alpha_2(R_1(t),R_2(t))}{\jump{2}{3}{\beta}R_2(t)},
    \end{split}
\end{align}
\end{subequations}
where $\dot{R}_i$ denotes the time derivative of $R_i$.
In order to make an accurate comparison between true and discrete solution,
we prefer to deal with a scalar ordinary differential equation for a single radius.
Thanks to \Proposition{prop:mass_p}, we see that the function
\begin{equation*}
    t\mapsto \pi\jump{2}{1}{\beta}{R_1(t)^2} - \pi\jump{2}{3}{\beta}{R_2(t)^2} + \beta_3\mL^{2}(\Omega),
\end{equation*}
is constant. This implies that $t\mapsto\jump{2}{1}{\beta}R_1(t)^2 - \jump{2}{3}{\beta}R_2(t)^2$ is also constant.
Thus, letting $D_2 := \jump{2}{1}{\beta}R_1(0)^2 - \jump{2}{3}{\beta}R_2(0)^2$, we obtain
\begin{equation}\label{eq:rd}
    \dot{R_2}(t) = -F(R_2(t))\qquad\mbox{with}\quad F(u) := \frac{\alpha_2\left(\sqrt{\frac{D_2 + \jump{2}{3}{\beta}u^2}{\jump{2}{1}{\beta}}},u\right)}{\jump{2}{3}{\beta} u},
\end{equation}
where $F(u)$ is defined for every $u$ such that $(D_2 + \jump{2}{3}{\beta} u^2) / \jump{2}{1}{\beta} > 0$.
We can solve \eqref{eq:rd} in terms of the root-finding method as in \cite[Eq.(8.2)]{EtoGarckeNurnberg2024}, namely by finding $R_2$ such that
\begin{equation*}
    0 = t + \int_{R_2(0)}^{R_2}\frac{1}{F(u)}\,du.
\end{equation*}

\begin{figure}[H]
% ../../plotradii
% cp radii_clean.ps ~/tex/harald/toku2/figures/simul2c_r.ps && scpp radii_clean.ps e23:tex/harald/toku2/figures/simul2c_r.ps
% cp radii_clean.ps ~/tex/harald/toku2/figures/simul2c_long_r.ps && scpp radii_clean.ps e23:tex/harald/toku2/figures/simul2c_long_r.ps
\centering
\includegraphics[angle=-90, keepaspectratio, scale=0.25]{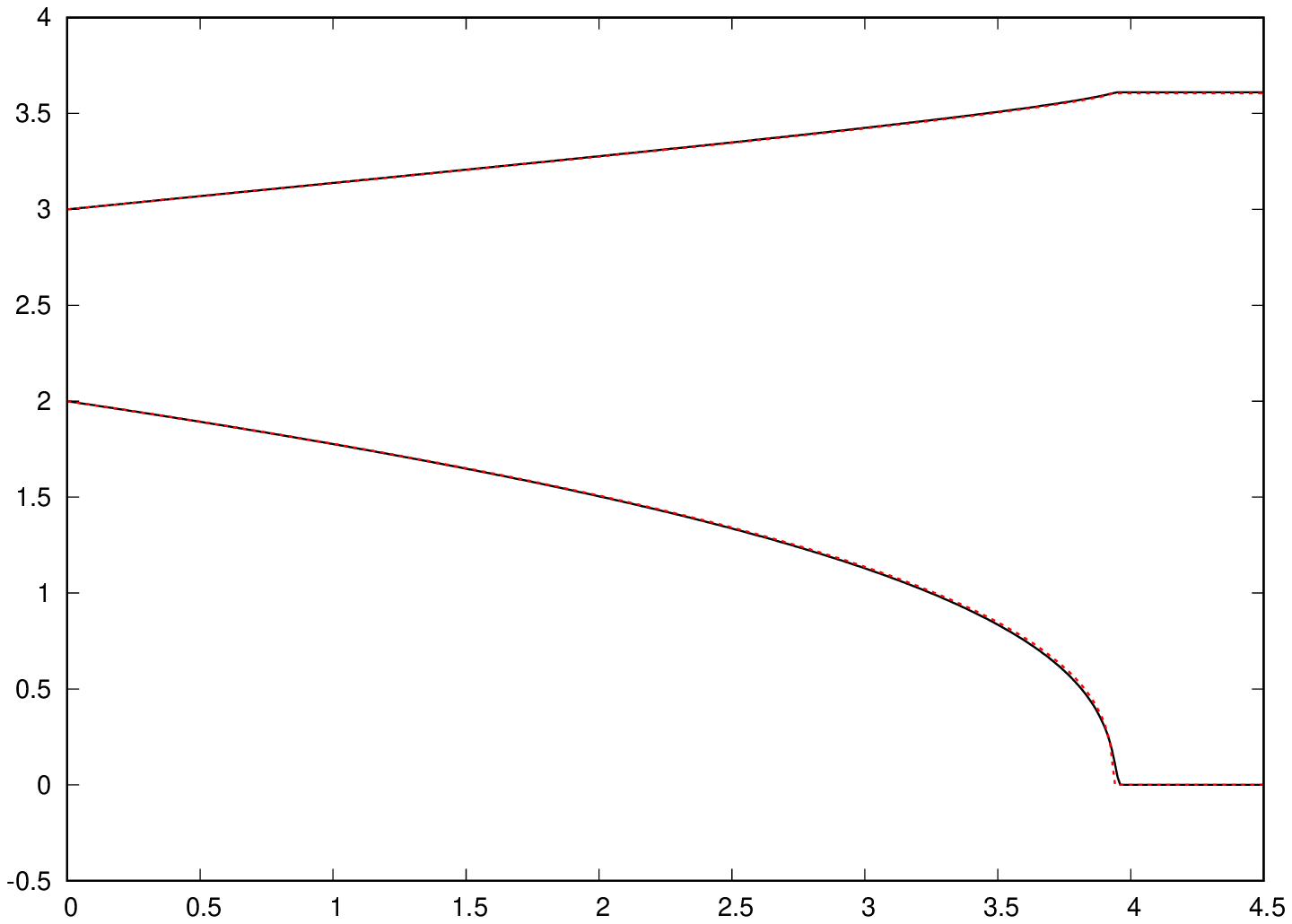}
\setlength{\abovecaptionskip}{20pt}
\caption{Comparison of the discrete (black, solid)\\ and exact (red, dashed) 
radii for $(\beta_1,\beta_2,\beta_3) = (-1,0,1)$.}
% ~/hpc_cluster/data/alberta/egn2/2d.convergence_experiment_NoC2/plot/run1
% ~/hpc_cluster/data/alberta/egn2/2d.convergence_experiment_NoC2/long/run1
\label{fig:simul2c}
\end{figure}

We show in Figure~\ref{fig:simul2c} the evolutions of the radii $R_1(t)$ and $R_2(t)$ satisfying \eqref{eq:rd2},
and compare them to the evolutions of the radii of the discrete solution for our scheme \eqref{eq:fma}.
We note the perfect agreement between the two.
For this choice of $\phaseContent$, the inner
phase shrinks and eventually disappears, leaving a trivial stationary solution
for the two-phase problem.

\begin{figure}[H]
% ../../plotradii
% cp radii_clean.ps ~/tex/harald/toku2/figures/simul2cb_r.ps && scpp radii_clean.ps e23:tex/harald/toku2/figures/simul2cb_r.ps
% cp radii_clean.ps ~/tex/harald/toku2/figures/simul2cb_long_r.ps && scpp radii_clean.ps e23:tex/harald/toku2/figures/simul2cb_long_r.ps
\center
\includegraphics[angle=-90, keepaspectratio, scale=0.25]{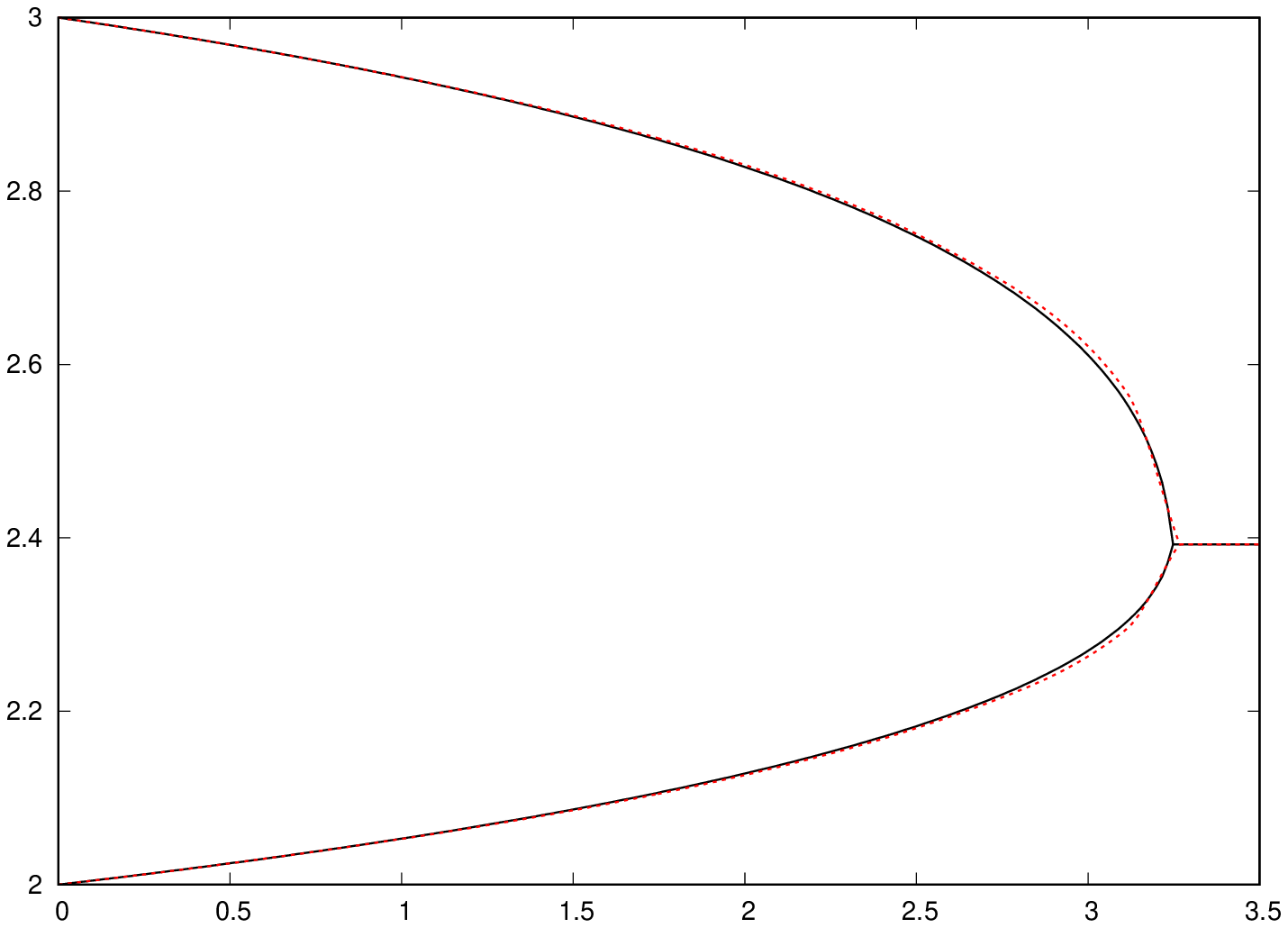}
\setlength{\abovecaptionskip}{20pt}
\caption{Comparison of the discrete (black, solid)\\
and exact (red, dashed) radii for $(\beta_1,\beta_2,\beta_3) = (-1.6,0.3,1.3)$.}
% ~/hpc_cluster/data/alberta/egn2/2d.convergence_experiment_NoC2b/lin/run0
% ~/hpc_cluster/data/alberta/egn2/2d.convergence_experiment_NoC2b/plot/run1
\label{fig:simul2cb}
\end{figure}

The evolution for a different choice of $\phaseContent$ is shown Figure~\ref{fig:simul2cb},
and we observe that the second phase vanishes as inner and outer circles approach each
other to finally only leave a single circular interface that separates phases $1$ and $3$.

In summary, we see from \Figure{fig:simul2c} and \Figure{fig:simul2cb}, that
setting $\beta_1$ to be very small leads to a growth of the inner circle and a 
shrinking of the outer circle.
In fact, it can be seen from \eqref{eq:rd2} that changing the value of 
$\beta_1$ may cause the radius $R_1(t)$ to evolve differently due to the change of the sign of $\dot{R}_1(t)$. It is also easy to check that
$\dot{R}_2(t)$ has a different sign to $\dot{R}_1(t)$, and so the outer circle
will move in the opposite direction to the inner one.

We perform a convergence experiment for the true solution of \eqref{eq:truesol2d} on the time interval $[0,1]$. 
For the discretization parameters, for $i=0\to 4$, we set 
$N_f = \frac12 K^0_\Gamma = 2^{7+i}$, $N_c = 4^{i}$ and $\tau= 4^{3-i}\times10^{-3}$. 
For the discrete solutions of \eqref{eq:fma} and \eqref{eq:fma2} we 
define the errors 
\begin{equation*}
    \errorXx =
    \max_{m=1,\ldots, M} 
    \max_{i=1,\ldots,I_S}
    \max_{k=1,\ldots,K^m_i}
    \operatorname{dist}(\vec{q}^m_{i,k}, \Gamma_i(t_m))
\end{equation*}
and
\begin{equation*}
    \errorUu = \max_{m=1,\ldots, M} \|\largeChemical^m - I^m \smallChemical(\cdot,t_m)\|_{L^\infty},
\end{equation*}
where $I^m : C(\overline\Omega) \to \femspaceChemical{m}$ denotes the standard 
interpolation operator.
We also let
\begin{equation*}
    h^m_\Gamma = \max_{i=1,\ldots,I_S}
    \max_{j = 1,\ldots,\numSimplex{}{i}} |\simplex{m}{i,j}|
\end{equation*}
and
\begin{equation*}
    v_\Delta^M = v^M - v^0,
\end{equation*}
recall \eqref{eq:vm}.

We show the errors for the nonlinear scheme \eqref{eq:fma2} 
and for the linear scheme \eqref{eq:fma} in Tables~\ref{tab:2dvol} and
\ref{tab:2d}, respectively.
Here, we recall that $\numDomVertex{m}$ denotes the degree of freedom of $\femspaceChemical{m}$, while $K^m_\Gamma$ counts the total number of vertices on $\Gamma^m$.
\begin{table}[H]
\center
\begin{tabular}{c|c|c|c|c|c|c}
% grep "LaTeX:" run*/specs.used -A 2 | grep '\\\\' | sed 's/.*specs.used- //'
 $h_{f}$ & $h^M_\Gamma$ & $\errorUu$ & $\errorXx$ & $K^M_\Omega$ & $K^M_\Gamma$
 & $|v_\Delta^M|$ \\ \hline
6.2500e-02 & 1.5401e-01 & 4.2628e-03 & 2.0103e-03 & 2721  & 256 & $<10^{-10}$\\
3.1250e-02 & 7.7011e-02 & 1.4983e-03 & 1.9699e-03 & 5417  & 512 & $<10^{-10}$\\
1.5625e-02 & 3.8502e-02 & 1.4098e-03 & 1.3554e-03 & 10837 & 1024& $<10^{-10}$\\
7.8125e-03 & 1.9249e-02 & 7.8872e-04 & 7.5991e-04 & 22745 & 2048& $<10^{-10}$\\
3.9062e-03 & 9.6240e-03 & 7.7764e-04 & 4.7532e-04 & 50417 & 4096& $<10^{-10}$\\
\end{tabular}
\caption{Convergence test for \eqref{eq:truesol2d} over the time interval 
$[0,1]$ for the scheme \eqref{eq:fma2}.}\label{tab:2dvol}
% /home/rn/hpc_cluster/data/alberta/egn2/2d.convergence_experiment_NoC2/vol
\end{table}%
\begin{table}[H]
\center
\begin{tabular}{c|c|c|c|c|c|c}
% grep "LaTeX:" run*/specs.used -A 2 | grep '\\\\' | sed 's/.*specs.used- //'
% grep -A1 "TEC" run*/specs.used 
 $h_{f}$ & $h^M_\Gamma$ & $\errorUu$ & $\errorXx$ & $K^M_\Omega$ & $K^M_\Gamma$
 & $|v_\Delta^M|$ \\ \hline
6.2500e-02 & 1.5394e-01 & 7.5877e-03 & 1.8445e-03 & 2761 & 256  & 1.33e-02\\
3.1250e-02 & 7.7002e-02 & 1.5180e-03 & 1.0028e-03 & 5433 & 512  & 3.44e-03\\
1.5625e-02 & 3.8501e-02 & 1.3572e-03 & 1.1123e-03 & 10853 & 1024& 8.64e-04\\
7.8125e-03 & 1.9249e-02 & 7.7496e-04 & 6.9937e-04 & 22769 & 2048& 2.15e-04\\
3.9062e-03 & 9.6240e-03 & 7.7408e-04 & 4.6021e-04 & 50353 & 4096& 5.37e-05\\
\end{tabular}
\caption{Convergence test for \eqref{eq:truesol2d} over the time interval 
$[0,1]$ for the scheme \eqref{eq:fma}.}\label{tab:2d}
% /home/rn/hpc_cluster/data/alberta/egn2/2d.convergence_experiment_NoC2/lin
\end{table}%
\revised{%
As a comparison, we show the results of the same convergence experiment
for the schemes \eqref{eq:fma2} and \eqref{eq:fma} with mass-lumping
in Tables~\ref{tab:2dvolml} and \ref{tab:2dml}, respectively.}
\begin{table}[H]
\center
\begin{tabular}{c|c|c|c|c|c|c}
% grep "LaTeX:" run*/specs.used -A 2 | grep '\\\\' | sed 's/.*specs.used- //'
 $h_{f}$ & $h^M_\Gamma$ & $\errorUu$ & $\errorXx$ & $K^M_\Omega$ & $K^M_\Gamma$
 & $|v_\Delta^M|$ \\ \hline
6.2500e-02 &1.5398e-01 & 6.0357e-03 & 1.0375e-03 & 2721  & 256 & $<10^{-10}$\\ 
3.1250e-02 &7.7003e-02 & 2.0175e-03 & 1.4508e-03 & 5425  & 512 & $<10^{-10}$\\ 
1.5625e-02 &3.8500e-02 & 1.5274e-03 & 1.0748e-03 & 10869 & 1024& $<10^{-10}$\\ 
7.8125e-03 &1.9248e-02 & 8.4669e-04 & 6.1128e-04 & 22761 & 2048& $<10^{-10}$\\ 
3.9062e-03 &9.6238e-03 & 8.0877e-04 & 3.9738e-04 & 50321 & 4096& $<10^{-10}$\\ 
\end{tabular}
\caption{\revised{Convergence test for \eqref{eq:truesol2d} over the time interval 
$[0,1]$ for the scheme \eqref{eq:fma2} with mass-lumping.}}\label{tab:2dvolml}
% /home/rn/hpc_cluster/data/alberta/egn2/2d.ml.convergence_experiment_NoC2/vol
\end{table}%
\begin{table}[H]
\center
\begin{tabular}{c|c|c|c|c|c|c}
% grep "LaTeX:" run*/specs.used -A 2 | grep '\\\\' | sed 's/.*specs.used- //'
% grep -A1 "TEC" run*/specs.used 
 $h_{f}$ & $h^M_\Gamma$ & $\errorUu$ & $\errorXx$ & $K^M_\Omega$ & $K^M_\Gamma$
 & $|v_\Delta^M|$ \\ \hline
6.2500e-02 &1.5392e-01 & 8.8775e-03 & 3.0074e-03 & 2713  & 256  & 1.3143e-02\\
3.1250e-02 &7.6994e-02 & 2.7763e-03 & 4.8992e-04 & 5433  & 512  & 3.4191e-03\\
1.5625e-02 &3.8499e-02 & 1.4789e-03 & 8.3211e-04 & 10877 & 1024 & 8.6125e-04\\
7.8125e-03 &1.9248e-02 & 8.3425e-04 & 5.5075e-04 & 22785 & 2048 & 2.1483e-04\\
3.9062e-03 &9.6238e-03 & 8.0541e-04 & 3.8228e-04 & 50329 & 4096 & 5.3629e-05\\
\end{tabular}
\caption{\revised{Convergence test for \eqref{eq:truesol2d} over the time interval 
$[0,1]$ for the scheme \eqref{eq:fma} with mass-lumping.}}\label{tab:2dml}
% /home/rn/hpc_cluster/data/alberta/egn2/2d.ml.convergence_experiment_NoC2/lin
\end{table}%

\begin{figure}[H]
    \centering
    \includegraphics[keepaspectratio, scale=0.15]{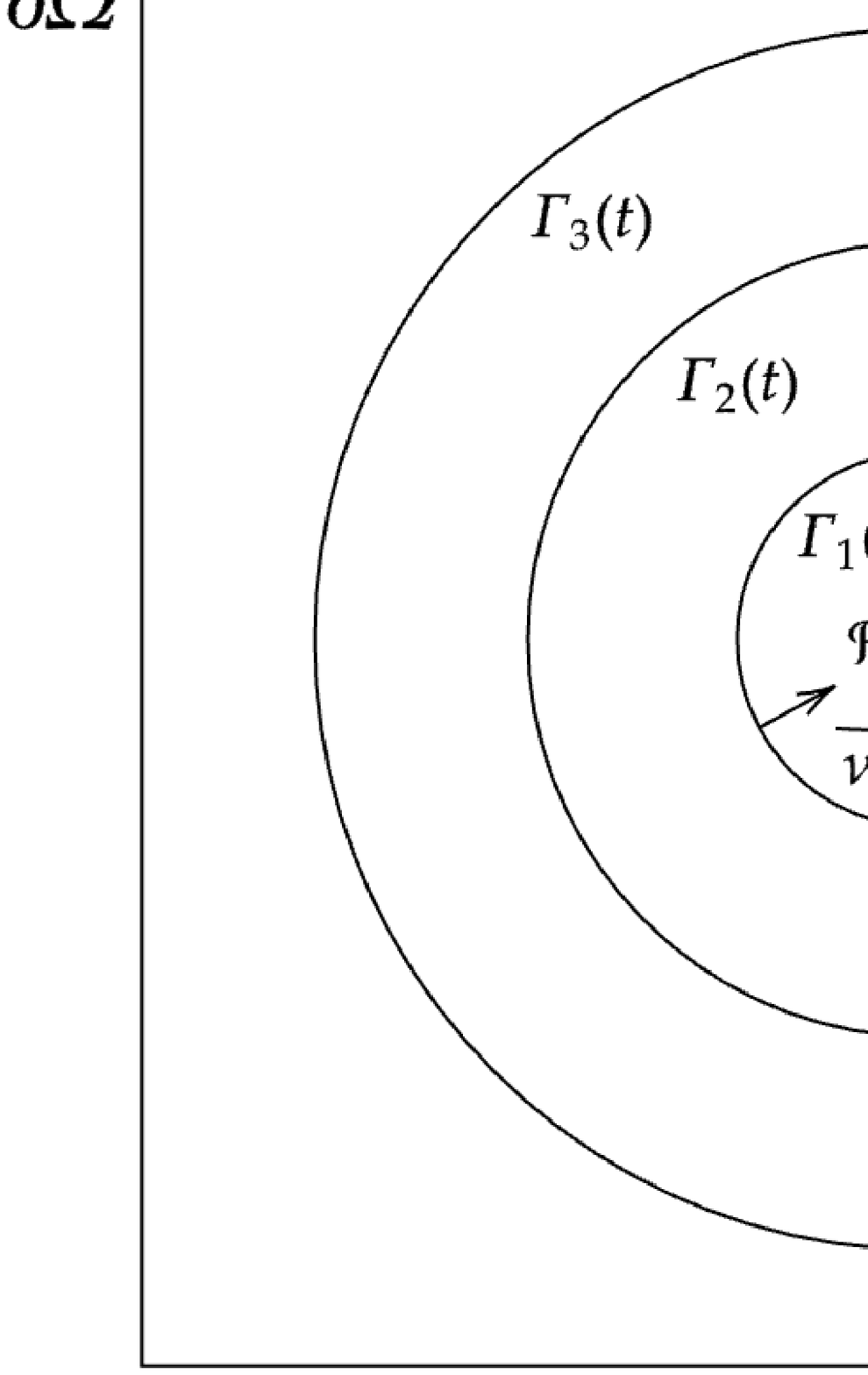}
    \caption{Three concentric circles.}\label{fig:3c}
\end{figure}

We proceed with another curve network consisting of three concentric circles in $\mathbb{R}^2$.
Let $R_1(t)$, $R_2(t)$ and $R_3(t)$ be the radii of the three circles at time $t$ with $R_1(t) < R_2(t) < R_3(t)$.
See \Figure{fig:3c} for the setting.

We note that the region $\subdomains{1}$ is not connected in contrast to the previous case (see \Figure{fig:2c}).
In this case, we see that $(\cpindexp{1},\cpindexm{1}) = (1,2)$, $(\cpindexp{2},\cpindexm{2})=(3,2)$ and $(\cpindexp{3},\cpindexm{3})=(3,1)$,
and we have a rigorous solution to \eqref{eq:sharp_p} with $\Gamma(t) = \bigcup_{i=1}^3\partial B(0,R_i(t))$ as follows:
\begin{equation*}
    w(x,t) := \begin{cases}
        -\frac{1}{\jump{2}{1}{\beta}R_1(t)}&\qquad\mbox{if}\quad |x| < R_1(t),\\
        \frac{\frac{1}{\jump{2}{3}{\beta}R_2(t)} + \frac{1}{\jump{2}{1}{\beta}R_1(t)}}{\log{\frac{R_2(t)}{R_1(t)}}}\log{\frac{|x|}{R_1(t)}} - \frac{1}{\jump{2}{1}{\beta}R_1(t)}&\qquad\mbox{if}\quad R_1(t)\leq |x| < R_2(t),\\
        \frac{\frac{1}{\jump{2}{3}{\beta}R_2(t)} + \frac{1}{\jump{1}{3}{\beta}R_3(t)}}{\log{\frac{R_2(t)}{R_3(t)}}}\log{\frac{|x|}{R_3(t)}} - \frac{1}{\jump{1}{3}{\beta}R_3(t)}&\qquad\mbox{if}\quad R_2(t)\leq |x| < R_3(t),\\
        -\frac{1}{\jump{1}{3}{\beta}R_3(t)}&\qquad\mbox{if}\quad R_3(t)\leq |x|,
    \end{cases}
\end{equation*}
with $R_1(t)$, $R_2(t)$ and $R_3(t)$ satisfying the following ordinary differential equations:
\begin{align}\label{eq:ode}
    \begin{split}
        \dot{R}_1(t) &= - \frac{\frac{1}{\jump{2}{3}{\beta} R_2(t)} + \frac{1}{\jump{2}{1}{\beta}R_1(t)}}{\jump{2}{1}{\beta} R_1(t)\log{\frac{R_2(t)}{R_1(t)}}},\\
        \dot{R}_2(t) &= \frac{1}{\jump{2}{3}{\beta}R_2(t)}\left(\frac{\frac{1}{\jump{2}{3}{\beta}R_2(t)} + \frac{1}{\jump{1}{3}{\beta}R_3(t)}}{\log{\frac{R_2(t)}{R_3(t)}}} - \frac{\frac{1}{\jump{2}{3}{\beta}R_2(t)} + \frac{1}{\jump{2}{1}{\beta}R_1(t)}}{\log{\frac{R_2(t)}{R_1(t)}}}\right),\\
        \dot{R}_3(t) &= \frac{\frac{1}{\jump{2}{3}{\beta} R_2(t)} + \frac{1}{\jump{1}{3}{\beta}R_3(t)}}{\jump{1}{3}{\beta} R_3(t)\log{\frac{R_2(t)}{R_3(t)}}}.
    \end{split}
\end{align}
    % From the diffuse equation of \Eqref{eq:sharp_p}, $\smallChemical$ can be written as
    % $\smallChemical\mid_{\subdomains{2}} = c_2 + \alpha_2\log{|x|}$ and $w\mid_{\subdomains{3}} = c_3 + \alpha_3\log{|x|}$
    % for some $c_2$, $c_3$, $\alpha_2$ and $\alpha_3$ which are independent of $x$.
    % Moreover, we may assume that $w$ is constant in $\subdomains{1}$. Then, the Gibbs--Thomson law yields the values of $c_2$, $c_3$, $\alpha_2$ and $\alpha_3$, and the value of $\smallChemical$  in $\subdomains{1}$ is determined according to the values of $\smallChemical$ on $\partial\subdomains{1}$ and $\partial\subdomains{3}\backslash\partial\subdomains{2}$.
    % Finally, the motion law yields the ordinary differential equations for $R_1$, $R_2$ and $R_3$.

\begin{figure}[H]
% ../../plotradii
% cp radii_clean.ps ~/tex/harald/toku2/figures/simul3c_r.ps && scpp radii_clean.ps e23:tex/harald/toku2/figures/simul3c_r.ps
% cp radii_clean.ps ~/tex/harald/toku2/figures/simul3c_long_r.ps && scpp radii_clean.ps e23:tex/harald/toku2/figures/simul3c_long_r.ps
% \includegraphics[angle=-90,width=0.75,scale=0.5\textwidth]{figures/simul3c_r}
\centering
\includegraphics[angle=-90, keepaspectratio, scale=0.25]{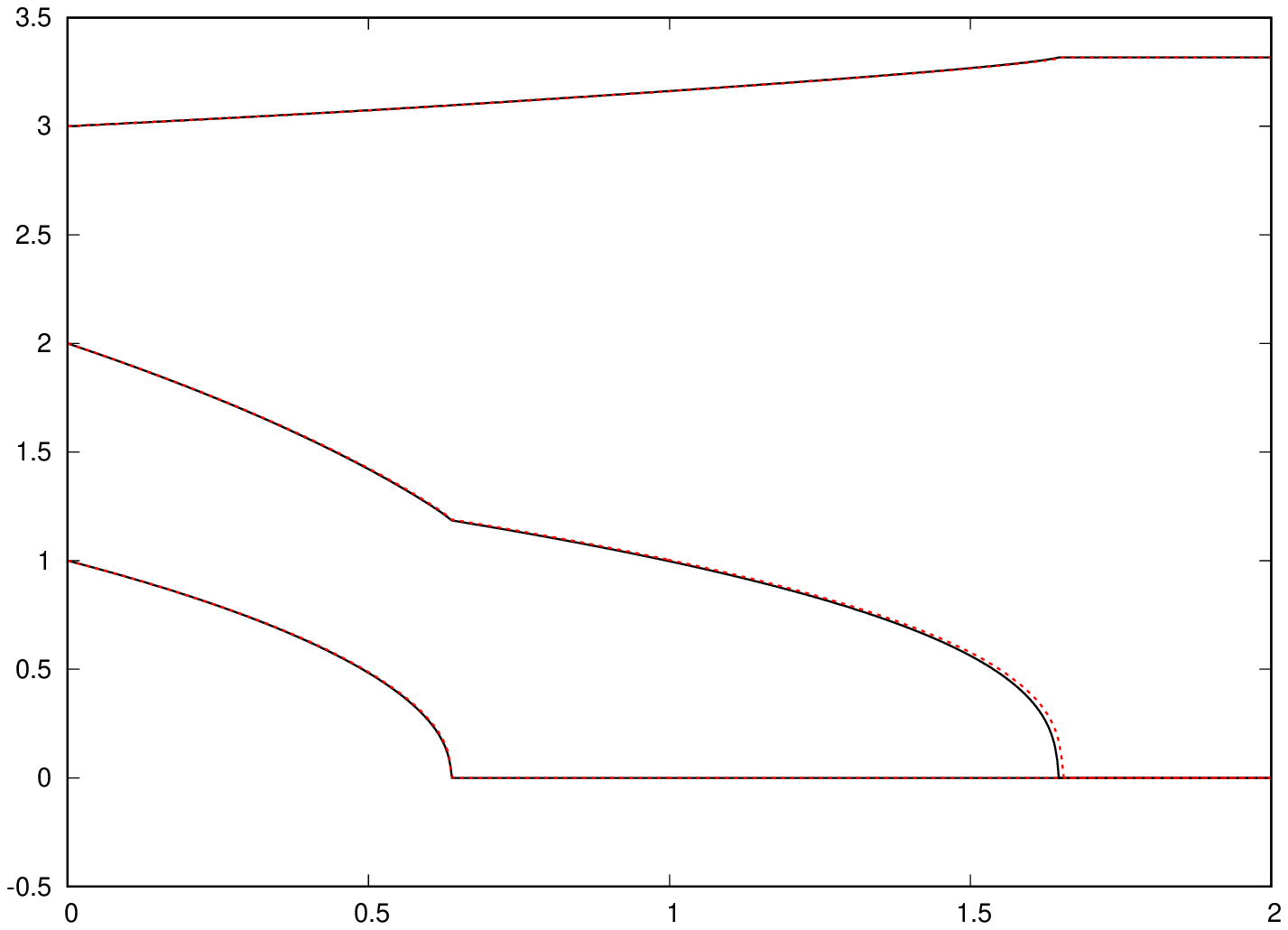}
\setlength{\abovecaptionskip}{20pt}
\caption{Comparison of the discrete (black, solid) \\
and exact (red, dashed) radii for $(\beta_1,\beta_2,\beta_3) = (-1,1,0)$.}
% ~/hpc_cluster/data/alberta/egn2/2d.convergence_experiment_NoP3/lin/run2
% ~/hpc_cluster/data/alberta/egn2/2d.convergence_experiment_NoP3/plot/run3
\label{fig:simul3c}
\end{figure}

We show in Figure~\ref{fig:simul3c} a comparison between the solution
to \eqref{eq:ode}, and the discrete finite element solution to \eqref{eq:fma}
that approximates the system \eqref{eq:sharp_p}. 
Once again, we note the excellent agreement between the two. During the evolution, the two
innermost disks disappear, leaving a trivial two-phase stationary solution, 
involving only phases 1 and 3. Observe that once the inner
disk vanishes, the remaining radii $R_2(t)$ and $R_3(t)$ satisfy a
system of ODEs analogous to \eqref{eq:rd2},
involving three phases and two concentric circles.

\subsection{Numerical simulations in 2d}
We shall carry out several numerical simulations for the degenerate multi-phase Stefan problem \eqref{eq:sharp_p} in the case $d = 2$.

\begin{figure}[H]
% ../plotradii; mv radii.ps 2d_multiple049.ps && extractdata 0 1 2 3 4 5 uv.out && plotcurve; mv extracted.ps 2d_multiple049_ts.ps
% ../plotradii; mv radii.ps 2d_multiple05.ps && extractdata 0 1 uv.out && plotcurve; mv extracted.ps 2d_multiple05_ts.ps
% ../plotradii; mv radii.ps 2d_multiple051.ps && extractdata 0 1 2 3 4 5 6 7 uv.out && plotcurve; mv extracted.ps 2d_multiple051_ts.ps
% fixbb 2d_*.ps && enlargepsfont_all 2d_*.ps && cp 2d_*.ps ~/tex/harald/toku2/figures && scpp 2d_*.ps e23:tex/harald/toku2/figures/ 
\centering
\includegraphics[angle=-90,width=0.3\textwidth]{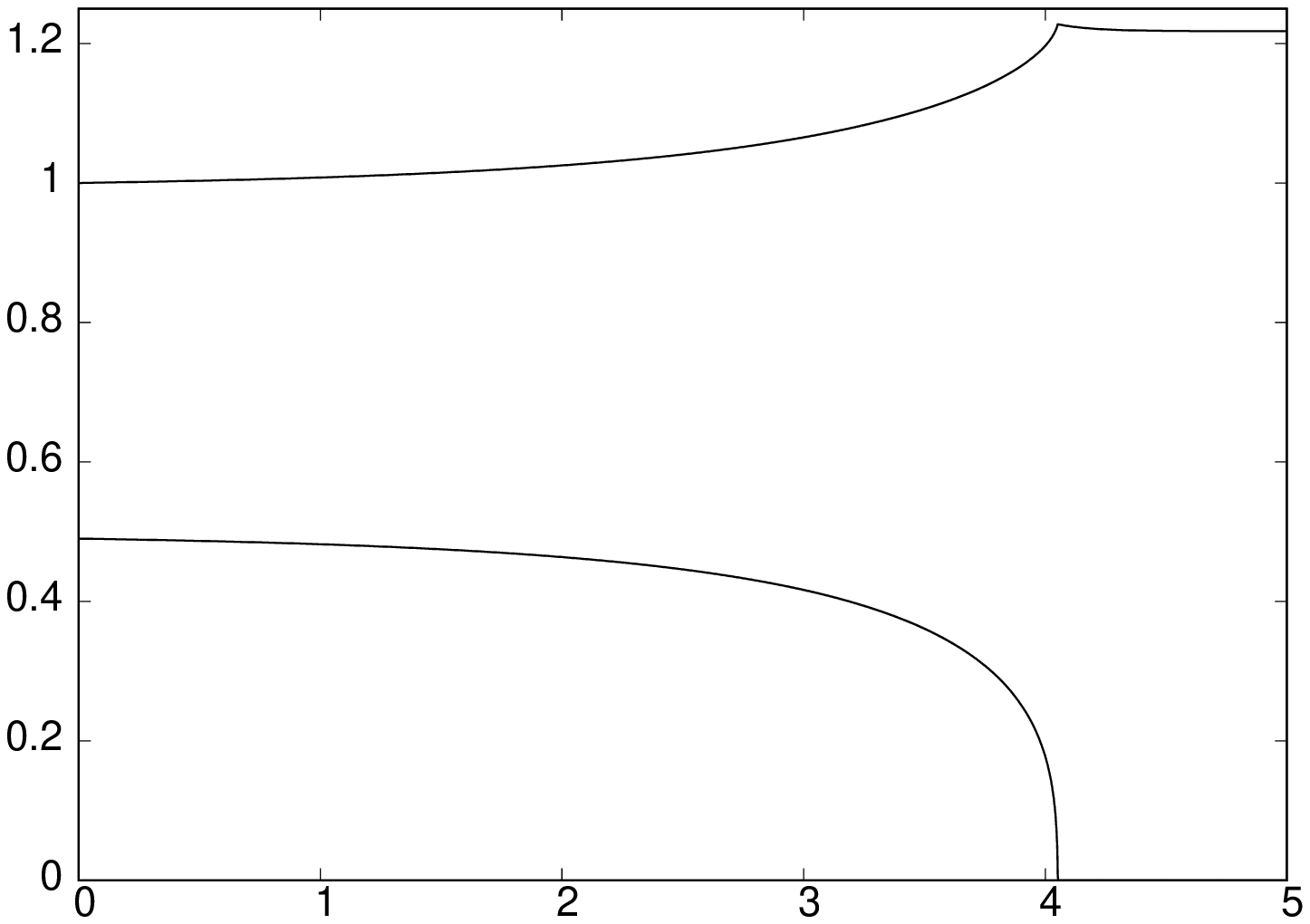}
\includegraphics[angle=-90,width=0.3\textwidth]{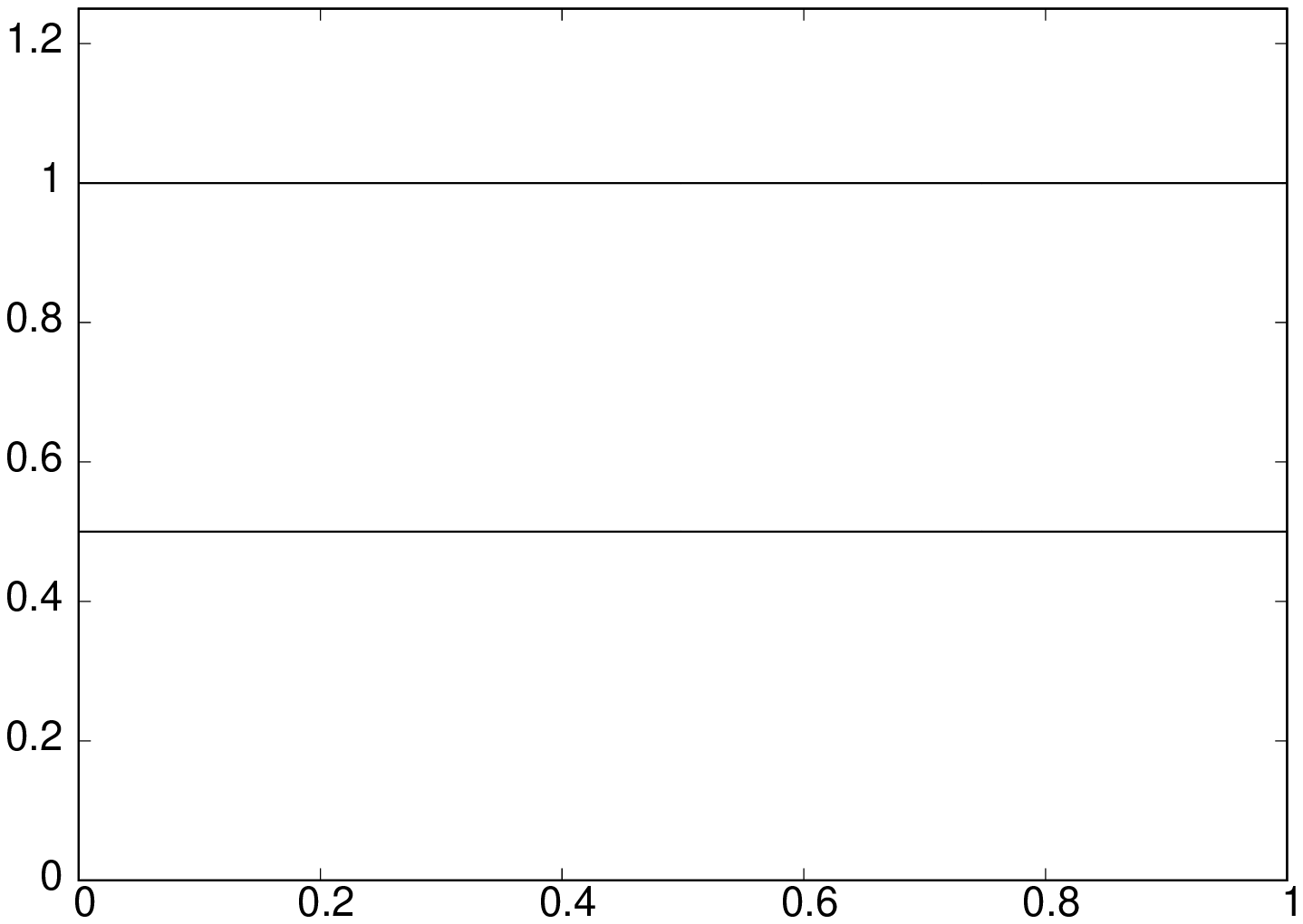}
\includegraphics[angle=-90,width=0.3\textwidth]{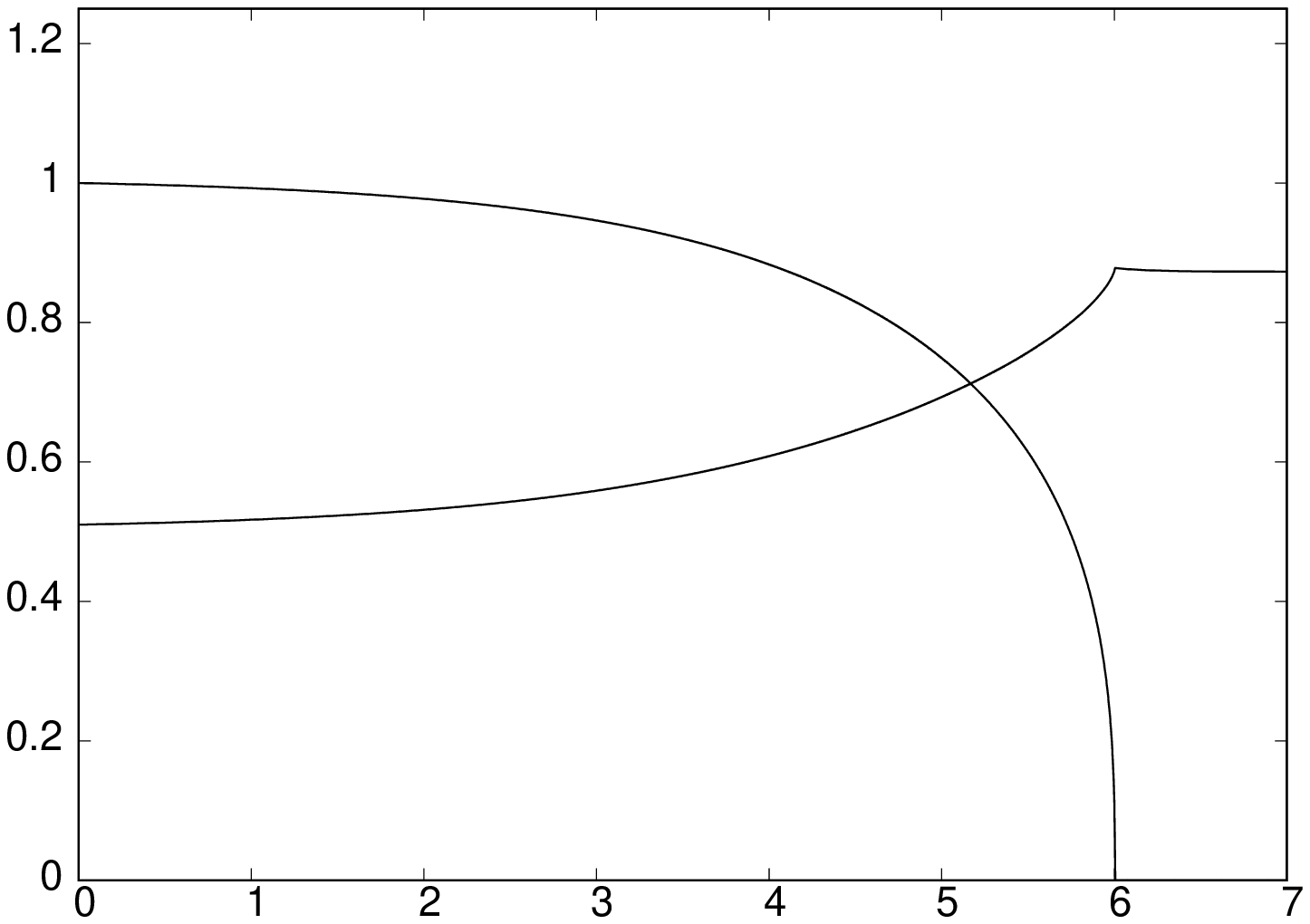}
\includegraphics[angle=-90,width=0.3\textwidth]{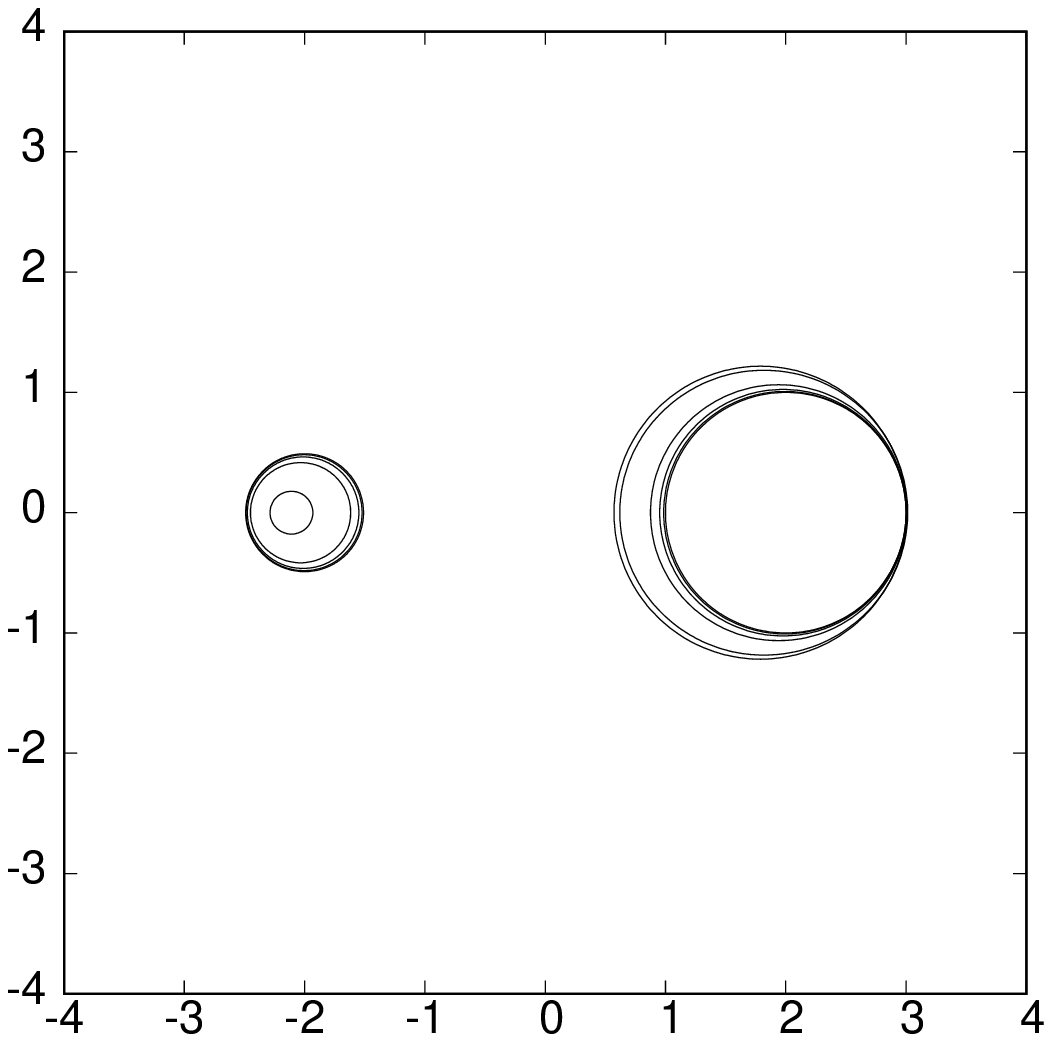}
\includegraphics[angle=-90,width=0.3\textwidth]{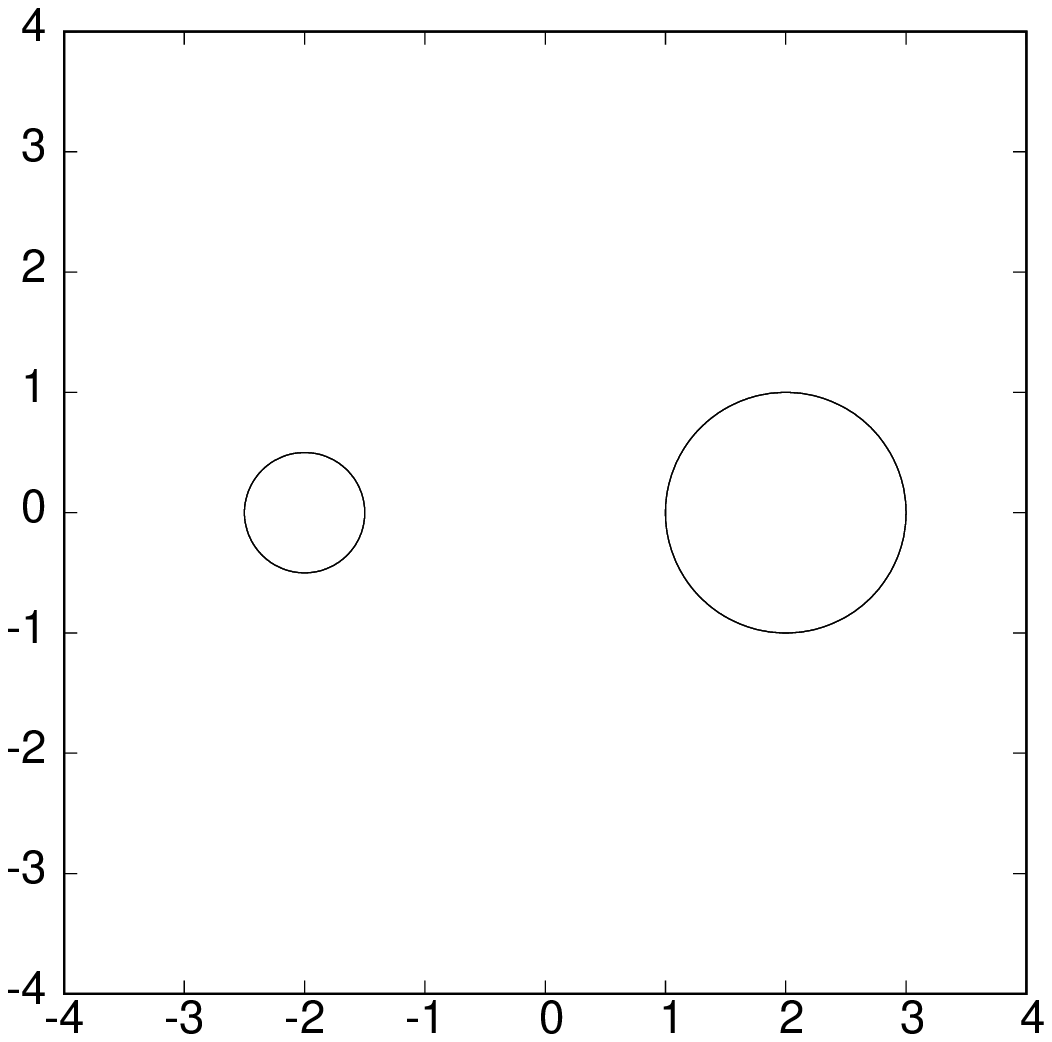}
\includegraphics[angle=-90,width=0.3\textwidth]{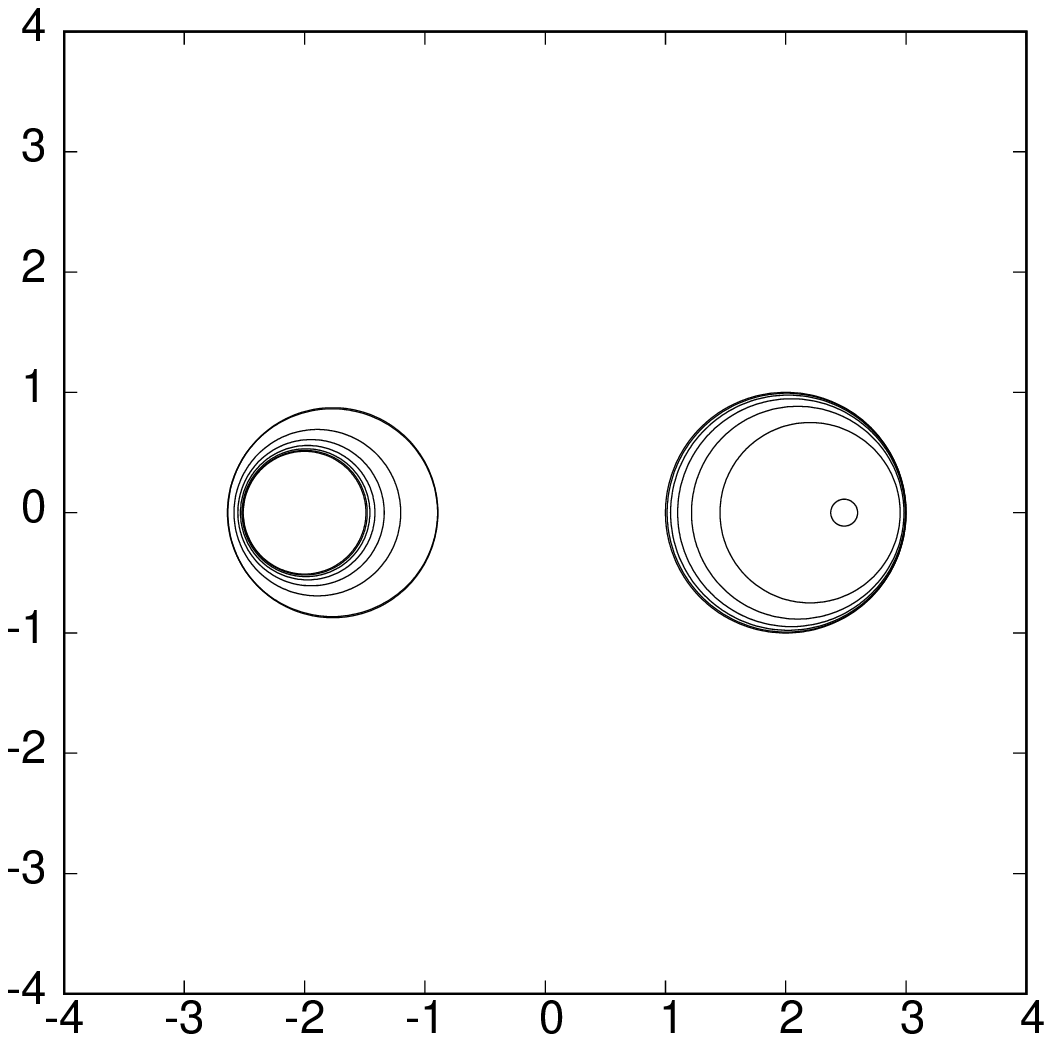}
\setlength{\abovecaptionskip}{20pt}
\caption{%
The evolution of half of the maximum extension in $x$-direction of the
two phases over time. For the initial radius $0.49$, $0.5$ or $0.51$ (from left
to right) of the initially smaller disk. Below we show plots of 
the evolutions at times $t=0,1,\ldots,T$ with $T=5$, $1$ and $7$, respectively.}
\label{fig:multiple}
% /home/rn/hpc_cluster/data/alberta/egn2/2d.multiple_049
% /home/rn/hpc_cluster/data/alberta/egn2/2d.multiple_05
% /home/rn/hpc_cluster/data/alberta/egn2/2d.multiple_051
\end{figure}

We begin with a simulation for $(\beta_1,\beta_2,\beta_3)=(0, -1, 1)$ for a stationary solution,
that is two separate disks with radii $1$ and $0.5$, respectively.
When the smaller initial radius is smaller than $0.5$, the smaller phase disappears.
When it is slightly larger than $0.5$ it grows and the initially larger phase
eventually disappears. We visualize the three different behaviors in
Figure~\ref{fig:multiple}, where we choose the initially smaller disk to have
radius $0.49$, $0.5$ or $0.51$, respectively.

\begin{figure}[H]
% rsync -av "e23:$PWD/*.out" .
% extractdata 0 uv.out && plotcurve; mv extracted.ps 2d_db_t0.ps; extractdata 0.1 uv.out && plotcurve; mv extracted.ps 2d_db_t01.ps; extractdata 0.2 uv.out && plotcurve; mv extracted.ps 2d_db_t02.ps; extractdata 0.5 uv.out && plotcurve; mv extracted.ps 2d_db_t05.ps; extractdata 2 uv.out && plotcurve; mv extracted.ps 2d_db_t2.ps; plotenergy; mv energy.ps 2d_db_e.ps
% Matlab: plottriangulation
% mv mesh.ps 2d_db_mesh0.ps
% mv mesh.ps 2d_db_mesh1.ps
% fixbb 2d_*.ps && enlargepsfont_all 2d_*.ps && cp 2d_*.ps ~/tex/harald/toku2/figures && scpp 2d_*.ps e23:tex/harald/toku2/figures/ 
\centering
\includegraphics[angle=-90,width=0.18\textwidth]{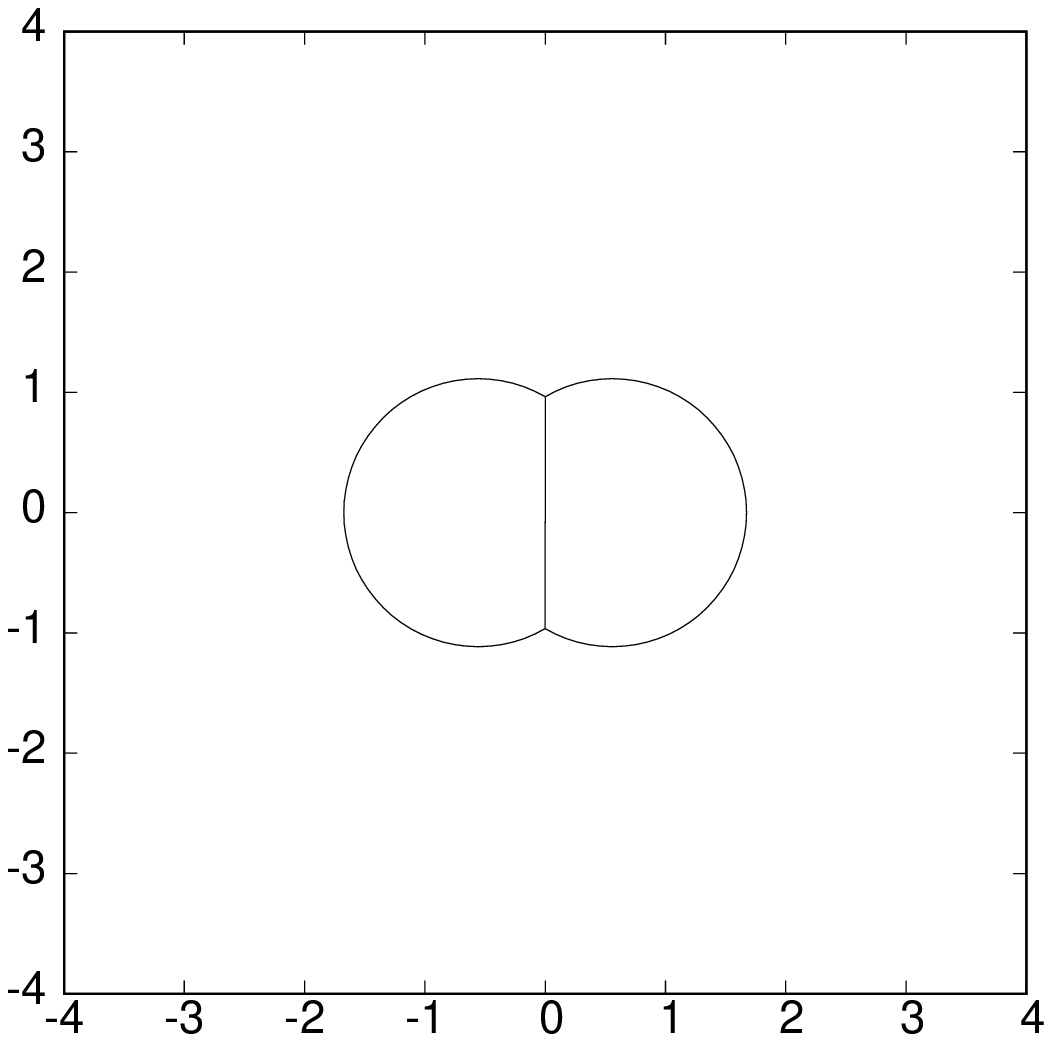}
\includegraphics[angle=-90,width=0.18\textwidth]{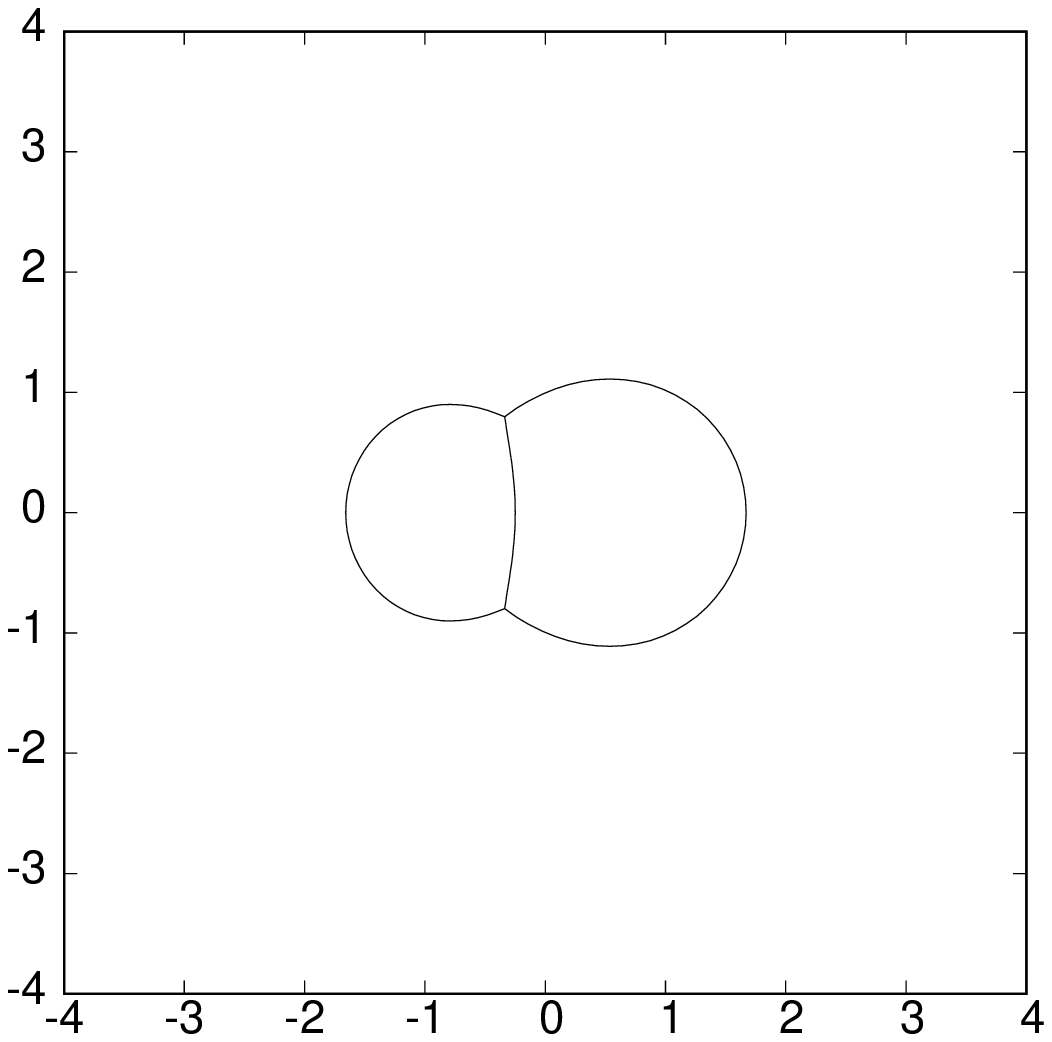}
\includegraphics[angle=-90,width=0.18\textwidth]{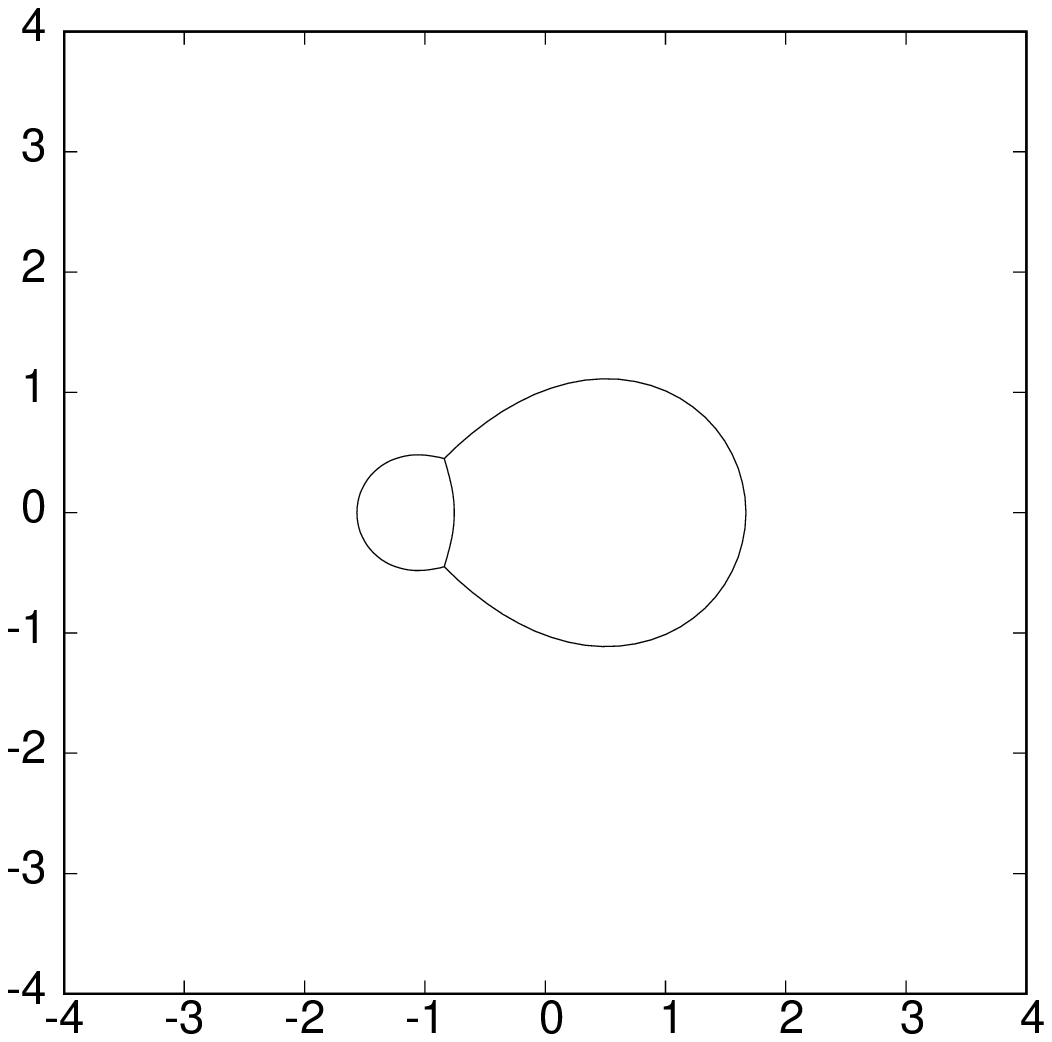}
\includegraphics[angle=-90,width=0.18\textwidth]{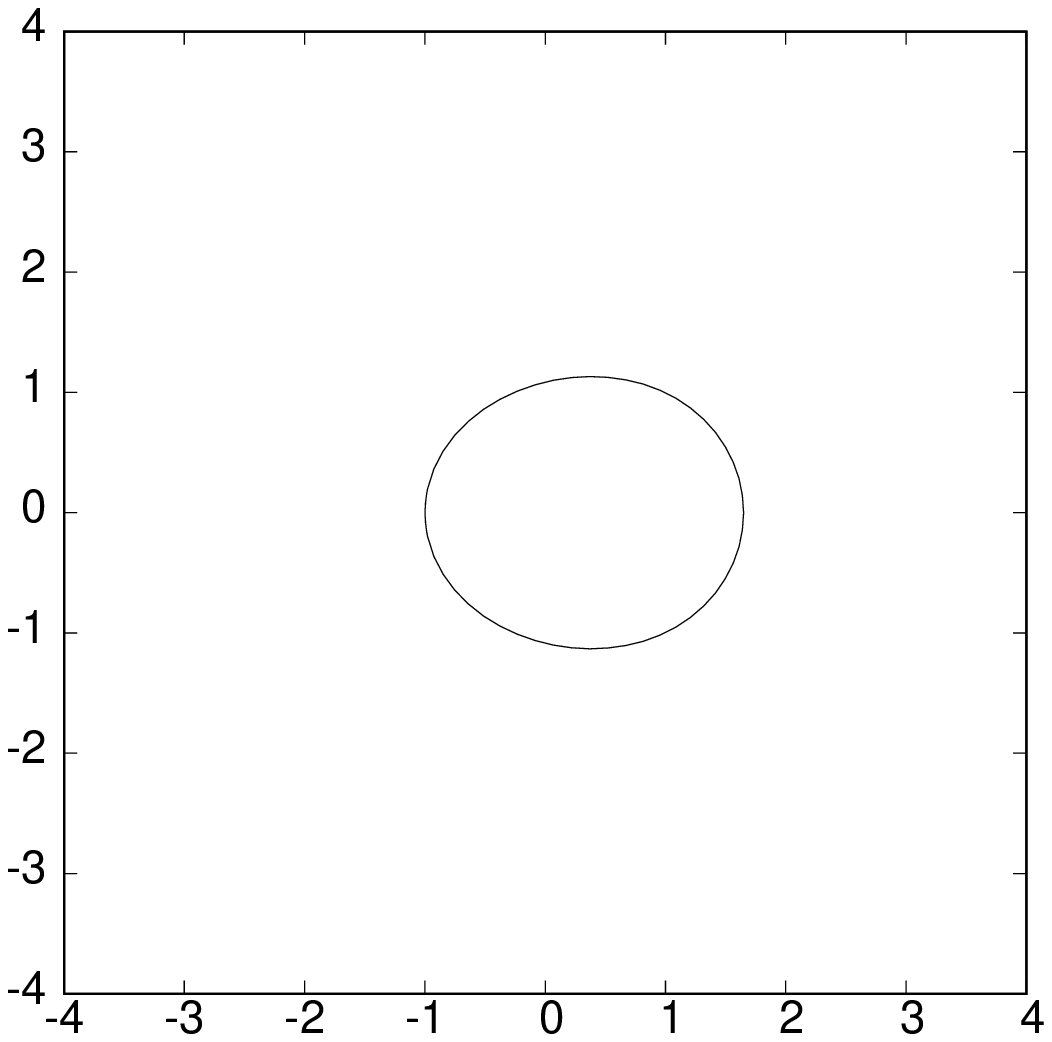}
\includegraphics[angle=-90,width=0.18\textwidth]{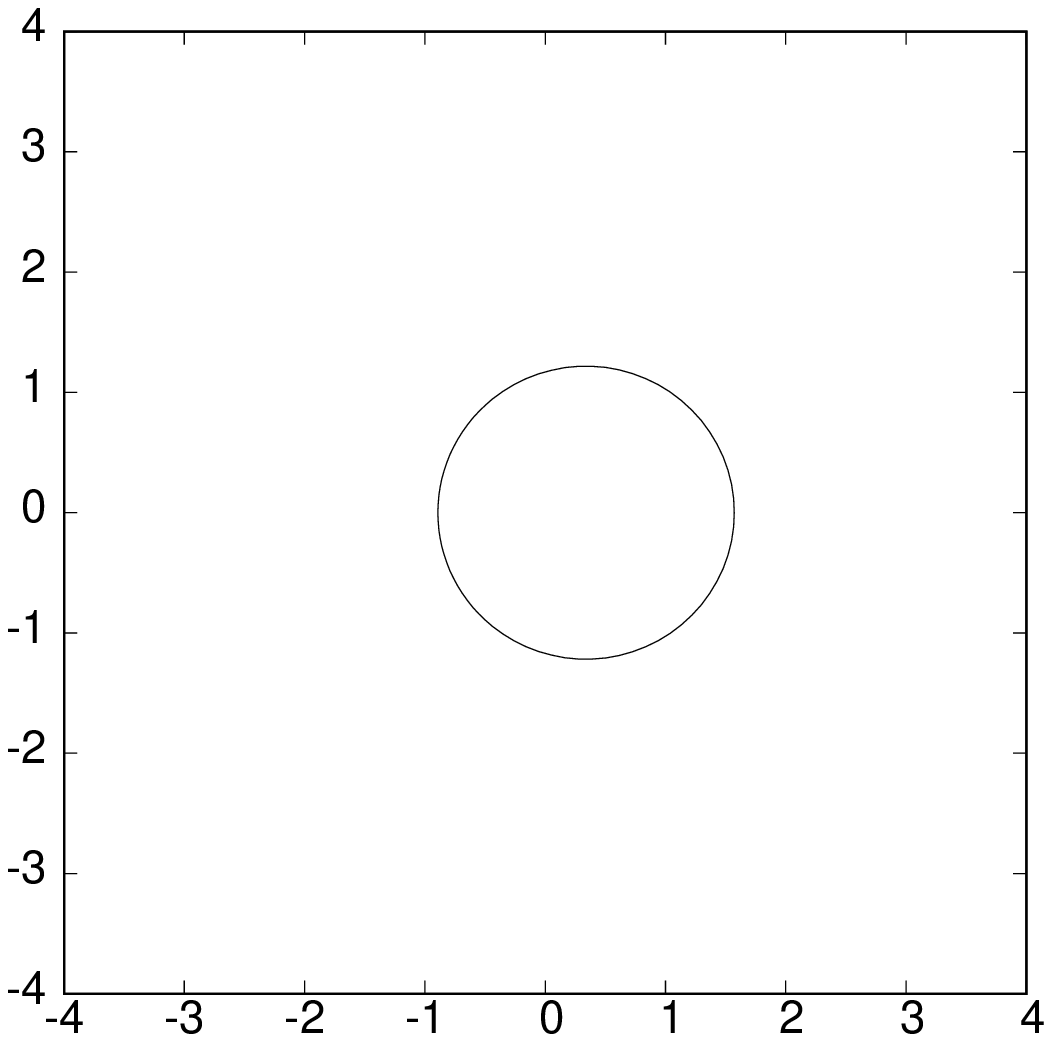}
\includegraphics[angle=-90,width=0.3\textwidth]{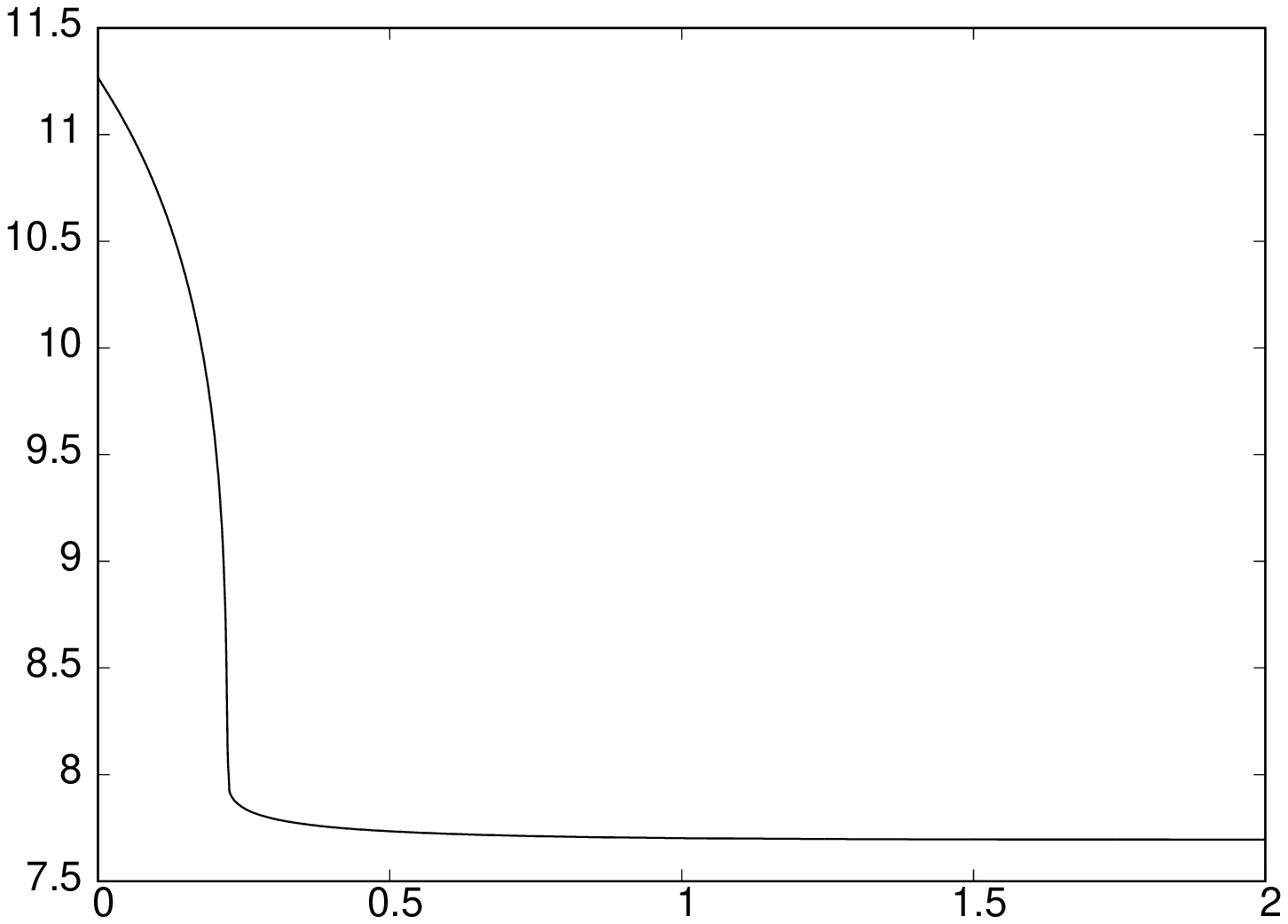} 
%\\[5mm]
%\includegraphics[angle=-0,width=0.18\textwidth]{figures/2d_db_mesh0}
%\includegraphics[angle=-0,width=0.18\textwidth]{figures/2d_db_mesh1}
\setlength{\abovecaptionskip}{20pt}
\caption{(($\beta_1, \beta_2,\beta_3) = (-1,0,1) =$ right bubble, left bubble, outer phase)\\
The solution at times $t=0, 0.1, 0.2, 0.5, 2$, 
and a plot of the discrete energy over time. 
For this computation, we have $\discreteVol{0} = 54.6$.
%Below we show the adaptive bulk mesh at times $t=0$ and $t=1$.
}
\label{fig:2d_db}
% /home/rn/hpc_cluster/data/alberta/egn2/2d.db
% /home/rn/hpc_cluster/data/alberta/egn2/2d.db/vol3
\end{figure}

We proceed with a simple example of a symmetric double bubble in 2d.
In this experiment, we set $I_R = 3$, $I_S = 3$, $I_T=2$, $(\beta_1, \beta_2,\beta_3) = (-1,0,1)$,
$(\cpindexp{1},\cpindexm{1}) = (3,1)$, $(\cpindexp{2}, \cpindexm{2}) = (2,3)$, and $(\cpindexp{3}, \cpindexm{3}) = (1,2)$.
We show its result in Figure~\ref{fig:2d_db} and observe that the left bubble shrinks and the right bubble grows.
Eventually, the left bubble vanishes, and the right bubble and the outer phase survive.
We note that without heuristic changes of topology our numerical methods 
would not be able to integrate beyond the vanishing of a (part of) a phase.
Hence, we need to perform some simply topological surgery to allow our scheme
to continue with the approximation the evolution. To this end, and similarly to
our previous work \cite{EtoGarckeNurnberg2024}, we discard (parts of) phases
that have become too small. If such a vanishing part is enclosed by a single
curve, it is simply discarded. Otherwise, one of the bounding curves is removed,
and the triple junctions it was part of are changed to intersections of the 
two remaining curves.

We next carry out a simulation for the same double bubble, but with a different choice of surface tensions.
Then, we observe a different behavior of the evolution of the double bubble. In spite of the same choice of $\phaseContent$, now
the left bubble grows, and the right bubble shrinks and eventually vanishes.
We note also that the three contact angles at the triple junctions in this
experiment are different from each other, as is to expected from Young's law, see the fourth condition in \eqref{eq:sharp_p}, in contrast to the previous case.
\begin{figure}[H]
% rsync -av "e23:$PWD/*.out" .
% extractdata 0 uv.out && plotcurve; mv extracted.ps 2d_db_a_t0.ps; extractdata 0.2 uv.out && plotcurve; mv extracted.ps 2d_db_a_t02.ps; extractdata 0.5 uv.out && plotcurve; mv extracted.ps 2d_db_a_t05.ps; extractdata 1 uv.out && plotcurve; mv extracted.ps 2d_db_a_t1.ps; extractdata 2 uv.out && plotcurve; mv extracted.ps 2d_db_a_t2.ps; plotenergy; mv energy.ps 2d_db_a_e.ps
% fixbb 2d_*.ps && enlargepsfont_all 2d_*.ps && cp 2d_*.ps ~/tex/harald/toku2/figures && scpp 2d_*.ps e23:tex/harald/toku2/figures/ 
\center
\includegraphics[angle=-90,width=0.18\textwidth]{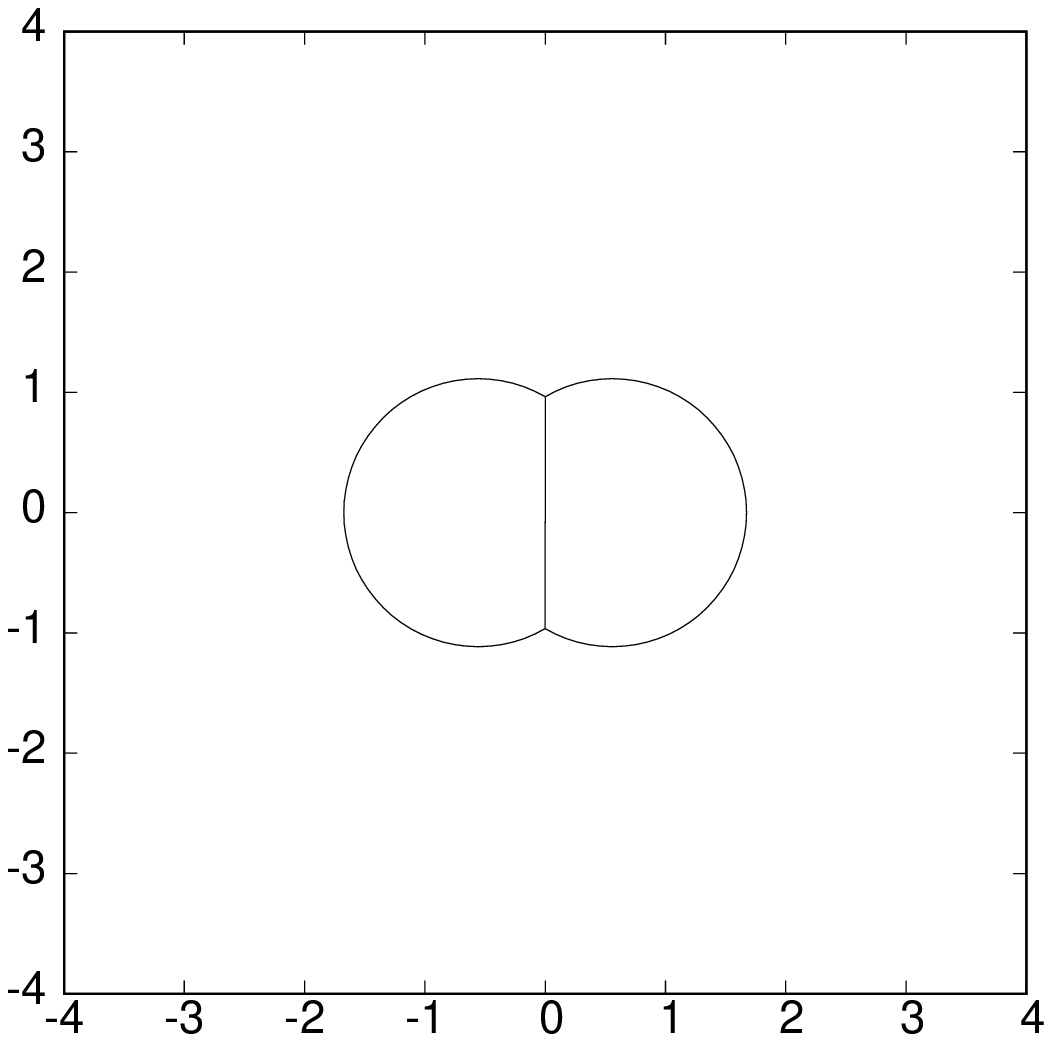}
\includegraphics[angle=-90,width=0.18\textwidth]{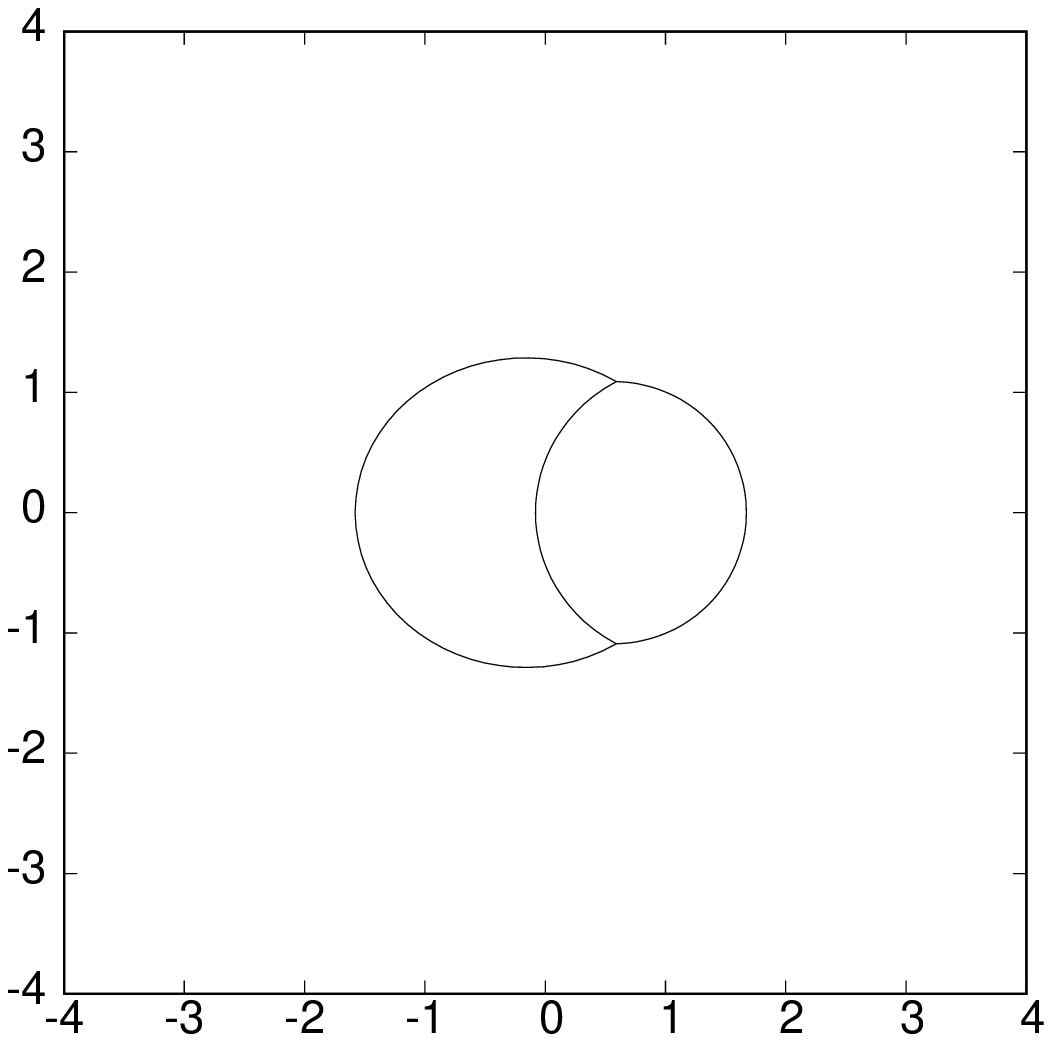}
\includegraphics[angle=-90,width=0.18\textwidth]{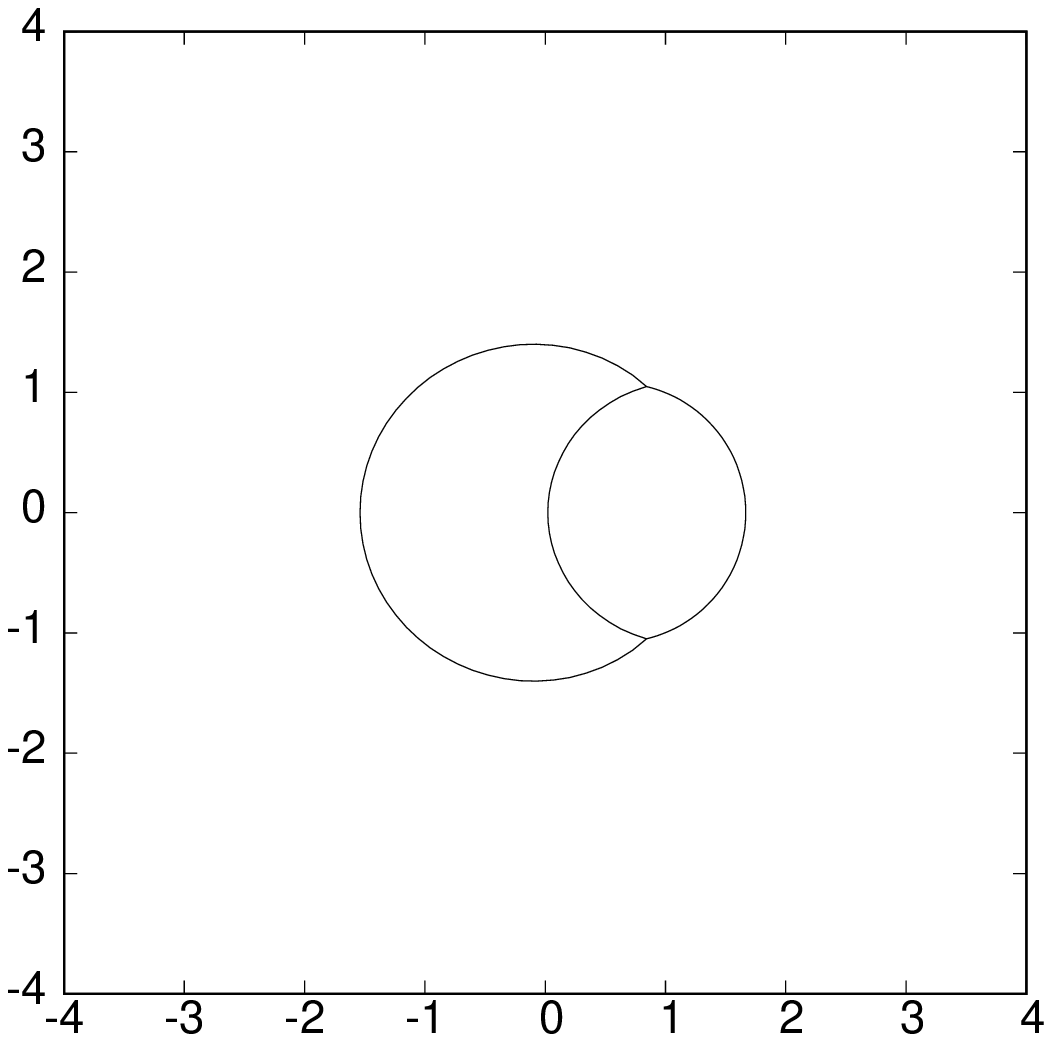}
\includegraphics[angle=-90,width=0.18\textwidth]{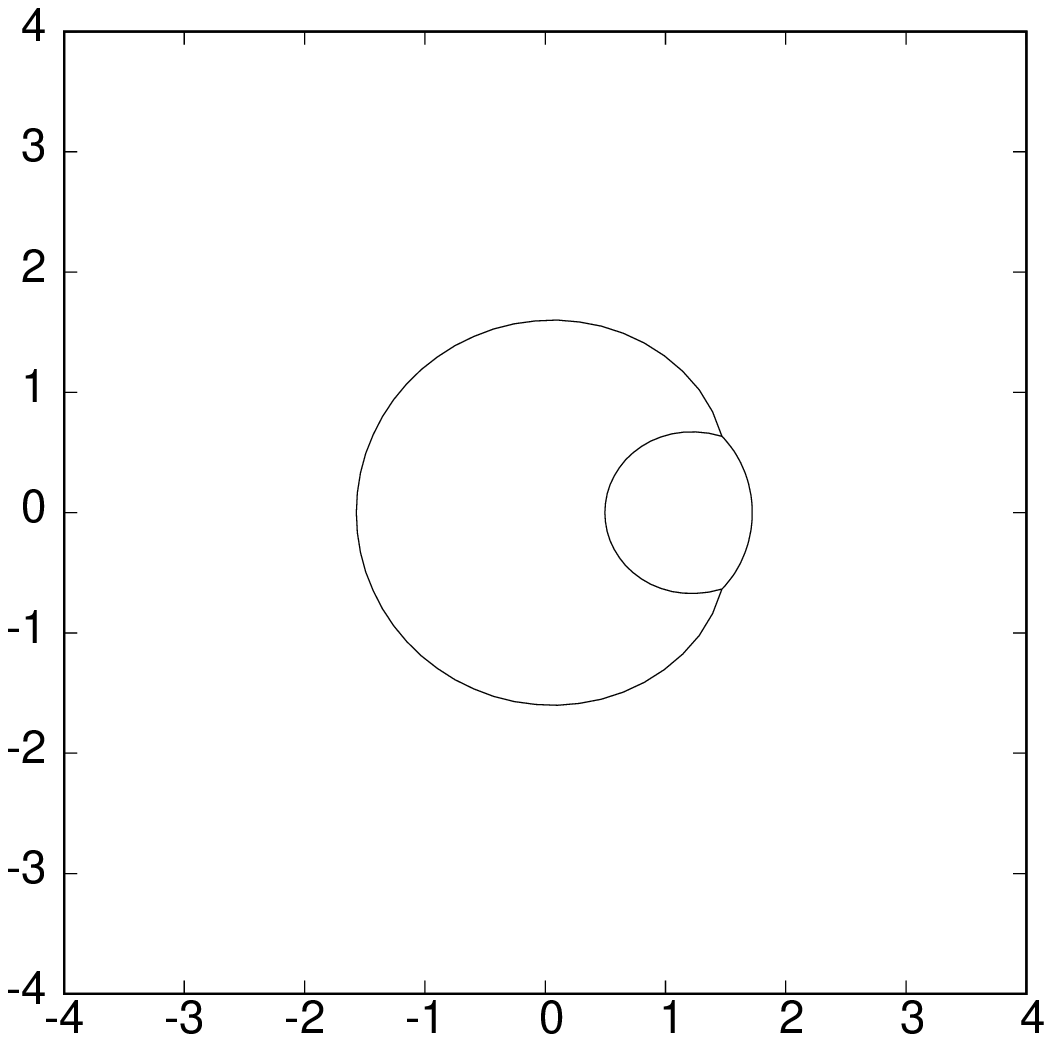}
\includegraphics[angle=-90,width=0.18\textwidth]{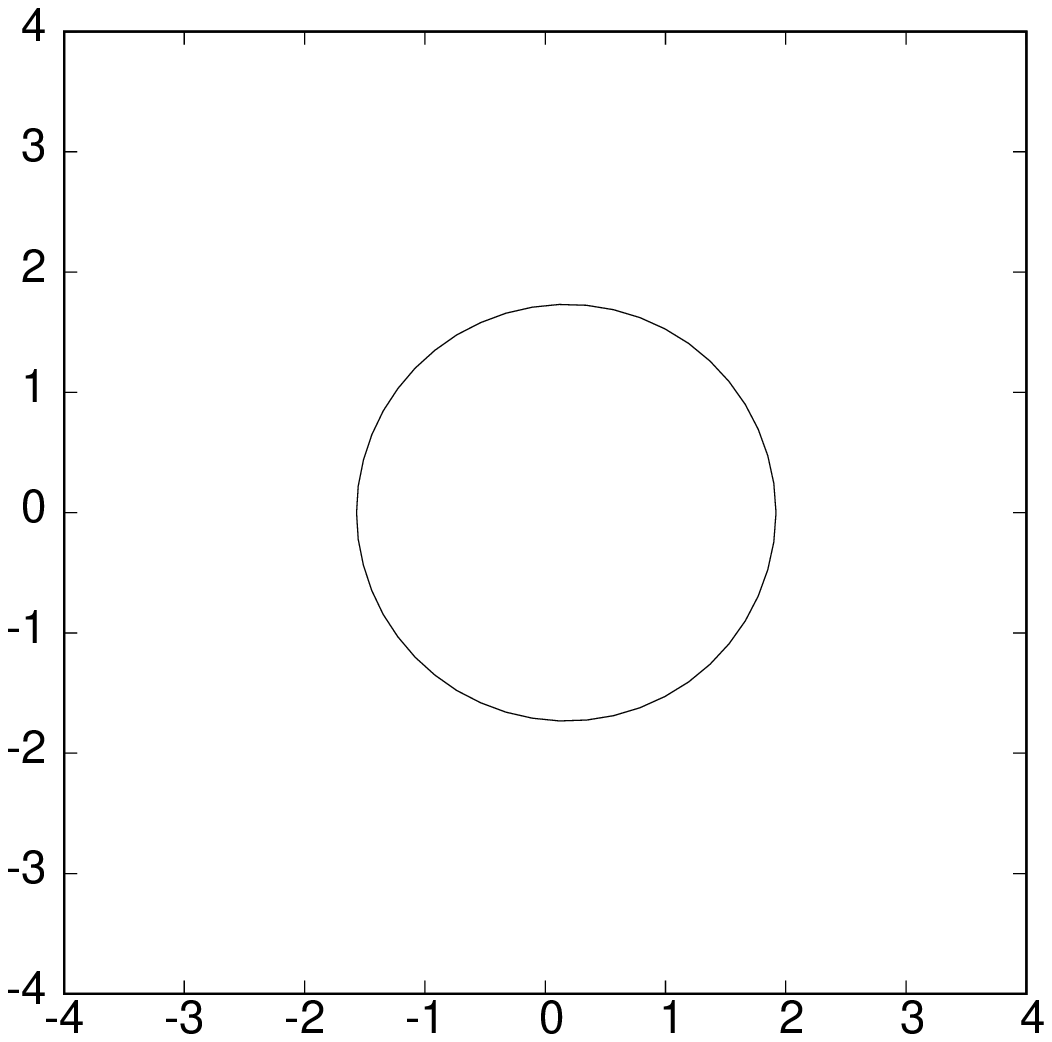}
\includegraphics[angle=-90,width=0.3\textwidth]{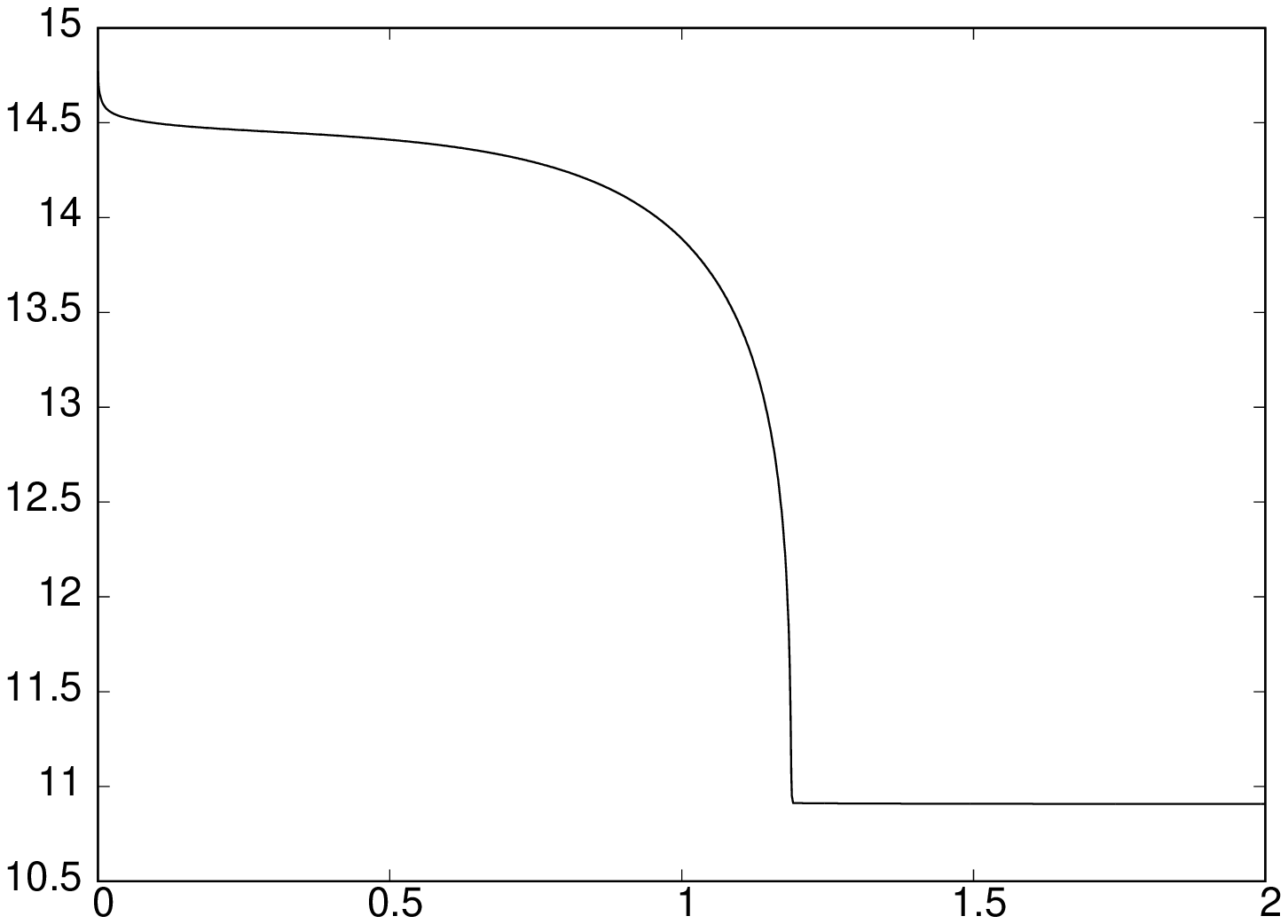} 
\setlength{\abovecaptionskip}{20pt}
\caption{
(($\beta_1, \beta_2,\beta_3) = (-1,0,1) =$ right bubble, left bubble,
outer phase, and with $\tension{} = (1.75, 1,1)$).
The solution at times $t=0, 0.2, 0.5, 1, 2$, 
and a plot of the discrete energy over time.
For this computation, we have $\discreteVol{0} = 54.6$.
}
\label{fig:2d_db_a}
% /home/rn/hpc_cluster/data/alberta/egn2/2d.db_a
\end{figure}%

\begin{figure}[H]
% rsync -av "e23:$PWD/*.out" .
% extractdata 0 uv.out && plotcurve; mv extracted.ps 2d_db_plus_onet0.ps; extractdata 0.015 uv.out && plotcurve; mv extracted.ps 2d_db_plus_onet0015.ps; extractdata 0.02 uv.out && plotcurve; mv extracted.ps 2d_db_plus_onet002.ps; extractdata 0.1 uv.out && plotcurve; mv extracted.ps 2d_db_plus_onet01.ps; extractdata 0.2 uv.out && plotcurve; mv extracted.ps 2d_db_plus_onet02.ps; plotenergy; mv energy.ps 2d_db_plus_onee.ps
% Matlab: plottriangulation
% mv mesh.ps 2d_db_plus_onemesh0.ps
% mv mesh.ps 2d_db_plus_onemesh1.ps
% fixbb 2d_*.ps && enlargepsfont_all 2d_*.ps && cp 2d_*.ps ~/tex/harald/toku2/figures && scpp 2d_*.ps e23:tex/harald/toku2/figures/ 
\centering
\includegraphics[angle=-90,width=0.18\textwidth]{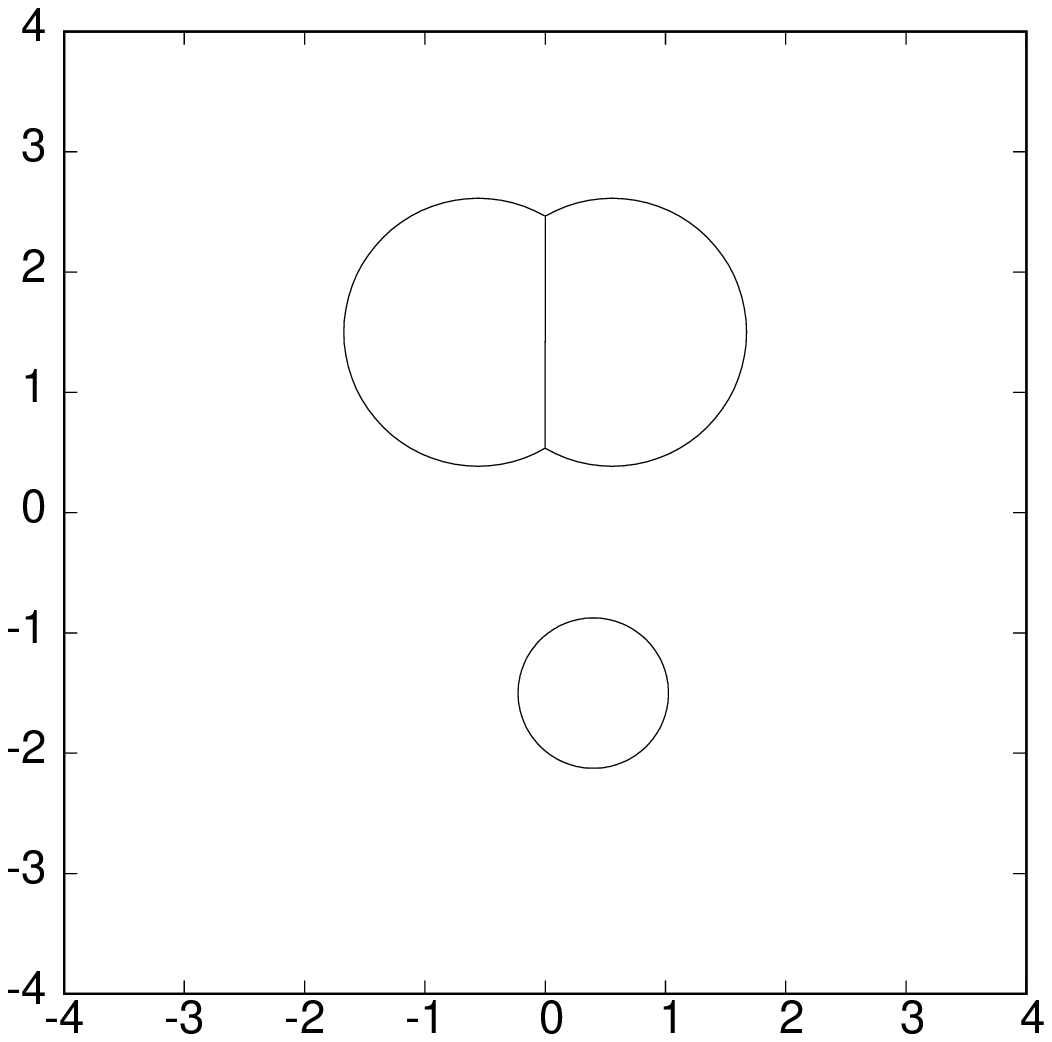}
\includegraphics[angle=-90,width=0.18\textwidth]{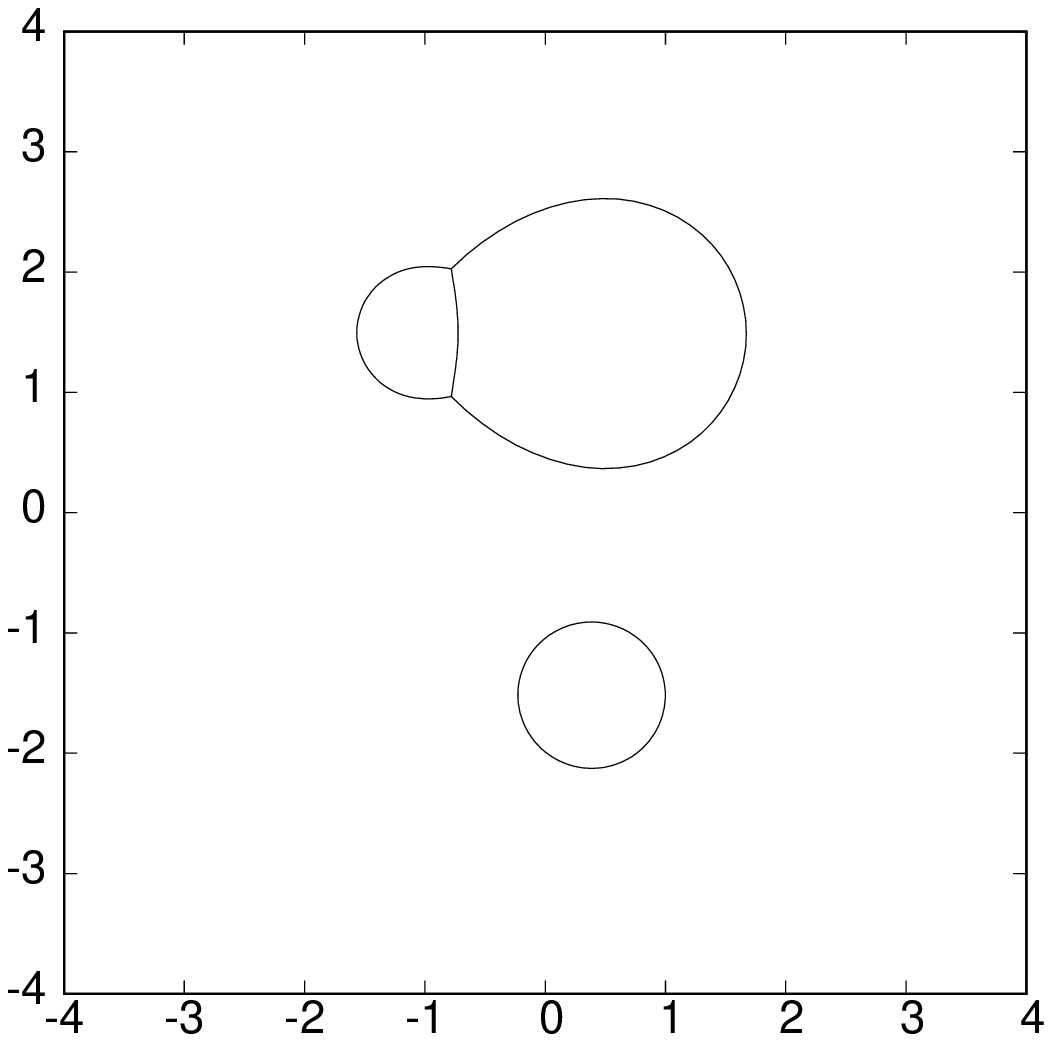}
\includegraphics[angle=-90,width=0.18\textwidth]{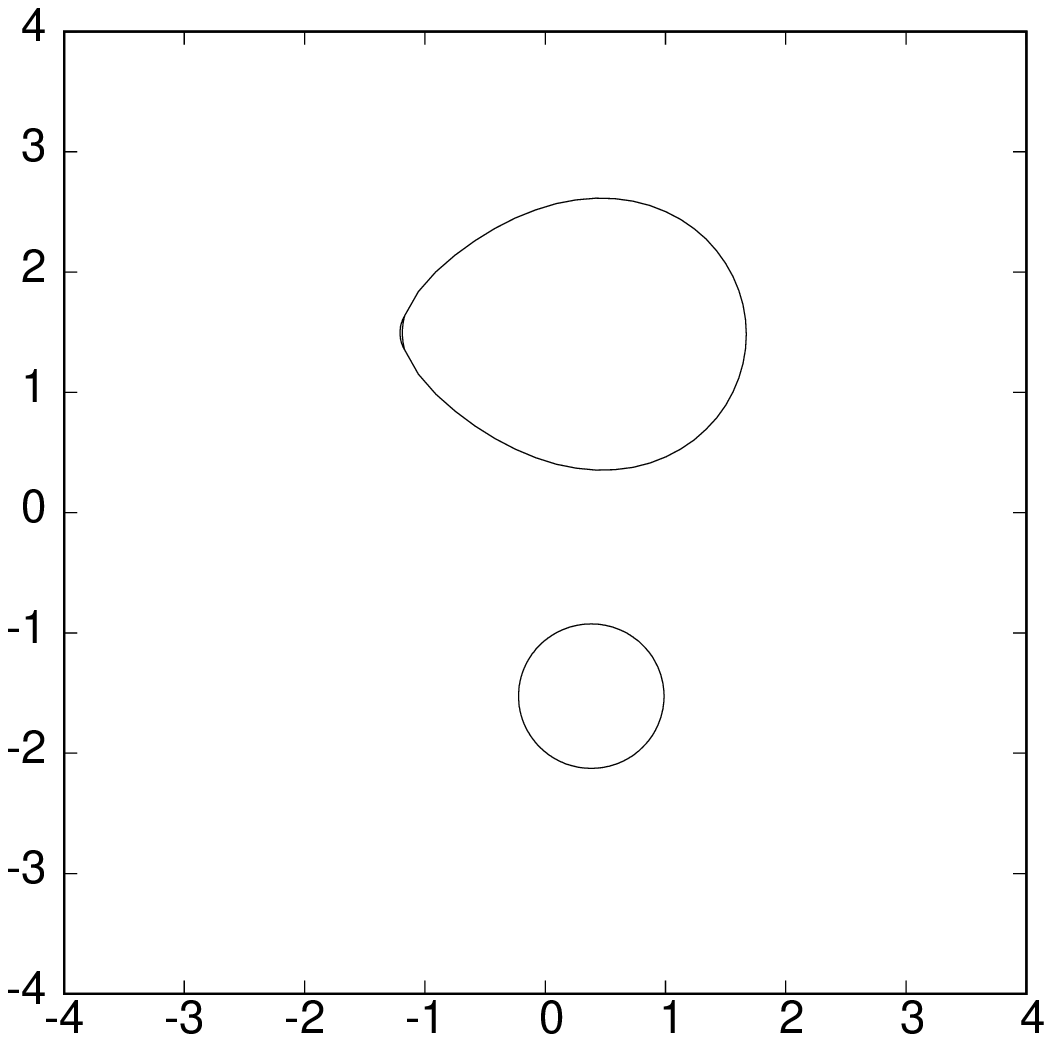}
\includegraphics[angle=-90,width=0.18\textwidth]{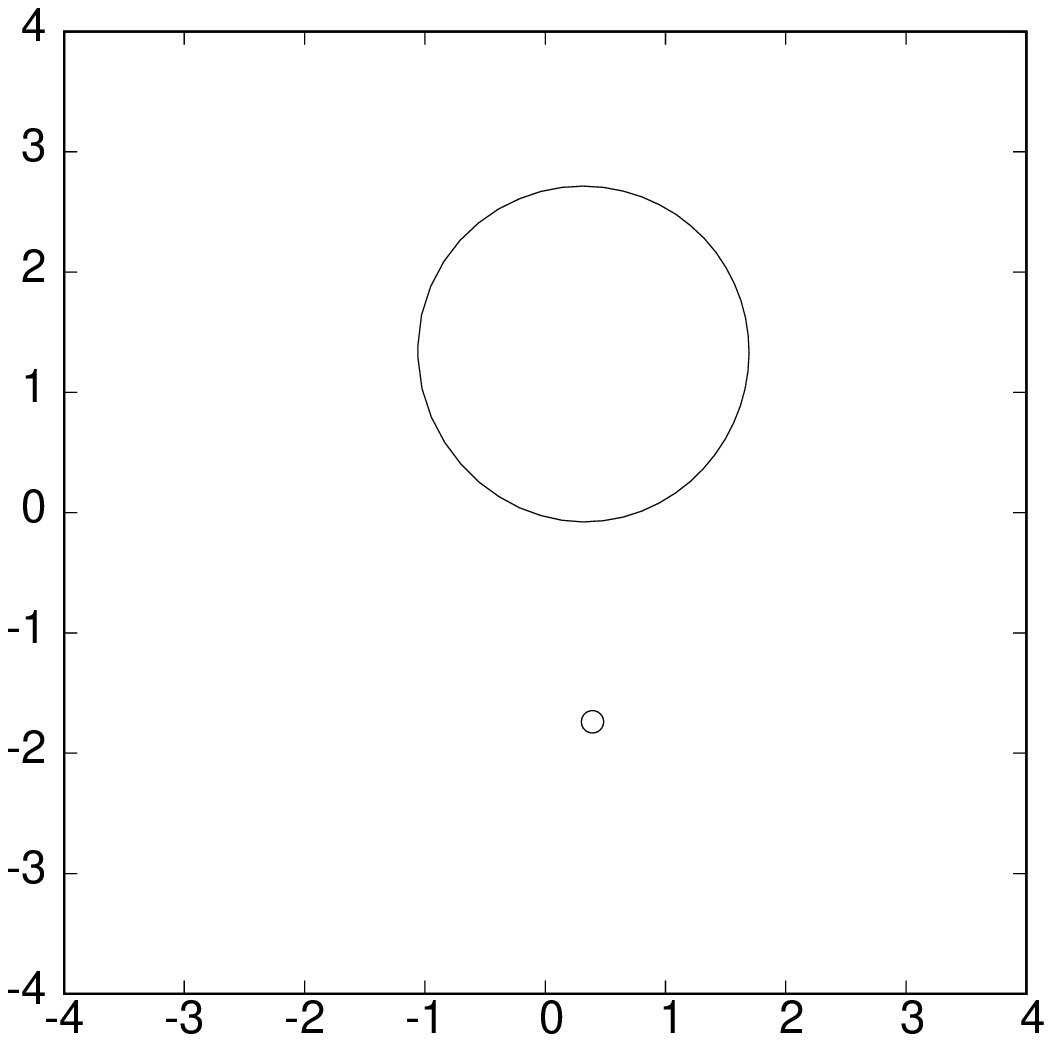}
\includegraphics[angle=-90,width=0.18\textwidth]{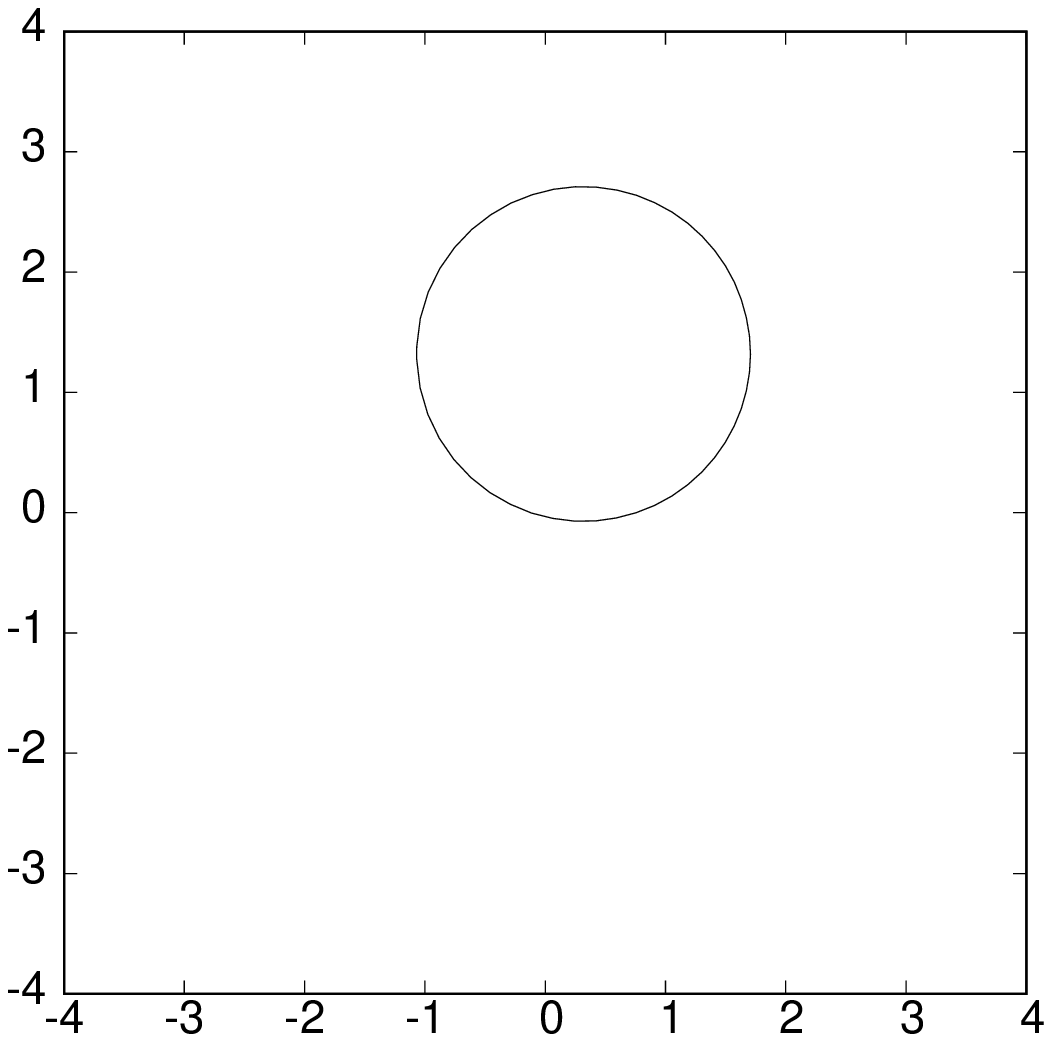}
\includegraphics[angle=-90,width=0.3\textwidth]{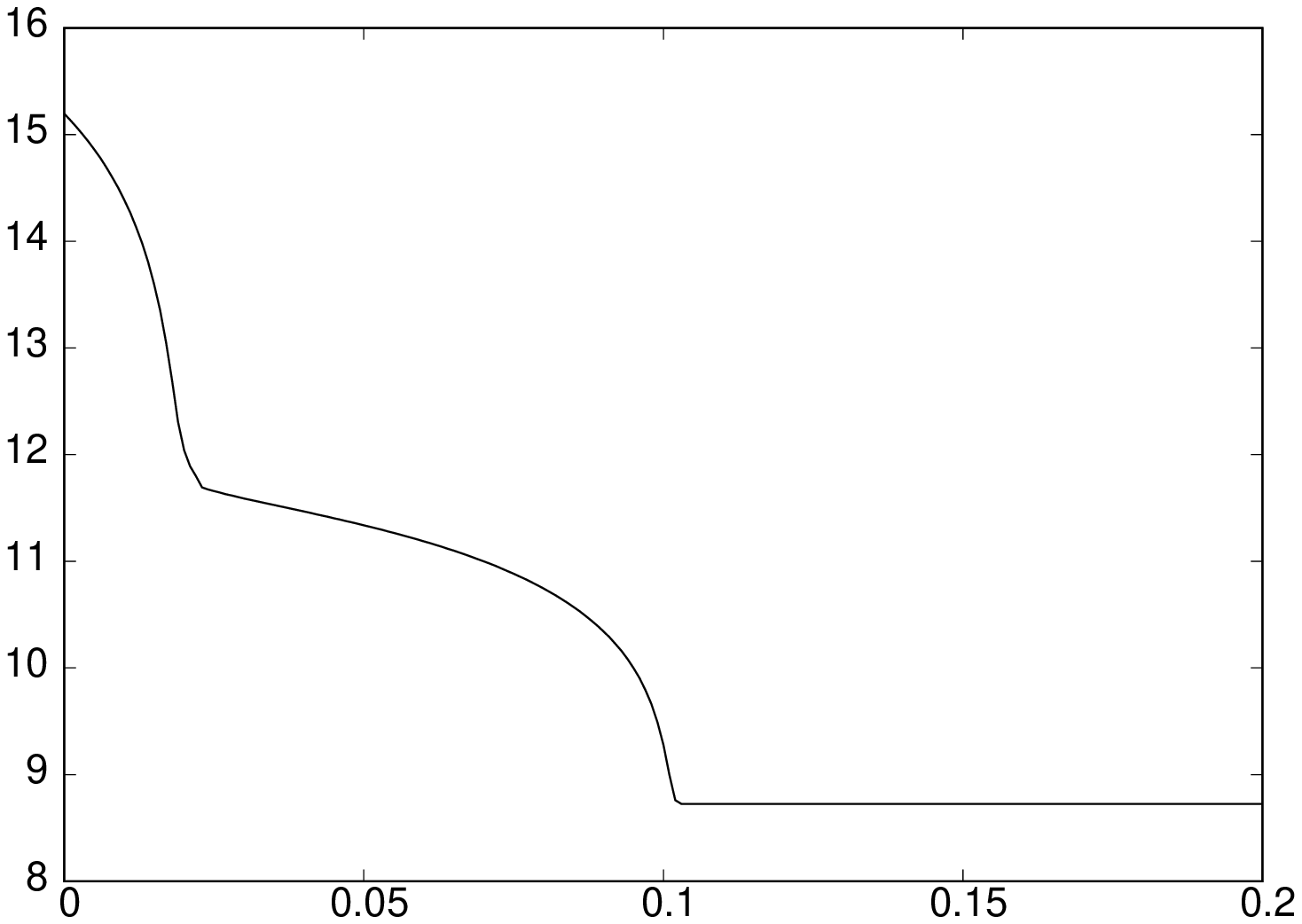} 
%\\[5mm]
%\includegraphics[angle=-0,width=0.18\textwidth]{figures/2d_db_plus_onemesh0}
%\includegraphics[angle=-0,width=0.18\textwidth]{figures/2d_db_plus_onemesh1}
\setlength{\abovecaptionskip}{20pt}
\caption{(($\beta_1, \beta_2,\beta_3) = (0.25,0,-0.25) =$ right bubble plus disk, left bubble, outer phase)\\
The solution at times $t=0, 0.015, 0.02, 0.1, 0.2$,
and a plot of the discrete energy over time.
For this computation, we have $\discreteVol{0} = -13$.
%Below we show the adaptive bulk mesh at times $t=0$ and $t=1$.
}
\label{fig:2d_db_plus_one}
% /home/rn/hpc_cluster/data/alberta/egn2/2d.db_plus_one
% /home/rn/hpc_cluster/data/alberta/egn2/2d.db_plus_one/vol3
\end{figure}

In the next set of experiments, we investigate simulations 
for a standard double bubble, with one of the bubbles making a phase
with a separate disk. In particular, we set $I_R = 3$, $I_S = 4$, $I_T = 2$, $(\beta_1, \beta_2,\beta_3) = (0.25,0,-0.25)$,
$(\cpindexp{1},\cpindexm{1}) = (3,1)$, $(\cpindexp{2},\cpindexm{2}) = (2,3)$, $(\cpindexp{3},\cpindexm{3}) = (1,2)$, and $(\cpindexp{4},\cpindexm{4}) = (1,3)$.
The two bubbles of the double bubble enclose an area about 3.139 each, while the disk has an initial radius of $\frac{8}{5}$,
meaning it initially encloses an area of $\frac{25\pi}{64} \approx 1.227$.
During the evolution, the left bobble shrinks, while the right bubble grows correspondingly.
Eventually, the left bobble vanishes, and the right bubble and the disk survive.
These survivors enclose the same phase, and the disk is absorbed by the right bubble according to the energy minimization.
See Figure~\ref{fig:2d_db_plus_one} for the results.

\begin{figure}[H]
% rsync -av "e23:$PWD/*.out" .
% extractdata 0 uv.out && plotcurve; mv extracted.ps 2d_db_plus_bigonet0.ps; extractdata 0.015 uv.out && plotcurve; mv extracted.ps 2d_db_plus_bigonet0015.ps; extractdata 0.02 uv.out && plotcurve; mv extracted.ps 2d_db_plus_bigonet002.ps; extractdata 1 uv.out && plotcurve; mv extracted.ps 2d_db_plus_bigonet1.ps; extractdata 2 uv.out && plotcurve; mv extracted.ps 2d_db_plus_bigonet2.ps; plotenergy; mv energy.ps 2d_db_plus_bigonee.ps
% Matlab: plottriangulation
% mv mesh.ps 2d_db_plus_bigonemesh0.ps
% mv mesh.ps 2d_db_plus_bigonemesh1.ps
% fixbb 2d_*.ps && enlargepsfont_all 2d_*.ps && cp 2d_*.ps ~/tex/harald/toku2/figures && scpp 2d_*.ps e23:tex/harald/toku2/figures/ 
\center
\includegraphics[angle=-90,width=0.18\textwidth]{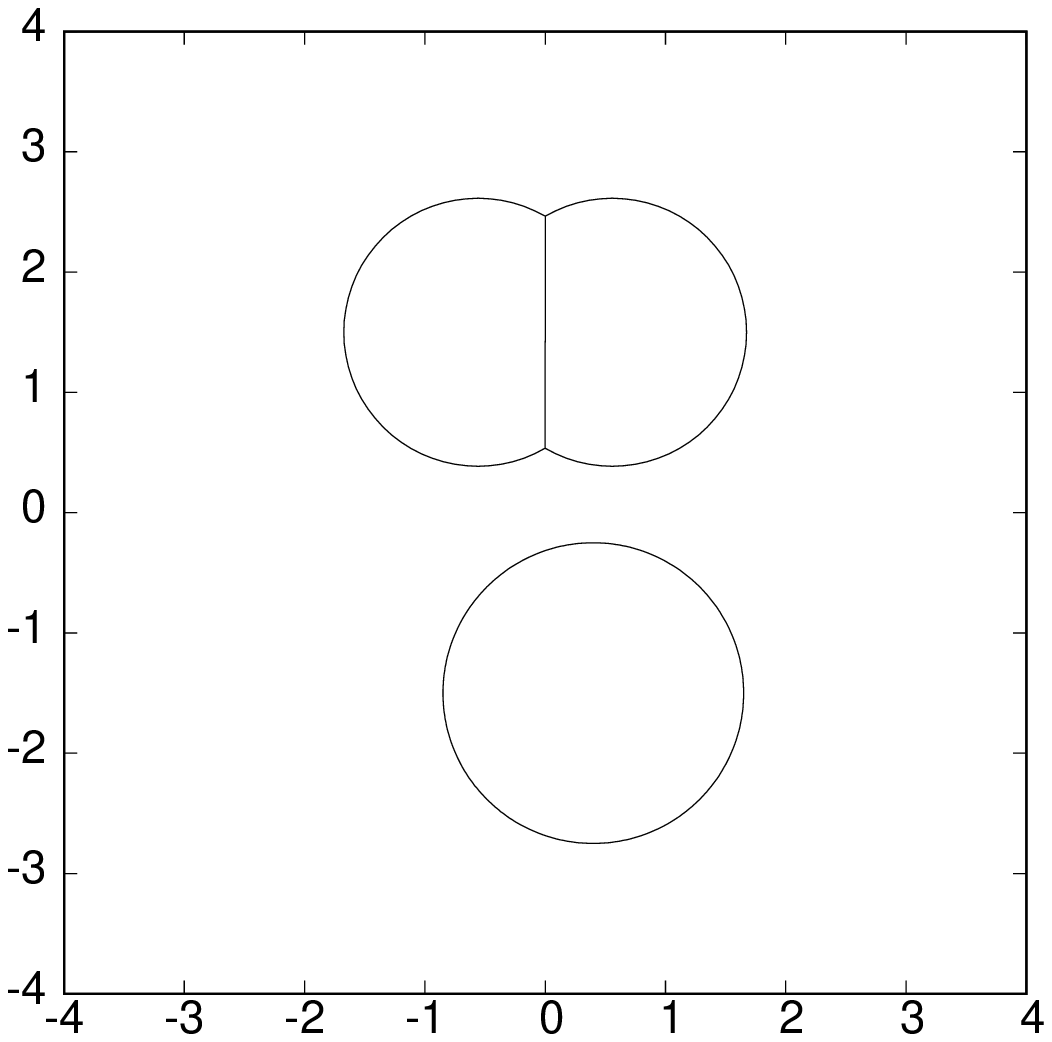}
\includegraphics[angle=-90,width=0.18\textwidth]{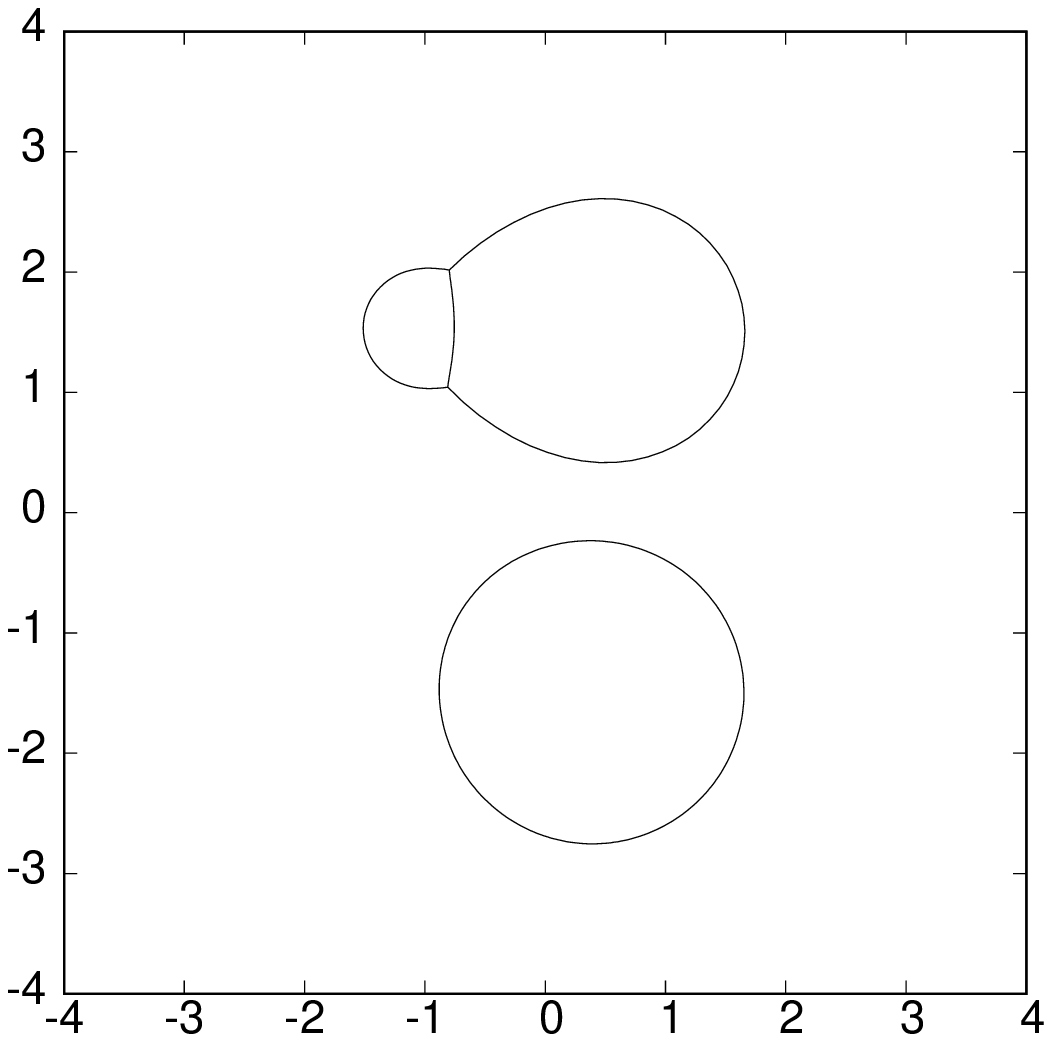}
\includegraphics[angle=-90,width=0.18\textwidth]{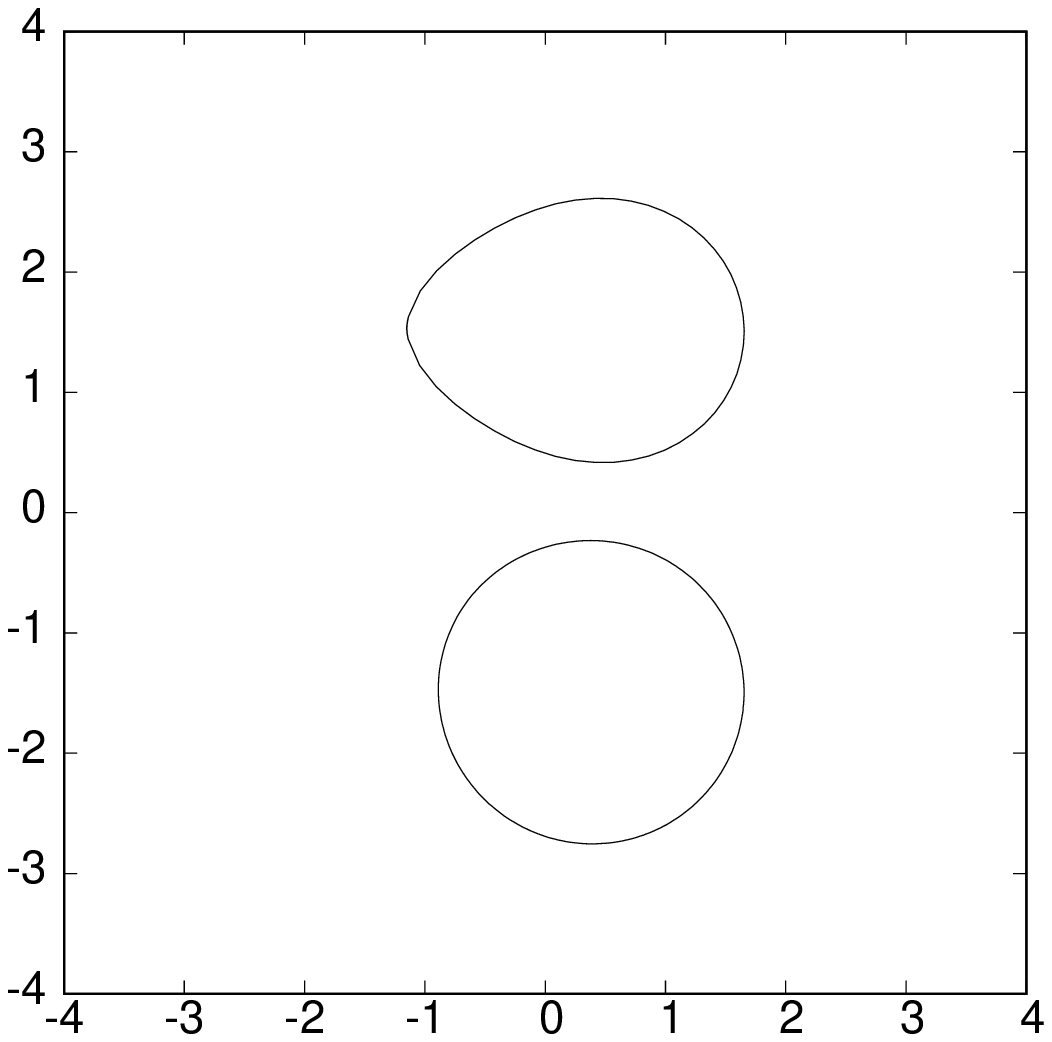}
\includegraphics[angle=-90,width=0.18\textwidth]{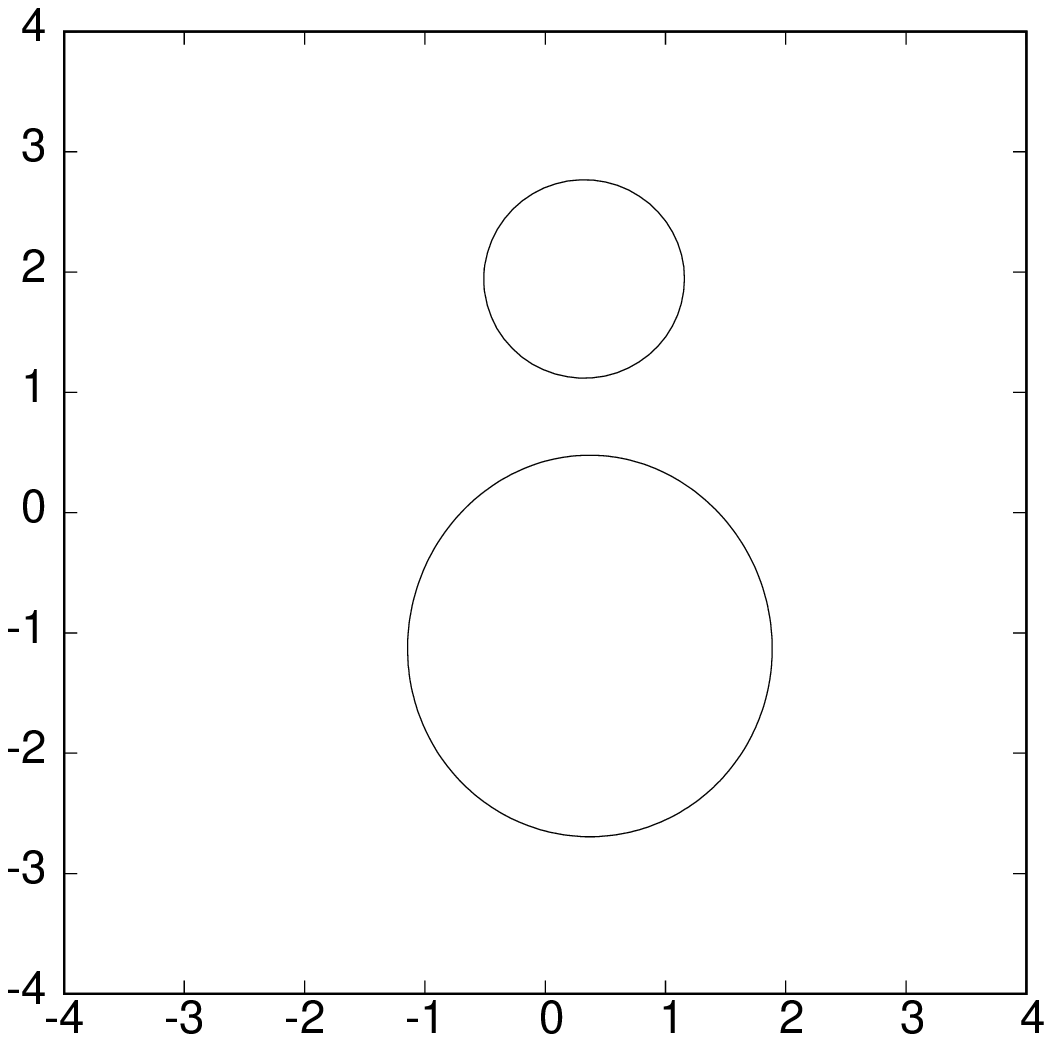}
\includegraphics[angle=-90,width=0.18\textwidth]{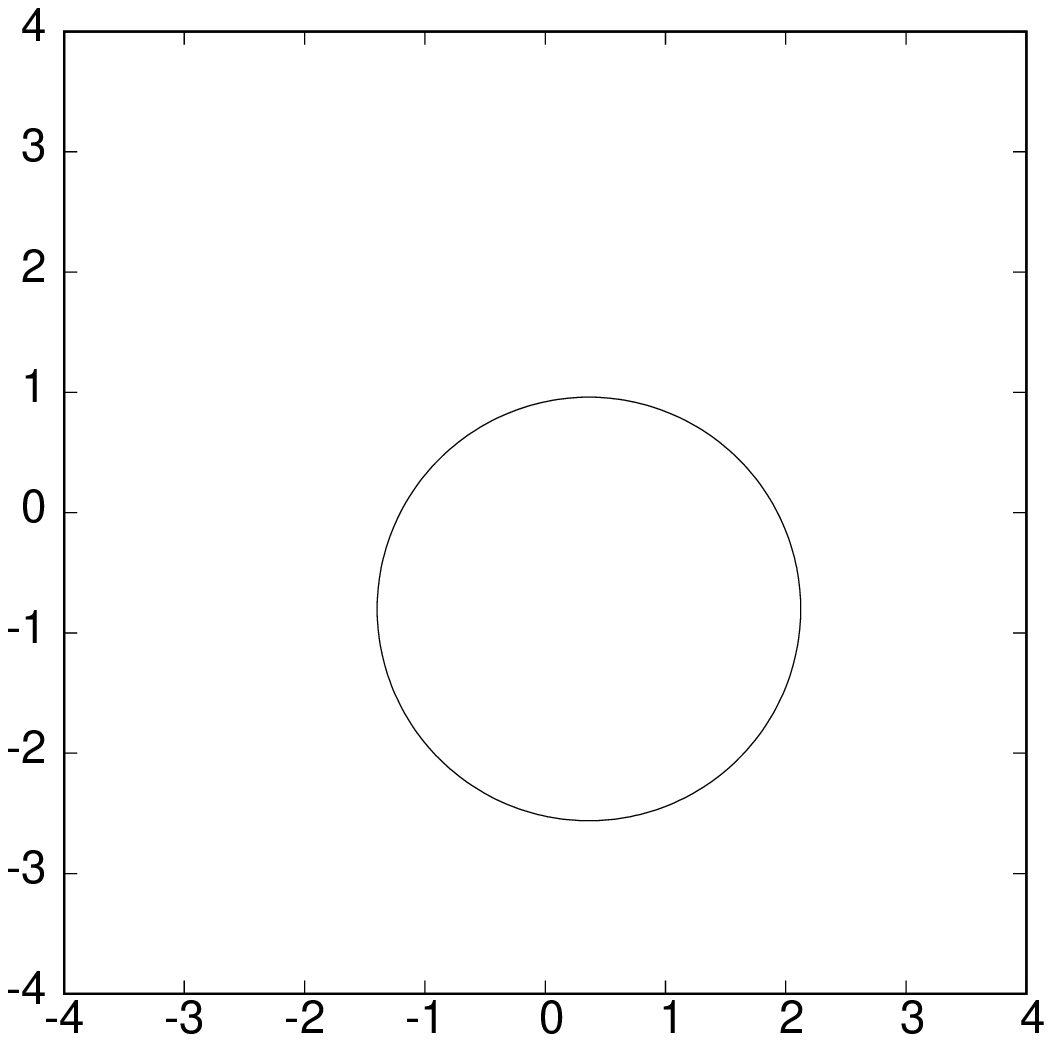}
\includegraphics[angle=-90,width=0.3\textwidth]{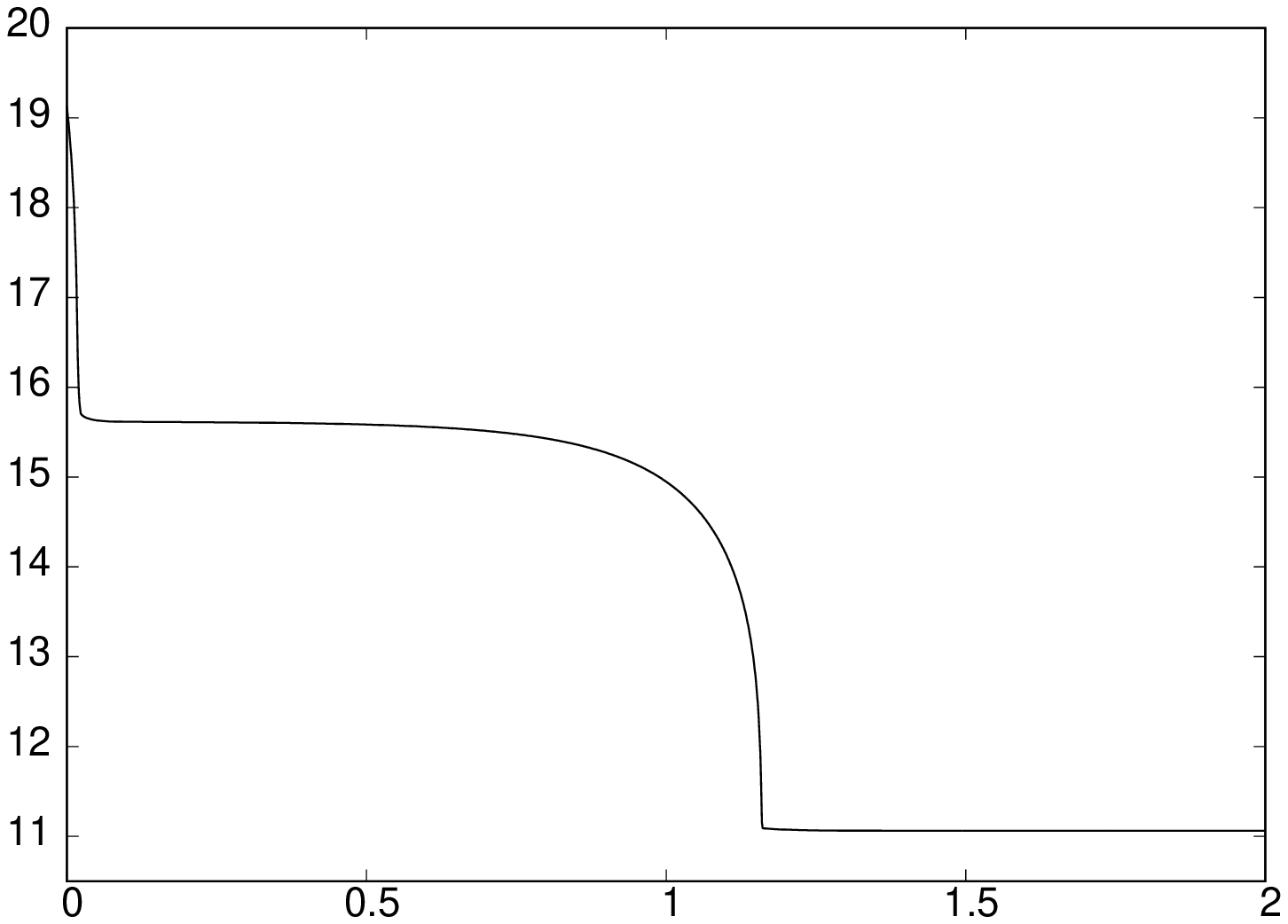} 
%\\[5mm]
%\includegraphics[angle=-0,width=0.18\textwidth]{figures/2d_db_plus_bigonemesh0}
%\includegraphics[angle=-0,width=0.18\textwidth]{figures/2d_db_plus_bigonemesh1}
\setlength{\abovecaptionskip}{20pt}
\caption{(($\beta_1, \beta_2,\beta_3) = (0.25,0,-0.25) =$ right bubble plus disk, left bubble, outer phase)\\
The solution at times $t=0, 0.015, 0.02, 1, 2$,
and a plot of the discrete energy over time.
For this computation, we have $\discreteVol{0} = -11.2$.
%Below we show the adaptive bulk mesh at times $t=0$ and $t=1$.
}
\label{fig:2d_db_plus_bigone}
% /home/rn/hpc_cluster/data/alberta/egn2/2d.db_plus_bigone
% /home/rn/hpc_cluster/data/alberta/egn2/2d.db_plus_bigone/vol3
\end{figure}

We repeat the same experiment with the disk which encloses a larger area, namely its radius is set as $\frac{5}{4}$ so that it encloses an area of $\frac{25\pi}{16} \approx 4.909$.
Then, as in the previous case, the left bubble shrinks and vanishes, and the right bubble and the disk survive.
However, to obey the energy minimization, the disk eventually absorbs the right bubble. See Figure~\ref{fig:2d_db_plus_bigone} for the results.

\subsection{Convergence experiment in 3d}

\begin{figure}[H]
% ../../plotradii
% cp radii_clean.ps ~/tex/harald/toku2/figures/3dsimul2c_r.ps && scpp radii_clean.ps e23:tex/harald/toku2/figures/3dsimul2c_r.ps
% cp radii_clean.ps ~/tex/harald/toku2/figures/3dsimul2c_long_r.ps && scpp radii_clean.ps e23:tex/harald/toku2/figures/3dsimul2c_long_r.ps
\centering
\includegraphics[angle=-90,keepaspectratio,scale=0.25]{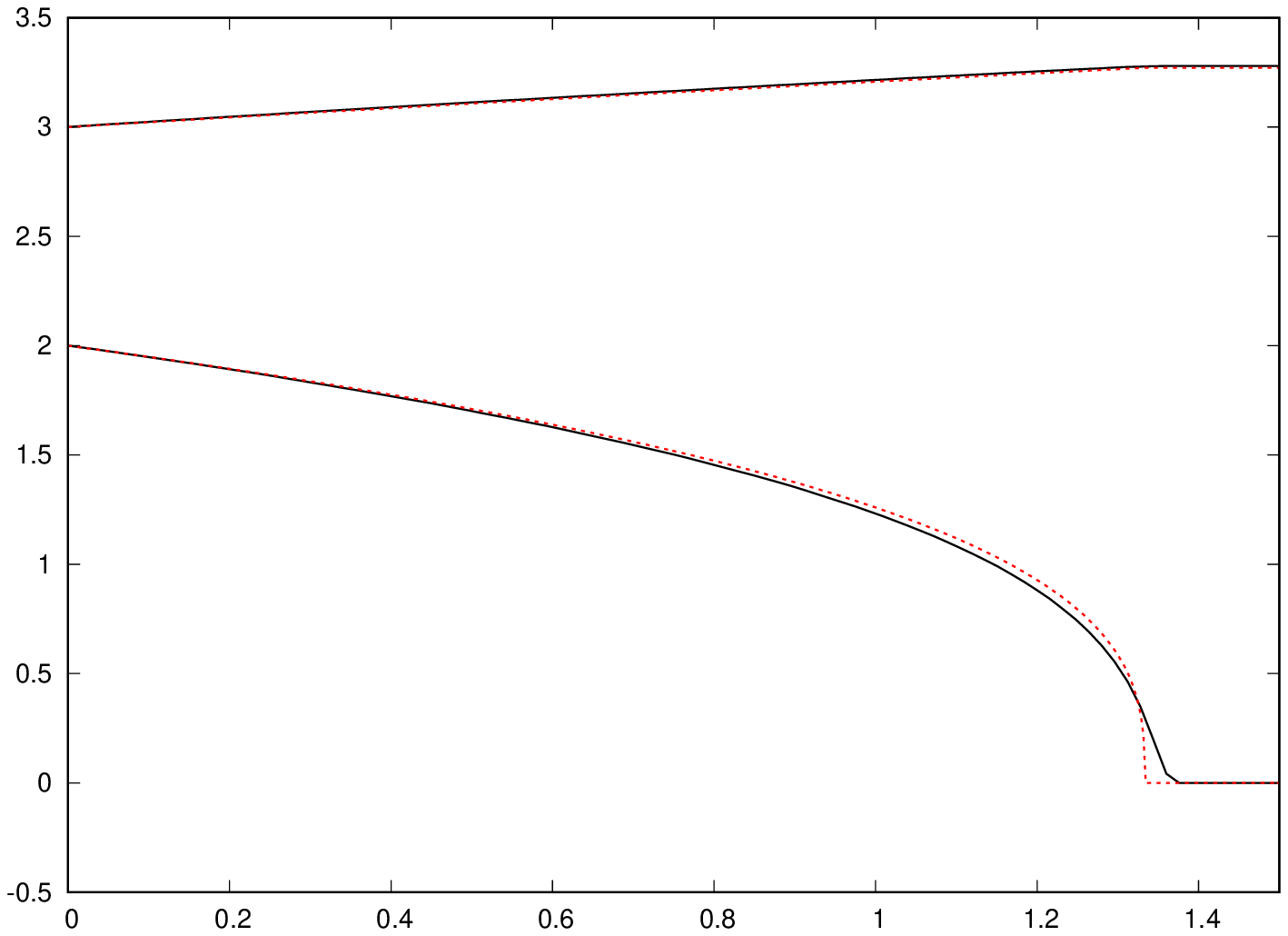}
\setlength{\abovecaptionskip}{20pt}
\caption{Comparison of the discrete (black, solid)\\
and exact (red, dashed) radii for $(\beta_1,\beta_2,\beta_3) = (-1,0,1)$.}
% ~/hpc_cluster/data/alberta/egn2/3d.convergence_experiment_NoC2/plot/run0
% ~/hpc_cluster/data/alberta/egn2/3d.convergence_experiment_NoC2/long/run1
\label{fig:3dsimul2c}
\end{figure}

We consider the 3d equivalent of the setup displayed in Figure~\ref{fig:2c}.
For the case $d = 3$, we have the rigorous solution to the three-phase system as follows:
% for general $d$.
% \begin{equation*}
%     w(x) =
%     \begin{cases}
%         -\frac{d-1}{\jump{2}{1}{\beta}R_1} & \qquad\mbox{if}\quad |x| < R_1,\\
%         -\frac{d-1}{\jump{2}{1}{\beta}R_1} - \frac{(d-1)\alpha(R_1,R_2)}{R_1^{d-2}} + \frac{(d-1)\alpha(R_1,R_2)}{|x|^{d-2}} & \qquad\mbox{if}\quad R_1\leq |x| < R_2,\\
%         \frac{d-1}{\jump{2}{3}{\beta}R_2} &\qquad\mbox{if}\quad R_2 \leq |x|,
%     \end{cases}
% \end{equation*}
\begin{equation*}
    w(x,t) =
    \begin{cases}
        -\frac{2}{\jump{2}{1}{\beta}R_1(t)} & \qquad\mbox{if}\quad |x| < R_1(t),\\
        -\frac{2}{\jump{2}{1}{\beta}R_1(t)} - \frac{2\alpha_3(R_1(t),R_2(t))}{R_1(t)} + \frac{2\alpha_3(R_1(t),R_2(t))}{|x|} & \qquad\mbox{if}\quad R_1(t)\leq |x| < R_2(t),\\
        \frac{2}{\jump{2}{3}{\beta}R_2(t)} &\qquad\mbox{if}\quad R_2(t) \leq |x|,
    \end{cases}
\end{equation*}
where
% for general $d$.
% \begin{equation*}
%     \alpha(R_1,R_2) := \frac{\frac{1}{\jump{2}{1}{\beta}R_1} + \frac{1}{\jump{2}{3}{\beta}R_2}}{\frac{1}{R_2^{d-2}} - \frac{1}{R_1^{d-2}}}.
% \end{equation*}
\begin{equation*}
    \alpha_3(R_1,R_2) := \frac{\frac{1}{\jump{2}{1}{\beta}R_1} + \frac{1}{\jump{2}{3}{\beta}R_2}}{\frac{1}{R_2} - \frac{1}{R_1}}\qquad \mbox{for}\quad R_1, R_2 > 0.
\end{equation*}
The formulae for $R_1(t)$ and $R_2(t)$ are given as the solutions to the following system of ordinary differential equations:
% for general $d$.
% \begin{equation*}
%     \begin{cases}
%         \dot{R}_1 = \frac{(d-1)(d-2)\alpha(R_1,R_2)}{\jump{2}{1}{\beta} R_1^{d-1}},\\
%         \dot{R}_2 = \frac{(d-1)(d-2)\alpha(R_1,R_2)}{\jump{2}{3}{\beta} R_2^{d-1}}.
%     \end{cases}
% \end{equation*}
\begin{align*}
    \begin{split}
        \dot{R}_1(t) &= \frac{2\alpha_3(R_1(t),R_2(t))}{\jump{2}{1}{\beta} R_1(t)^2},\\
        \dot{R}_2(t) &= \frac{2\alpha_3(R_1(t),R_2(t))}{\jump{2}{3}{\beta} R_2(t)^2}.
    \end{split}
\end{align*}
\Proposition{prop:mass_p} again implies that the function $t\mapsto\jump{2}{1}{\beta}R_1(t)^3 - \jump{2}{3}{\beta}R_2(t)^3$ is constant,
and hence we can solve the above system as in the case $d = 2$.
Namely, we have the ordinary differential equation for $R_2(t)$ as
% for general $d$.
% \begin{equation}\label{eq:ODE3d}
%     \dot{R}_2 = -F(R_2),\qquad F(u) := -\frac{(d-1)(d-2)\,\operatorname{\alpha}\left(\left(\frac{D+\jump{2}{3}{\beta}u^d}{\jump{2}{1}{\beta}}\right)^{\frac{1}{d}}, u\right)}{\jump{2}{3}{\beta}u^{d-1}}
% \end{equation}
\begin{equation}\label{eq:ODE3d}
    \dot{R}_2(t) = -F(R_2(t)),\qquad F(u) := -\frac{2\,\operatorname{\alpha_3}\left(\sqrt[3]{\frac{D_3+\jump{2}{3}{\beta}u^3}{\jump{2}{1}{\beta}}}, u\right)}{\jump{2}{3}{\beta}u^2}
\end{equation}
with $D_3 := \jump{2}{1}{\beta}R_1(0)^3 - \jump{2}{3}{\beta}R_2(0)^3$.
We show in Figure~\ref{fig:3dsimul2c} the result of accuracy check for the formulae for $R_1(t)$ and $R_2(t)$.

We also perform a convergence experiment for the true solution
of \eqref{eq:ODE3d} for $d=3$ over the time interval $[0,0.5]$.
For $i=0\to 3$, we set $N_f = 2^{5+i}$, $N_c = 4^{i}$ and $\frac12 K^0_\Gamma=\hat K(i)$,
where $(\hat K(0), \hat K(1), \hat K(2), \hat K(3)) = (770, 3074, 12290, 49154)$,
and $\tau= 4^{3-i}\times10^{-3}$. 
The errors $\errorUu$ and $\errorXx$ are displayed in 
Tables~\ref{tab:3dvol} and \ref{tab:3d}.
\begin{table}[H]
% grep "LaTeX:" run*/specs.used -A 2 | grep '\\\\' | sed 's/.*specs.used- //'
\center
\begin{tabular}{c|c|c|c|c|c|c}
 $h_{f}$ & $h^M_\Gamma$ & $\errorUu$ & $\errorXx$ & $K^M_\Omega$ & $K^M_\Gamma$ 
 & $|v_\Delta^M|$ \\ \hline
% i = 0, N_f = 32, N_c = 1, tau = 0.064
% i = 1, N_f = 64, N_c = 4, tau = 0.016
% i = 2, N_f = 128, N_c = 16, tau = 0.004
% i = 3, N_f = 256, N_c = 64 , tau = 0.001
2.5000e-01 & 6.4741e-01 & 5.1264e-02 & 4.3538e-02 & 11043 & 1540 &$<10^{-10}$\\
1.2500e-01 & 3.2398e-01 & 3.3996e-02 & 1.9571e-02 & 45469 & 6148 &$<10^{-10}$\\
6.2500e-02 & 1.6204e-01 & 1.5014e-02 & 9.6263e-03 & 188931& 24580&$<10^{-10}$\\
\end{tabular}
\caption{Convergence test for \eqref{eq:ODE3d} over the time interval $[0,0.5]$
for the scheme \eqref{eq:fma}.}
% ~/hpc_cluster/data/alberta/egn2/3d.convergence_experiment_NoC2/vol
\label{tab:3dvol}
\end{table}%
\begin{table}[H]
% grep "LaTeX:" run*/specs.used -A 2 | grep '\\\\' | sed 's/.*specs.used- //'
% grep -A1 "TEC" run*/specs.used 
\center
\begin{tabular}{c|c|c|c|c|c|c}
 $h_{f}$ & $h^M_\Gamma$ & $\errorUu$ & $\errorXx$ & $K^M_\Omega$ & $K^M_\Gamma$ 
 & $|v_\Delta^M|$ \\ \hline
% i = 0, N_f = 32, N_c = 1, tau = 0.064
% i = 1, N_f = 64, N_c = 4, tau = 0.016
% i = 2, N_f = 128, N_c = 16, tau = 0.004
% i = 3, N_f = 256, N_c = 64 , tau = 0.001
2.5000e-01 & 6.4636e-01 & 2.4897e-02 & 1.3525e-02 & 11043 & 1540 & 3.17e-01 \\
1.2500e-01 & 3.2393e-01 & 2.7050e-02 & 1.2229e-02 & 45661 & 6148 & 8.18e-02 \\
6.2500e-02 & 1.6207e-01 & 1.1793e-02 & 7.8416e-03 & 189027 & 24580& 2.02e-02\\
3.1250e-02 & 8.1017e-02 & 7.1962e-03 & 4.2435e-03 & 787657 & 98308& 4.95e-03\\
\end{tabular}
\caption{Convergence test for \eqref{eq:ODE3d} over the time interval $[0,0.5]$
for the scheme \eqref{eq:fma2}.}
% ~/hpc_cluster/data/alberta/egn2/3d.convergence_experiment_NoC2/lin
\label{tab:3d}
\end{table}%

\subsection{Numerical simulations in 3d}

% ~/hpc_cluster/data/alberta/egn2/3d.db_d_fine_tau4

\begin{figure}[H]
% rsync -avz --include="*.vtk" --include="energy.out" --include="uv*.out" --exclude="*" e23:`pwd`/ .
% ../plotenergy && mv energy.ps 3d_db_e.ps
% paraview --data=kappa_h0_..vtk
% fixbb 3d_*.ps && enlargepsfont_all 3d_*.ps && cp 3d_*.ps ~/tex/harald/toku2/figures && scpp 3d_*.ps e23:tex/harald/toku2/figures/ 
\centering
\includegraphics[angle=-0,width=0.3\textwidth]{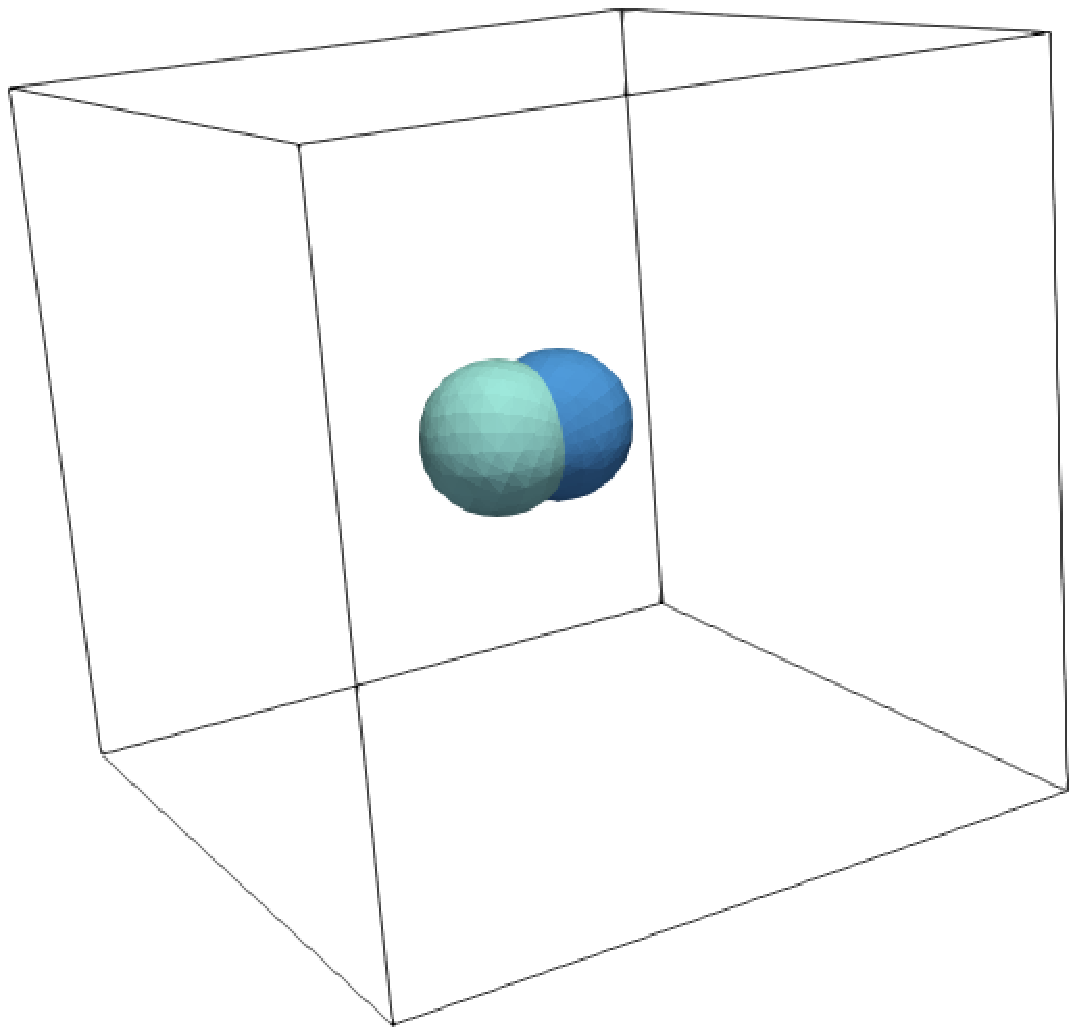}
\includegraphics[angle=-0,width=0.3\textwidth]{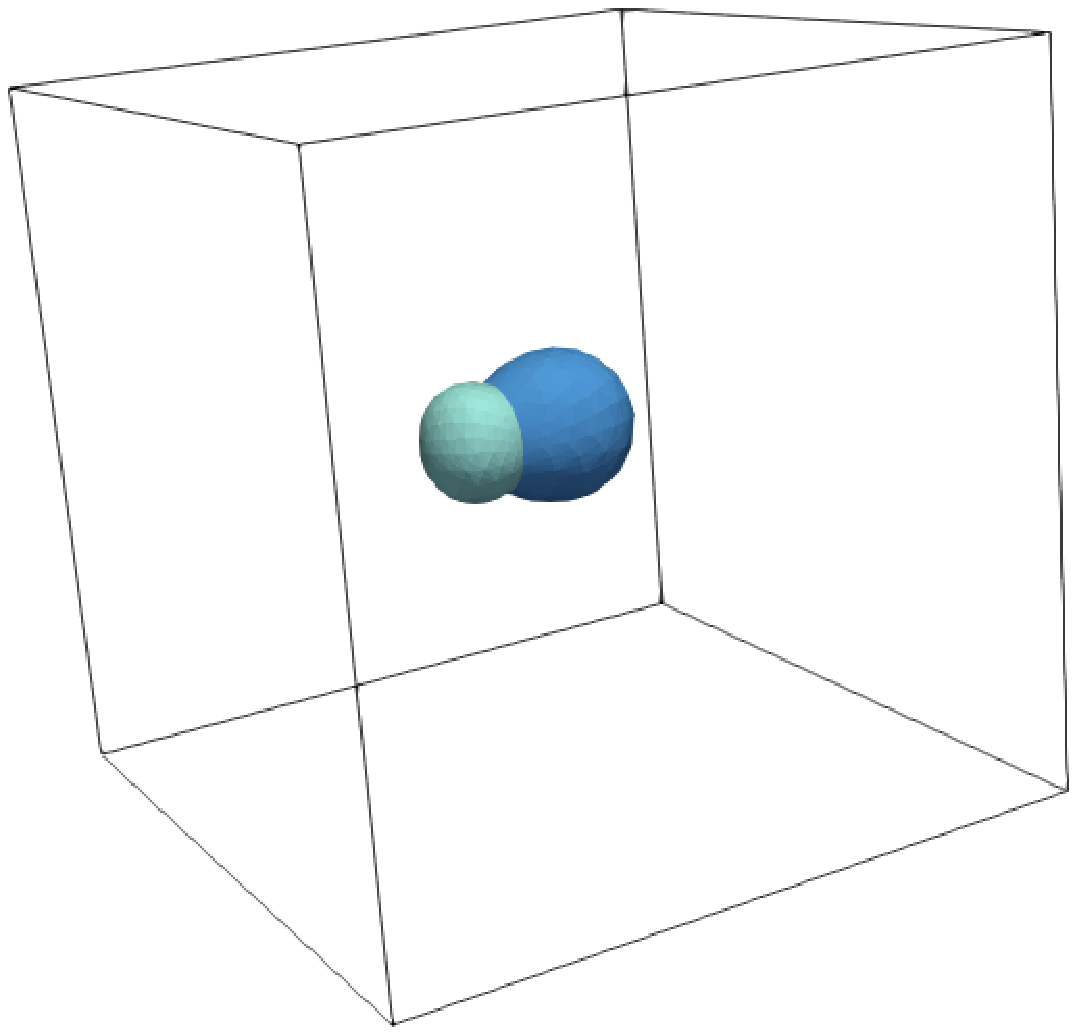}
\includegraphics[angle=-0,width=0.3\textwidth]{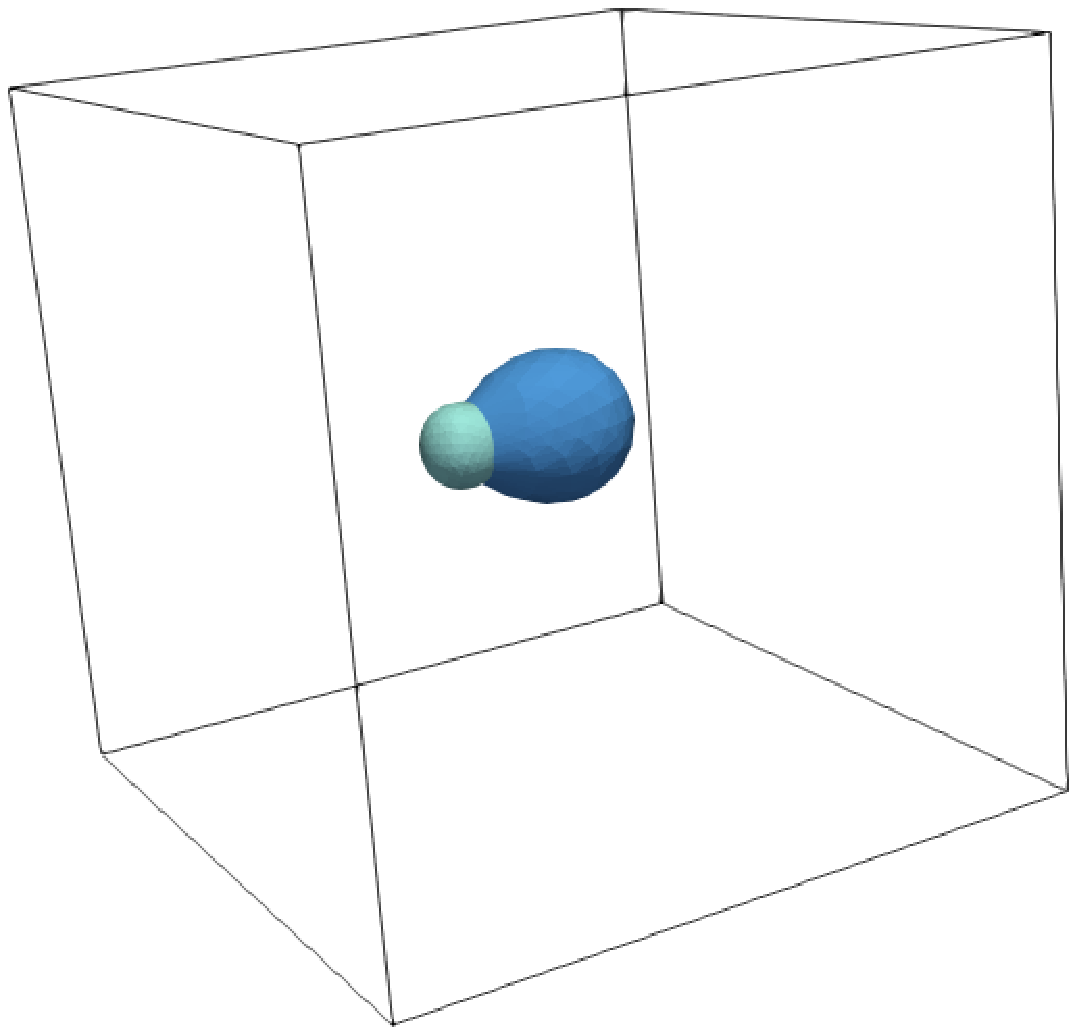}
\includegraphics[angle=-0,width=0.3\textwidth]{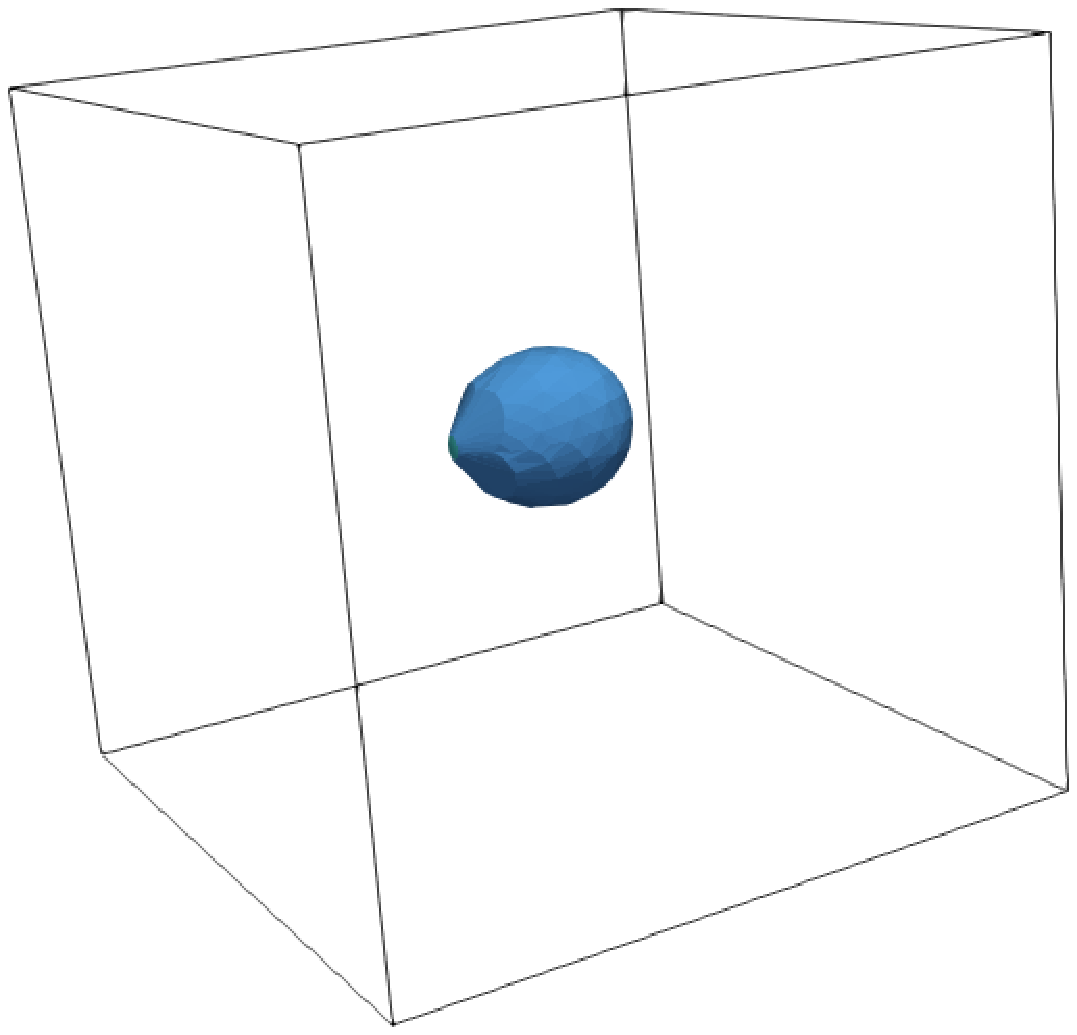}
\includegraphics[angle=-0,width=0.3\textwidth]{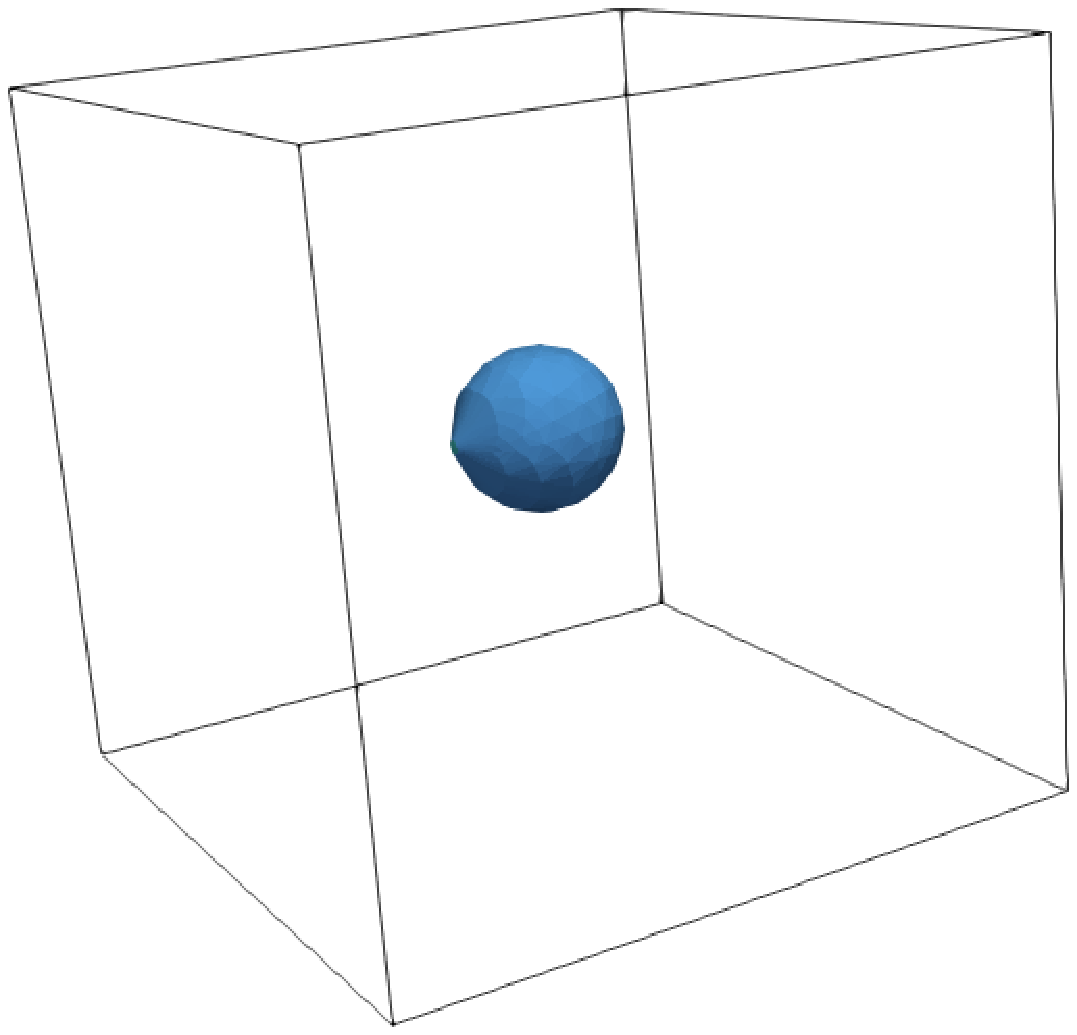} 
\includegraphics[angle=-0,width=0.3\textwidth]{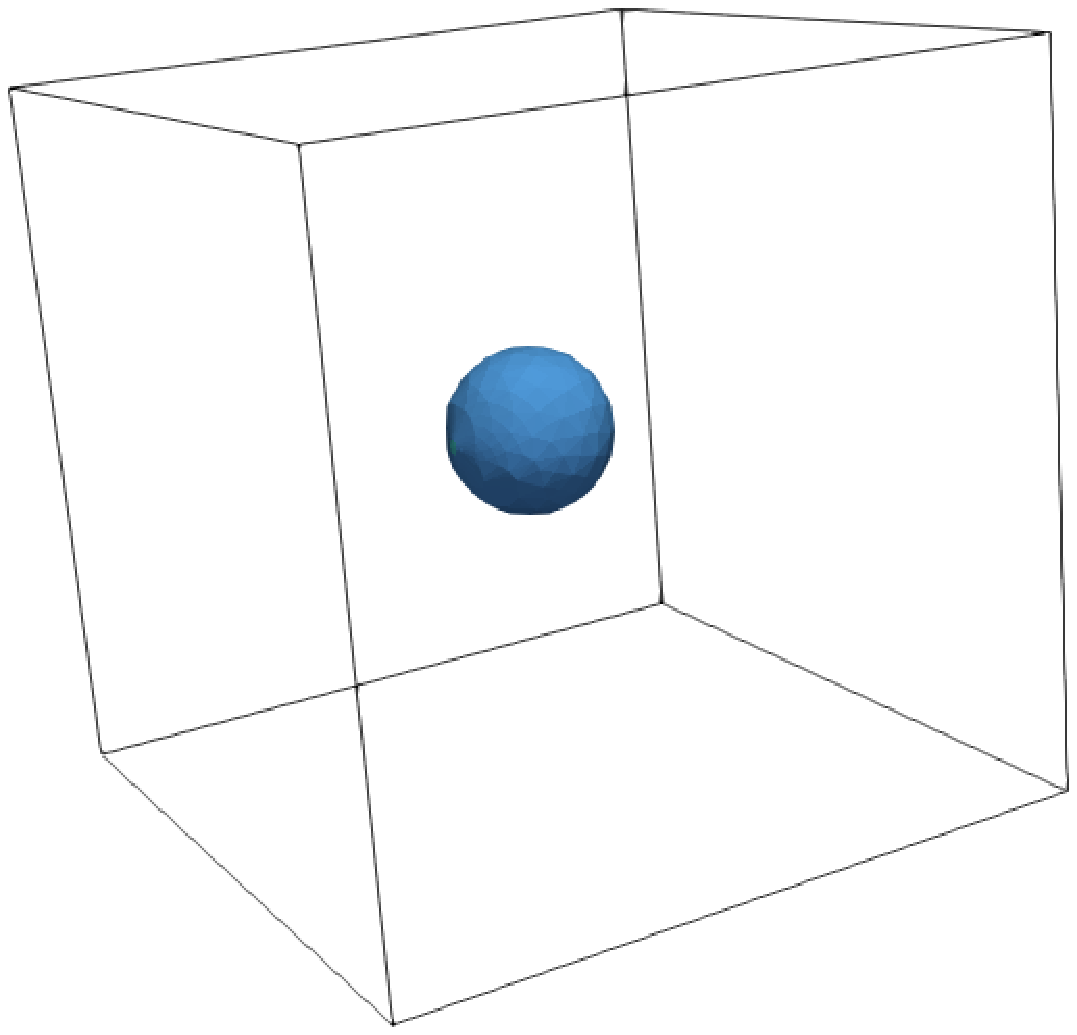} \\
\includegraphics[angle=-90,width=0.3\textwidth]{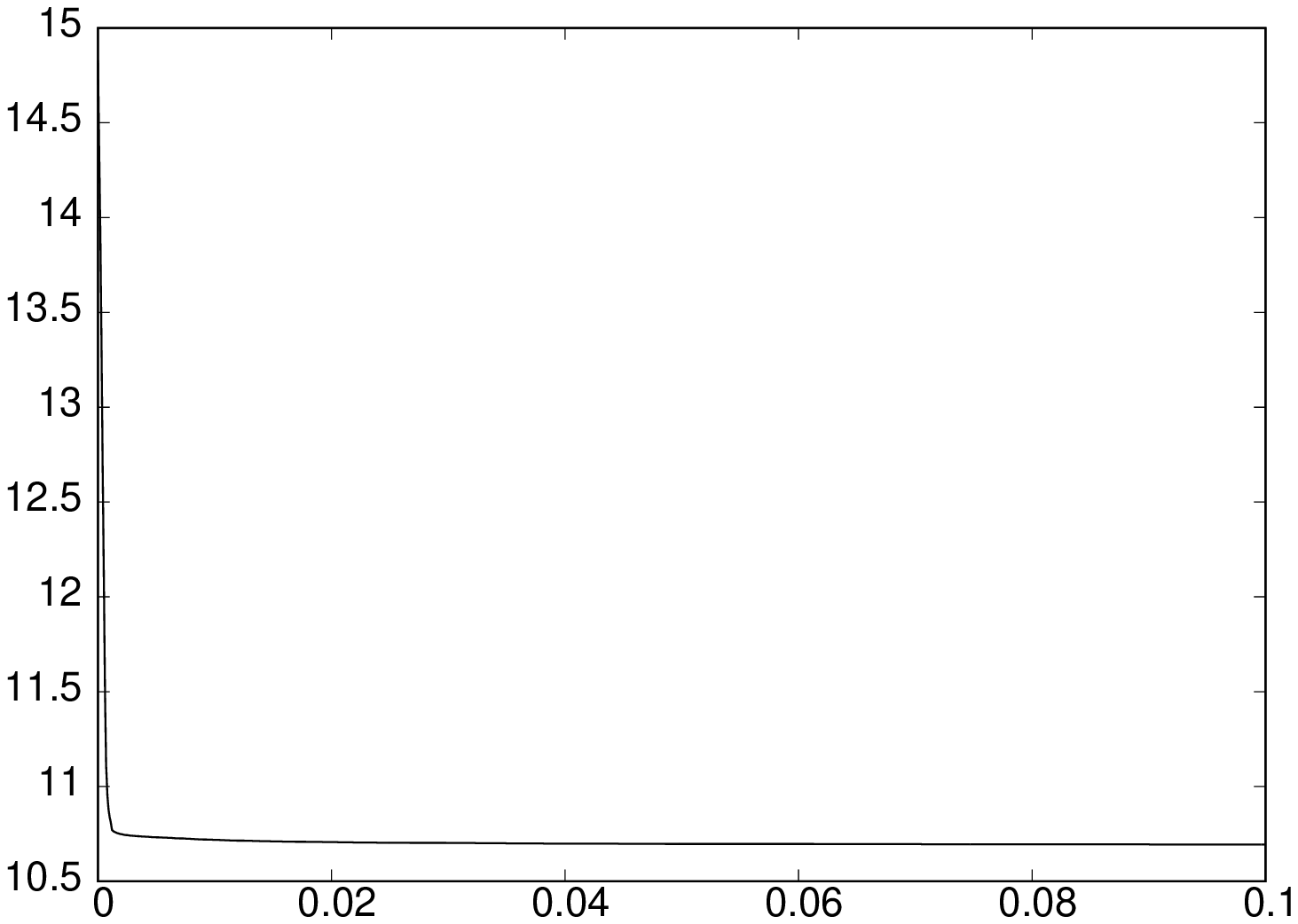} 
\setlength{\abovecaptionskip}{20pt}
\caption{(($\beta_1, \beta_2,\beta_3) = (-0.1,0,0.1) =$ right bubble, left bubble, outer phase)\\
The solution at times $t=0, 0.0002, 0.0003, 0.001, 0.01, 0.1$,
and a plot of the discrete energy over time. 
For this computation, we have $\discreteVol{0} = 50.6$.
}
\label{fig:3d_db}
% ~/hpc_cluster/data/alberta/egn2/3d.db_b_fine_tau4
% ~/hpc_cluster/data/alberta/egn2/3d.db_b_fine_short
\end{figure}

We consider a three-dimensional analogue of the simulation shown in 
Figure~\ref{fig:2d_db}, that is for a standard double bubble in 3d.
In this experiment, we set $I_R = 3$, $I_S = 3$, $I_T=1$, $(\beta_1, \beta_2,\beta_3) = (-0.1,0,0.1)$,
$(\cpindexp{1},\cpindexm{1}) = (3,1)$, $(\cpindexp{2}, \cpindexm{2}) = (2,3)$, and $(\cpindexp{3}, \cpindexm{3}) = (1,2)$.
See Figure~\ref{fig:3d_db} for the evolution, where we observe that the second
phases vanishes extremely quickly.

\section{Conclusion}\label{sec:con}
In this paper, we have proposed numerical methods for the degenerate multi-phase Stefan problem \eqref{eq:sharp_p} in two and three space dimensions.
The approach is based on the parametric finite element method and can deal
with several phases and the presence of triple junctions.
%The concerned problem is an extension of the two-phase Mullins--Sekerka problem to the multi-phase case,
%although the corresponding diffuse interface model is different from the one which was considered in the previous work \cite{EtoGarckeNurnberg2024}.
Overall we have introduced two schemes: a linear one and a nonlinear one.
Both schemes can be shown to be unconditionally stable, while the nonlinear one in addition also conserves the total energy content.
Moreover, we have carried out numerical experiments which show the practicality
and accuracy of the introduced methods.

\bibliographystyle{siam}
\bibliography{cite.bib}

\begin{thebibliography}{10}

\bibitem{AbelsMaxWilke}
{\sc H.~Abels, M.~Rauchecker, and M.~Wilke}, {\em Well-posedness and
  qualitative behaviour of the {M}ullins-{S}ekerka problem with ninety-degree
  angle boundary contact}, Math. Ann., 381 (2021), pp.~363--403.

\bibitem{BanschDeckelnickGarckePozzi2023}
{\sc E.~B\"{a}nsch, K.~Deckelnick, H.~Garcke, and P.~Pozzi}, {\em Interfaces:
  {M}odeling, {A}nalysis, {N}umerics}, Oberwolfach Seminars, Birkh\"{a}user
  Cham, 2023.

\bibitem{BaoGarckeNuernbergZhao2022}
{\sc W.~Bao, H.~Garcke, R.~N\"{u}rnberg, and Q.~Zhao}, {\em A
  structure-preserving finite element approximation of surface diffusion for
  curve networks and surface clusters}, Numer. Methods Partial Differ. Eq., 39
  (2023), pp.~759--794.

\bibitem{BaoLi2023}
{\sc W.~Bao and Y.~Li}, {\em A symmetrized parametric finite element method for
  anisotropic surface diffusion in three dimensions}, SIAM J. Sci. Comput., 45
  (2023), pp.~A1438--A1461.

\bibitem{BaoLi2024}
\leavevmode\vrule height 2pt depth -1.6pt width 23pt, {\em A
  structure-preserving parametric finite element method for geometric flows
  with anisotropic surface energy}, Numer. Math., 156 (2024), pp.~609--639.

\bibitem{BaoZhao2021}
{\sc W.~Bao and Q.~Zhao}, {\em A structure-preserving parametric finite element
  method for surface diffusion}, SIAM J. Numer. Anal., 59 (2021),
  pp.~2775--2799.

\bibitem{BarrettBlowey1995}
{\sc J.~W. Barrett and J.~F. Blowey}, {\em An error bound for the finite
  element approximation of the {C}ahn-{H}illiard equation with logarithmic free
  energy}, Numer. Math., 72 (1995), pp.~1--20.

\bibitem{BarrettBloweyGarcke1999}
{\sc J.~W. Barrett, J.~F. Blowey, and H.~Garcke}, {\em Finite element
  approximation of the {C}ahn-{H}illiard equation with degenerate mobility},
  SIAM J. Numer. Anal., 37 (1999), pp.~286--318.

\bibitem{BarrettBloweyGarcke2001}
\leavevmode\vrule height 2pt depth -1.6pt width 23pt, {\em On fully practical
  finite element approximations of degenerate {C}ahn-{H}illiard systems}, M2AN
  Math. Model. Numer. Anal., 35 (2001), pp.~713--748.

\bibitem{BarrettGarckeRobert2007}
{\sc J.~W. Barrett, H.~Garcke, and R.~N{\"u}rnberg}, {\em {O}n the variational
  approximation of combined second and fourth order geometric evolution
  equations}, SIAM J. Sci. Comput., 29 (2007), pp.~1006--1041.

\bibitem{BarrettGarckeNuernbergSurfaceDiffusion2007}
{\sc J.~W. Barrett, H.~Garcke, and R.~N{\"u}rnberg}, {\em A parametric finite
  element method for fourth order geometric evolution equations}, J. Comput.
  Phys., 222 (2007), pp.~441--462.

\bibitem{ejam3d}
{\sc J.~W. Barrett, H.~Garcke, and R.~N\"urnberg}, {\em Finite element
  approximation of coupled surface and grain boundary motion with applications
  to thermal grooving and sintering}, European J. Appl. Math., 21 (2010),
  pp.~519--556.

\bibitem{BarrettGarckeRobert2010}
{\sc J.~W. Barrett, H.~Garcke, and R.~N{\"u}rnberg}, {\em On stable parametric
  finite element methods for the {S}tefan problem and the {M}ullins--{S}ekerka
  problem with applications to dendritic growth}, J. Comput. Phys., 229 (2010),
  pp.~6270--6299.

\bibitem{clust3d}
{\sc J.~W. Barrett, H.~Garcke, and R.~N\"urnberg}, {\em Parametric
  approximation of surface clusters driven by isotropic and anisotropic surface
  energies}, Interfaces Free Bound., 12 (2010), pp.~187--234.

\bibitem{BarrettGarckeRobertBook2020}
{\sc J.~W. Barrett, H.~Garcke, and R.~N{\"u}rnberg}, {\em Parametric finite
  element approximations of curvature-driven interface evolutions}, Handb.
  Numer. Anal., 21 (2020), pp.~275--423.

\bibitem{Bates1995ANS}
{\sc P.~Bates, X.~Chen, and X.~Deng}, {\em A numerical scheme for the two phase
  {M}ullins-{S}ekerka problem}, Electronic Journal of Differential Equations,
  1995 (1995), pp.~1--27.

\bibitem{BatesBrown}
{\sc P.~W. Bates and S.~Brown}, {\em A numerical scheme for the
  {M}ullins-{S}ekerka evolution in three space dimensions}, in Differential
  equations and computational simulations ({C}hengdu, 1999), World Sci. Publ.,
  River Edge, NJ, 2000, pp.~12--26.

\bibitem{BloweyCopettiElliott1996}
{\sc J.~F. Blowey, M.~I.~M. Copetti, and C.~M. Elliott}, {\em Numerical
  analysis of a model for phase separation of a multi-component alloy}, IMA J.
  Numer. Anal., 16 (1996), pp.~111--139.

\bibitem{BloweyElliott1992}
{\sc J.~F. Blowey and C.~M. Elliott}, {\em The {C}ahn-{H}illiard gradient
  theory for phase separation with nonsmooth free energy. {II}. {N}umerical
  analysis}, European J. Appl. Math., 3 (1992), pp.~147--179.

\bibitem{BretinRolandMasnouSengersTerii2023}
{\sc E.~Bretin, R.~Denis, S.~Masnou, A.~Sengers, and G.~Terii}, {\em A
  multiphase {C}ahn-{H}illiard system with mobilities and the numerical
  simulation of dewetting}, ESAIM Math. Model. Numer. Anal., 57 (2023),
  pp.~1473--1509.

\bibitem{BronsardGarckeStoth1998}
{\sc L.~Bronsard, H.~Garcke, and B.~Stoth}, {\em A multi-phase
  {M}ullins-{S}ekerka system: matched asymptotic expansions and an implicit
  time discretisation for the geometric evolution problem}, Proc. Roy. Soc.
  Edinburgh Sect. A, 128 (1998), pp.~481--506.

\bibitem{ChambolleLaux2021}
{\sc A.~Chambolle and T.~Laux}, {\em Mullins-{S}ekerka as the {W}asserstein
  flow of the perimeter}, Proc. Amer. Math. Soc., 149 (2021), pp.~2943--2956.

\bibitem{ChenKublikTsai2017}
{\sc C.~Chen, C.~Kublik, and R.~Tsai}, {\em An implicit boundary integral
  method for interfaces evolving by {M}ullins-{S}ekerka dynamics}, in
  Mathematics for nonlinear phenomena---analysis and computation, vol.~215 of
  Springer Proc. Math. Stat., Springer, Cham, 2017, pp.~1--21.

\bibitem{ChenMettimanOsherSmereka1997}
{\sc S.~Chen, B.~Merriman, S.~Osher, and P.~Smereka}, {\em A simple level set
  method for solving {S}tefan problems}, J. Comput. Phys., 135 (1997),
  pp.~8--29.

\bibitem{ChenHongYi1996}
{\sc X.~Chen, J.~Hong, and F.~Yi}, {\em Existence, uniqueness, and regularity
  of classical solutions of the {M}ullins-{S}ekerka problem}, Comm. Partial
  Differential Equations, 21 (1996), pp.~1705--1727.

\bibitem{Davis04}
{\sc T.~A. Davis}, {\em Algorithm 832: {UMFPACK} {V}4.3---an
  unsymmetric-pattern multifrontal method}, ACM Trans. Math. Software, 30
  (2004), pp.~196--199.

\bibitem{Davis11}
{\sc T.~A. Davis}, {\em Algorithm 915, {SuiteSparseQR}: Multifrontal
  multithreaded rank-revealing sparse {QR} factorization}, ACM Trans. Math.
  Software, 38 (2011), pp.~1--22.

\bibitem{ElliottFrench1987}
{\sc C.~M. Elliott and D.~A. French}, {\em Numerical studies of the
  {C}ahn-{H}illiard equation for phase separation}, IMA J. Appl. Math., 38
  (1987), pp.~97--128.

\bibitem{EscherMatiocMatioc2024}
{\sc J.~Escher, A.-V. Matioc, and B.-V. Matioc}, {\em The {M}ullins-{S}ekerka
  problem via the method of potentials}, Math. Nachr., 297 (2024),
  pp.~1960--1977.

\bibitem{EScherSimonett}
{\sc J.~Escher and G.~Simonett}, {\em Classical solutions for {H}ele-{S}haw
  models with surface tension}, Adv. Differential Equations, 2 (1997),
  pp.~619--642.

\bibitem{EScherSimonett1998}
{\sc J.~Escher and G.~Simonett}, {\em A center manifold analysis for the
  {M}ullins-{S}ekerka model}, J. Differential Equations, 143 (1998),
  pp.~267--292.

\bibitem{Eto2023}
{\sc T.~Eto}, {\em A rapid numerical method for the {M}ullins-{S}ekerka flow
  with application to contact angle problems}, J. Sci. Comput., 98 (2024),
  p.~63.

\bibitem{EtoGarckeNurnberg2024}
{\sc T.~Eto, H.~Garcke, and R.~N\"urnberg}, {\em A structure-preserving finite
  element method for the multi-phase {M}ullins-{S}ekerka problem with triple
  junctions}, Numer. Math., 156 (2024), pp.~1479--1509.

\bibitem{Eyre1993}
{\sc D.~J. Eyre}, {\em Systems of {C}ahn-{H}illiard equations}, SIAM J. Appl.
  Math., 53 (1993), pp.~1686--1712.

\bibitem{FengProhl2004}
{\sc X.~Feng and A.~Prohl}, {\em Error analysis of a mixed finite element
  method for the {C}ahn-{H}illiard equation}, Numer. Math., 99 (2004),
  pp.~47--84.

\bibitem{FengProhl2005}
{\sc X.~Feng and A.~Prohl}, {\em Numerical analysis of the {C}ahn-{H}illiard
  equation and approximation of the {H}ele-{S}haw problem}, Interfaces Free
  Bound., 7 (2005), pp.~1--28.

\bibitem{FischerHenselLauxSimon2024}
{\sc J.~Fischer, S.~Hensel, T.~Laux, and T.~M. Simon}, {\em A weak-strong
  uniqueness principle for the {M}ullins-{S}ekerka equation}, arXiv preprint
  arXiv:2404.02682,  (2024).

\bibitem{Garcke2013}
{\sc H.~Garcke}, {\em Curvature driven interface evolution}, Jahresber. Dtsch.
  Math.-Ver., 115 (2013), pp.~63--100.

\bibitem{GarckeNuenbergZhao2023}
{\sc H.~Garcke, R.~N\"urnberg, and Q.~Zhao}, {\em Structure-preserving
  discretizations of two-phase {N}avier-{S}tokes flow using fitted and unfitted
  approaches}, J. Comput. Phys., 489 (2023), p.~112276.

\bibitem{GarckeNuernbergZhao2024}
{\sc H.~Garcke, R.~N\"{u}rnberg, and Q.~Zhao}, {\em A variational
  front-tracking method for multiphase flow with triple junctions}, Math.
  Comp.,  (2025), pp.~1--36.

\bibitem{GarckeRauchecker2022}
{\sc H.~Garcke and M.~Rauchecker}, {\em Stability analysis for stationary
  solutions of the {M}ullins-{S}ekerka flow with boundary contact}, Math.
  Nachr., 295 (2022), pp.~683--705.

\bibitem{GarckeSturzenhecker1998}
{\sc H.~Garcke and T.~Sturzenhecker}, {\em The degenerate multi-phase {S}tefan
  problem with {G}ibbs-{T}homson law}, Adv. Math. Sci. Appl., 8 (1998),
  pp.~929--941.

\bibitem{HenselStinson2024}
{\sc S.~Hensel and K.~Stinson}, {\em Weak solutions of {M}ullins-{S}ekerka flow
  as a {H}ilbert space gradient flow}, Arch. Ration. Mech. Anal., 248 (2024),
  p.~8.

\bibitem{IzsakDjebbar2023}
{\sc F.~Izsák and T.-E. Djebbar}, {\em Learning {D}ata for
  {N}eural-{N}etwork-{B}ased {N}umerical {S}olution of {PDE}s: {A}pplication to
  {D}irichlet-to-{N}eumann {P}roblems}, Algorithms, 16 (2023), p.~111.

\bibitem{LiZhao2024}
{\sc M.~Li and Q.~Zhao}, {\em Parametric finite element approximations for
  anisotropic surface diffusion with axisymmetric geometry}, J. Comput. Phys.,
  497 (2024), p.~112632.

\bibitem{LiChoiKim2016}
{\sc Y.~Li, J.~Choi, and J.~Kim}, {\em Multi-component {C}ahn-{H}illiard system
  with different boundary conditions in complex domains}, J. Comput. Phys., 323
  (2016), pp.~1--16.

\bibitem{LiLiuXiaHeLi2022}
{\sc Y.~Li, R.~Liu, Q.~Xia, C.~He, and Z.~Li}, {\em First- and second-order
  unconditionally stable direct discretization methods for multi-component
  {C}ahn-{H}illiard system on surfaces}, J. Comput. Appl. Math., 401 (2022),
  p.~113778.

\bibitem{LuckhausSturzenhecker}
{\sc S.~Luckhaus and T.~Sturzenhecker}, {\em Implicit time discretization for
  the mean curvature flow equation}, Calc. Var. Partial Differential Equations,
  3 (1995), pp.~253--271.

\bibitem{Mayer2000}
{\sc U.~F. Mayer}, {\em A numerical scheme for moving boundary problems that
  are gradient flows for the area functional}, European J. Appl. Math., 11
  (2000), pp.~61--80.

\bibitem{Nurnberg09}
{\sc R.~N\"urnberg}, {\em Numerical simulations of immiscible fluid clusters},
  Appl. Numer. Math., 59 (2009), pp.~1612--1628.

\bibitem{Nurnberg202203}
\leavevmode\vrule height 2pt depth -1.6pt width 23pt, {\em A structure
  preserving front tracking finite element method for the {M}ullins--{S}ekerka
  problem}, J. Numer. Math., 31 (2023), pp.~137--155.

\bibitem{Roger2005}
{\sc M.~R{\"o}ger}, {\em Existence of weak solutions for the
  {M}ullins-{S}ekerka flow}, SIAM J. Math. Anal., 37 (2005), pp.~291--301.

\bibitem{Alberta}
{\sc A.~Schmidt and K.~G. Siebert}, {\em Design of Adaptive Finite Element
  Software: The Finite Element Toolbox {ALBERTA}}, vol.~42 of Lecture Notes in
  Computational Science and Engineering, Springer-Verlag, Berlin, 2005.

\bibitem{Serfaty2011}
{\sc S.~Serfaty}, {\em Gamma-convergence of gradient flows on {H}ilbert and
  metric spaces and applications}, Discrete Contin. Dyn. Syst., 31 (2011),
  pp.~1427--1451.

\bibitem{ZhuChenHou}
{\sc J.~Zhu, X.~Chen, and T.~Hou}, {\em An {E}fficient {B}oundary {I}ntegral
  {M}ethod for the {M}ullins–{S}ekerka {P}roblem.}, J. Comput. Phys., 127
  (1996), pp.~246--267.

\end{thebibliography}
\end{document}